\documentclass[10pt,a4paper]{article}
\usepackage[utf8]{inputenc}
\usepackage[top = 1.4in, bottom = 1.4in, left = 1.2in, right = 1.2in]{geometry}

\usepackage[font=small, labelfont=bf]{caption}
\usepackage{graphicx}
\usepackage{amsmath}
\usepackage{amssymb}
\usepackage{lscape}
\usepackage{authblk}
\usepackage[colorlinks=true,urlcolor=blue]{hyperref} 
\usepackage{lineno}

\usepackage{float}            
\usepackage[section]{placeins}
\usepackage{microtype}        
\usepackage{epstopdf}         
\graphicspath{{images/}{Figs/}} 

\title{Memory reshapes stability landscapes: resilience–resistance tradeoffs and critical transitions}

\author[1,*]{Moein Khalighi}
\author[1]{Chandler Ross} 
\author[1]{Ville Laitinen}
\author[1,2]{Guilhem Sommeria-Klein}
\author[1]{Leo Lahti}

\affil[1]{Department of Computing, University of Turku; 20014 Turku, Finland}
\affil[2]{Inria, Univ. Bordeaux, INRAE; F-33400 Talence, France}
\affil[*]{Corresponding author: moein.khalighi@utu.fi}

\date{}

\begin{document}

\maketitle
\section*{Abstract}
Regime shifts in biology, ecology, and other complex systems are often interpreted through stability landscapes and early warning signals that implicitly assume dynamics without memory effects. Yet many real systems exhibit these effects, thus present dynamics depend on past states and past forcing. Here, we study how memory reshapes bistable stability landscapes and regime shifts using a minimal bistable model with a fractional derivative that controls memory strength. We connect landscape geometry to classical notions of resilience and resistance by quantifying basin curvature and the perturbation magnitude required to cross the unstable threshold, and we track how these quantities evolve in time after perturbations. Memory typically flattens basin floors, slowing recovery, while often increasing the perturbation threshold for stability transitions, revealing a tradeoff between resilience and resistance. Because the landscape becomes history-dependent and time-evolving, memory generates qualitative behaviors that do not appear in memory-free models, including delayed collapse or recovery after stress ends, rebound after apparently successful transition, and broadened hysteresis under gradual parameter change. Finally, we show that fitting a memory-free model to memory-driven data can reproduce trajectories while systematically shifting equilibrium branches and tipping locations, risking incorrect diagnosis and management of regime shifts. These results motivate a moving landscape view and provide practical guidance for interpreting observed anomalies and distinguishing memory-driven effects from noise.

\section*{Introduction}
\begin{figure}[ht!]
    \centering
    \includegraphics[width=\linewidth]{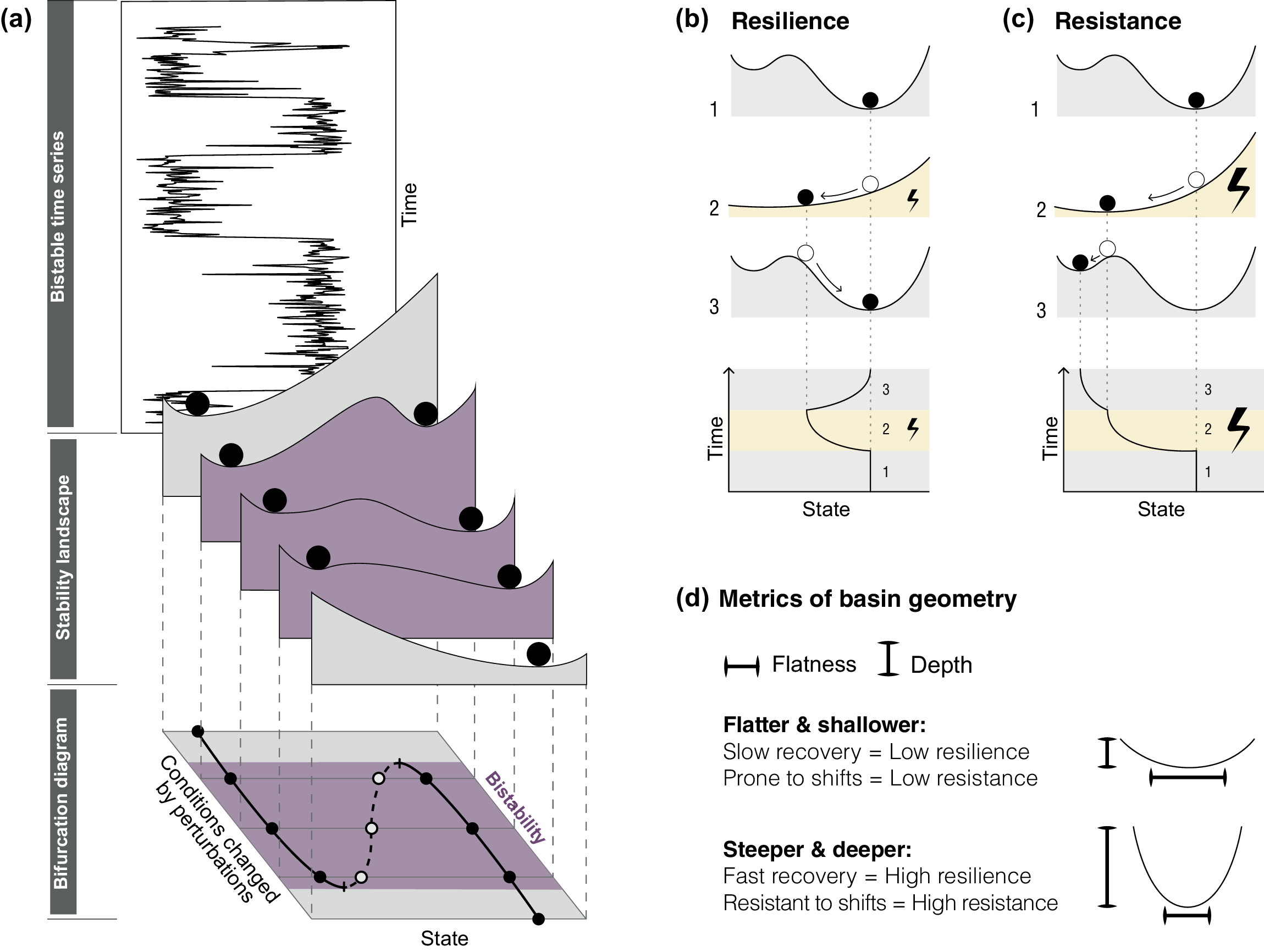}
    \caption{(a) The top panel shows a bistable time series of a univariate system. The middle panel illustrates its stability landscapes at different times as conditions change. The bottom panel projects the stability landscape onto a bifurcation diagram. In the stability landscape, valleys represent stable equilibria and peaks represent unstable states. The bifurcation diagram highlights equilibria between two tipping points, with filled circles indicating stable states and unfilled circles indicating unstable ones. Changes in system conditions can drive transitions between stable states. Here, the conditions are altered by endogenous perturbations that are tuned parameter values. The purple region denotes the parameter range in which bistability occurs.
    (b) This panel represents resilience, interpreted as the system’s rate of recovery (Fig.~\ref{fig: RL}). When the system is in a stable state (the ball is at the bottom of a basin of attraction), a mild or brief perturbation may deform the stability landscape and momentarily shift the system, but the state returns to its original equilibrium after the disturbance ends. 
    (c) This panel illustrates resistance, defined as the minimum strength (magnitude) of perturbation required to push the system into an alternative stable state. Suppose the system is initially in a stable state and experiences a perturbation of sufficient amplitude, regardless of its duration. In that case, the landscape is deformed enough that the system transitions to a different basin of attraction and does not return to the original state after the perturbation ends. The critical perturbation amplitude at which this transition occurs quantifies the system’s resistance.
    (d) Metrics of the basin of attraction. The two main properties are the flatness (inverse of curvature at the stable point) and the depth (distance from the valley bottom to the nearest hilltop, or unstable point). Flatter, shallower valleys correspond to slow recovery (low resilience) and lower resistance to shifts, while steeper, deeper valleys are associated with faster recovery (high resilience) and greater resistance to state shifts.
    }
    \label{fig 1}
\end{figure}

Stability landscapes and hysteresis loops are now standard tools for thinking about regime shifts in dynamic systems~\cite{scheffer2009critical}. In this view, a system’s state traverses a landscape of valleys and hills~\cite{scheffer2001catastrophic}: the shape of a valley encodes resilience as the recovery rate, the depth of a valley encodes resistance to switching, and fold bifurcations create hysteresis between forward and backward paths (see Fig.~\ref{fig 1}). Much empirical work interprets slow recovery, increased variance, or abrupt state changes in terms of this landscape metaphor and its associated early warning signals~\cite{dakos2022MainRef}.

However, real systems ranging from forests and lakes to microbial communities often behave in ways that strain this static picture. Recovery can be sluggish long after stress removal, with disturbances leaving ecological memory and legacies that slow dynamics for years or decades in vegetation and carbon cycling~\cite{egbert_scheffer2007, Anderegg2015}. Soil and host-associated microbiomes likewise retain altered composition and function long after environmental drivers return near their original range~\cite{jurburg2017legacy, martiny2016}. Collapses may still occur after conditions appear to have returned to normal, or altered states may persist, as legacy effects accumulate and push macroscopic ecosystems across thresholds with delay~\cite{Muller2022}, and similar delayed reorganizations are observed in microbial communities~\cite{bush2017oxic}. Hysteresis loops in slowly forced lake and plankton systems are often broader and more irregular than idealized bifurcation diagrams suggest~\cite{Schindler_Dilemma2012, JANSSEN2014813}, and model microbial communities show analogous multistability and hysteretic transitions~\cite{Dubinkina2019}. At the same time, applications of variance, autocorrelation, and related early warning indicators in theoretical and empirical ecosystems show mixed performance~\cite{Boettiger2012}, and in lake and microbial systems, these indicators can rise without an ensuing transition or remain weak even when a collapse occurs~\cite{Gsell2016, laitinen2021}. Together, these ecological and microbial examples suggest that something important is missing from the usual static stability landscape metaphor. 

One candidate for this missing ingredient is memory effects~\cite{khalighi2025impact}: the influence of past states on current dynamics. Our earlier work showed that memory, modeled via fractional time derivatives, can modulate resilience and resistance in interacting communities~\cite{Khalighi2022Ploscb}. Here, we build on that insight in a simpler setting, using one-dimensional bistable systems to ask how memory reshapes stability landscapes themselves, and how this in turn alters resilience, resistance, hysteresis, and the interpretation of early warning signals. Rather than treating the landscape as fixed, we reconstruct it from trajectories and show that memory makes the landscape history dependent and time evolving.

Specifically, in this paper, we address four points. First, we examine how adding long-term memory changes basic geometric features of stability landscapes, such as basin depth and flatness, and what tradeoffs between recovery speed and perturbation thresholds this implies. Second, we examine how memory alters the way exogenous fluctuations are transmitted through the landscape, from within-basin variability to switching statistics across time scales. Third, we examine what new critical outcome scenarios emerge in a moving landscape when slow internal parameter changes interact with external fluctuations (i.e., bifurcations), and how these scenarios differ from those in memory-free models. Finally, we examine how much the inferred equilibria, bifurcation structure, and tipping thresholds can be biased when a system with intrinsic memory is analyzed using a memory-free model. By addressing these points in univariate bistable systems, we aim to clarify how historical effects reshape stability landscapes and reinterpret familiar quantities such as resilience, resistance, and regime shifts.

\section*{Results}
Earlier work has shown that memory modulates resilience and resistance in ecological models~\cite{Khalighi2022Ploscb}. Here, we push that idea upstream and ask how memory reshapes the stability landscape itself. then trace how those geometric changes propagate to system recovery, regime shifting, hysteresis, and inference. We categorize drivers into two types~\cite{ashwin2012tipping}: \textbf{endogenous} perturbations that reshape the landscape and \textbf{exogenous} perturbations that nudge the state without altering the landscape parameters. This separation allows us to attribute each observed effect either to landscape deformation or to noise transmission through a landscape whose geometry and evolution are changes by memory.

\subsection*{Stability landscape perspective}
A method for understanding and visualizing the stability, resilience, and resistance of dynamical systems is through the use of potential energy landscapes, often referred to as \textbf{stability landscapes}~\cite{scheffer2009critical, scheffer2001catastrophic, dakos2022MainRef}.
This concept employs a metaphorical landscape composed of valleys and peaks, with a ball representing the system's state at any given time, Fig.~\ref{fig 1}a. The valleys denote stable states where the system naturally tends to settle. In contrast, the peaks represent unstable states, acting as barriers or thresholds that the system must overcome to transition between stable states. The movement of the ball across these landscapes illustrates shifts in stability, depicting how the system navigates through different regions of stability under various conditions.

Far from being purely conceptual, stability landscapes can be empirically reconstructed with appropriate data~\cite{dakos2022MainRef}. The key features of these landscapes, valley depth, flatness, and width, directly inform the dynamics of the system. Depth reflects resistance to state shifts, while flatness indicates the speed of recovery from perturbations, offering a tangible representation of resilience and resistance, Fig.~\ref{fig 1}b-d.

Therefore, to quantify these landscapes, we focus on two metrics (see also Supplementary S3): 1) \textbf{Potential depth}, the energy difference between a basin’s minimum and the adjacent peak, indicating resistance to transitions; and 2) \textbf{Potential flatness}, the inverse curvature at the basin minimum, reflecting slower recovery in flatter basins. We excluded other metrics, such as distance to the basin threshold~\cite{dakos2022MainRef}, due to strong correlations with depth (positive) and flatness (negative), which provided no additional insights.

\subsection*{Endogenous perturbation and memory effects}
To explore endogenous perturbations, we analyze pulses (temporary parameter shifts) through the lens of stability landscapes. This framework intuitively illustrates how internal changes reshape system dynamics. Below, we delve into the stability landscape and its connections to memory effects and system behavior, providing a perspective on their interplay.

\subsubsection*{Impact of memory on stability landscapes}

Sensitivity to initial conditions is crucial when interpreting stability landscapes. In memory-free systems, the reconstructed landscape is independent of the starting state: trajectories converge regardless of the initial conditions. With memory, by contrast, the dynamics are history-dependent so that identical parameters can yield different reconstructed landscapes for various initial conditions. Consequently, our examinations show that initial conditions from distinct values trace different trajectories and produce altered landscape profiles (see Fig.~\ref{fig: initial conditions}). To ensure consistency, all comparisons were performed using uniform initial conditions to yield identical landscapes and eliminate potential discrepancies.

We found that memory alters both stability landscape metrics, potential depth, and flatness (Figs.~\ref{fig S1} and \ref{fig: Dp_Curv}). Depending on parameter values, it can deepen or shallow the basins, but it consistently flattens them (Fig.~\ref{fig: Cor_Mem_dp_Fl}). The systematic flattening arises because memory slows the approach to equilibrium, so trajectories spend more time near the basin floor, and the reconstructed slope becomes gentler. In contrast, changes in depth depend on the baseline geometry: in flat and shallow basins, memory often reduces depth, whereas in steep and deep basins it more frequently increases it (Fig.~\ref{fig: Mem_fl_dp}). Using the basin sharpness index introduced in the supplementary, we show that the effect of memory on flatness is strongest for steep and deep basins, while the impact on depth is largest for flat and shallow ones. These patterns are not model-specific but reflect a general interaction between memory strength and landscape shape (see Supplementary S3; Figs.~\ref{fig: detailed results}-\ref{fig: summary results}).

\begin{figure}[ht!]
    \centering    \includegraphics[width=\linewidth]{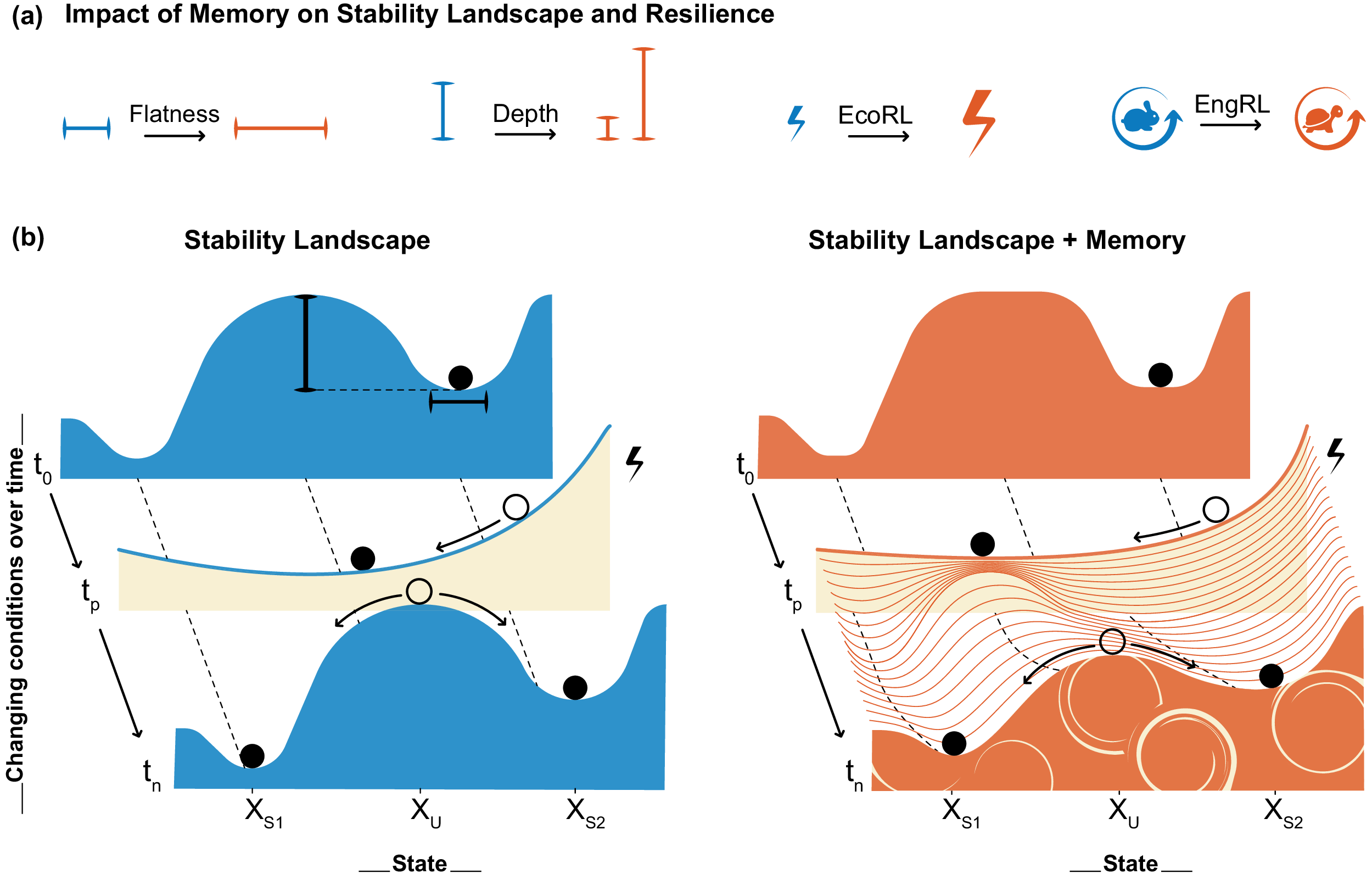}
    \caption{
    (a) Summary of the main findings, the impact of memory on stability metrics, resilience, and resistance of 1000 randomly parameterized polynomial models: memory can either increase or decrease potential depth, tends to flatten the bottom of the basin of attraction, slows recovery (decreasing resilience), and raises the endogenous perturbation threshold for switching to an alternative stable state (increasing resistance). 
    (b) Conceptual comparison of stability landscapes in one-dimensional bistable systems, with and without memory. Three phases are shown: before, during, and after a perturbation. $X_{S1}$ (set to 0 in this study) and $X_{S2}$ (positive; initial state) are stable states, while $X_{U}$ is an unstable state. $t_0$ marks the initial time or endogenous perturbation onset (the derivative is zero at this point and does correspond to the systematic flattening of the basin floor under memory, whereas the size of the resistance gain depends on the turns to the original stable state or transitions to an alternative one. In memory-free systems, the landscape remains unchanged before and after perturbation. In contrast, systems with memory exhibit evolving landscapes during and after the perturbation, which alter the system’s response.
    }
    \label{fig 2}
\end{figure}

\subsubsection*{Impact of memory on resilience and resistance}
Changes in depth and flatness of the stability landscape directly shape the response to perturbations. Using short endogenous pulses, we quantify resilience and resistance. Based on numerical results, we found that memory consistently slows recovery (lower resilience) while requiring larger or longer pulses to tip the system (higher resistance). Fig.~\ref{fig 2}b summarizes the ensemble-level effects on depth, flatness, resilience, and resistance. The universal slowdown tracks the systematic flattening of the basin floor under memory, whereas the size of the resistance gain depends on baseline geometry (Supplementary S3).

Geometry clarifies when effects are strongest. In steep and deep basins, memory most strongly increases flatness, which translates into the largest loss of resilience; in flat and shallow basins, memory most strongly changes depth and yields the largest gain in resistance. Thus, memory induces an inverse relation between resilience and resistance at fixed parameter values: as memory grows, recovery becomes slower even as tipping becomes harder.

For definitions and measurement, see the resilience index schematic in Supplementary S3 (Fig.~\ref{fig: RL}). Case-by-case demonstrations are given for the ecological model in Figs.~\ref{fig: mem_herb_RL1}–\ref{fig: land_herb_RS}. The dependence on landscape geometry is established with the ensemble summary in Fig.~\ref{fig: mem_RL_RS}, and the basin sharpness analysis, with continuous trends detailed in Figs.~\ref{fig: detailed results}-\ref{fig: Mem_RL_RS_Basin} and \ref{fig: summary results}-\ref{fig: dp_fl_RL_RS}. Together, these results show that the resilience loss is a robust consequence of memory-driven flattening, while the resistance gain is largest in basins that are initially flat and shallow.

\subsubsection*{Dynamic stability landscapes: Memory-induced aftereffects}
A key implication of the landscape view is that memory makes the landscape itself dynamic (Fig.~\ref{fig 2}a; also Supplementary S2). Depth and flatness can continue to evolve after the pulse ends, causing the tipping point to shift and settle asymptotically back to its original position, so the familiar fixed surface intuition can fail. We highlight two deterministic, repeatable behaviors that follow from this moving landscape picture and that can be misread as noise-driven effects.

\textit{Rollback after near crossing:}
A short pulse can push the state just past the instantaneous unstable point. In a memory-free model, stopping the pulse would restore the original landscape and typically trigger an immediate commitment to the new basin. With memory, the landscape keeps changing after the parameter is reset. The barrier can reform behind the state while the state is still drifting near the ridge, allowing the system to roll back into the original valley (Figs.~\ref{fig: dyn_pert1}-\ref{fig: dyn_pert2}; also animations in Videos S1-S2). This rollback is slow and paced by landscape evolution rather than rapid relaxation, which can make the state appear stuck near the threshold.

\textit{Delayed collapse and delayed recovery:}
Because landscape reshaping persists after stress removal, collapse or recovery can occur long after the control parameter returns to baseline (Figs.~\ref{fig: dyn_pert1}-\ref{fig: dyn_pert2}; also animations in Videos S1-S2). Whether the system ultimately collapses or recovers depends on the relative timing of the slow landscape motion and the state position near the barrier. These delayed outcomes are predictable signatures of memory-driven landscape drift and motivate tracking not only state variables but also slow variables or cumulative stress that reveal landscape motion.

\subsection*{Exogenous perturbation and memory effects}
We then study exogenous perturbations, modeled as external fluctuations that nudge the state without directly changing the parameters that shape the landscape (Supplementary S4). Memory changes how these fluctuations are transmitted, producing time-dependent signatures in variability and correlation (Figs.~\ref{fig:traj}-\ref{fig:persist}).

\subsubsection*{Within basin variability across time scales}
Holding noise level fixed, adding memory makes the time series appear rougher at fine resolution (i.e., larger point-to-point changes over short time steps; Fig.~\ref{fig:traj}).
The system retains the impact of fast disturbances instead of smoothing them out, which increases variance and reduces short-term similarity (Figs.~\ref{fig:traj} and \ref{fig:acf}). In the frequency domain, this corresponds to weaker attenuation of high-frequency components under memory (Fig.~\ref{fig:psd}). When the data are averaged over short windows, slower dynamics dominate, and the memory case appears more persistent than the memory-free case. This scale crossover is visible in the variance after coarse graining (Fig.~\ref{fig:acf}) and in the windowed integrated autocorrelation time, which becomes larger under memory once modest averaging removes the fastest fluctuations (Fig.~\ref{fig:persist}). This shows that the apparent effect of memory on variability depends on the temporal scale at which the time series is observed (sampling resolution and averaging window; Figs.~\ref{fig:traj} and \ref{fig:persist}).

\subsubsection*{Committed transitions and dwell time}
Under exogenous noise, trajectories may cross the unstable boundary repeatedly. If every boundary crossing is recorded, transitions are overcounted, especially with memory, which promotes near-threshold excursions that hover around the boundary equilibrium without settling into the opposite basin (Fig.~\ref{fig:trajswitch}). 

\begin{figure}[ht!]
    \centering
    \includegraphics[width=0.8\linewidth]{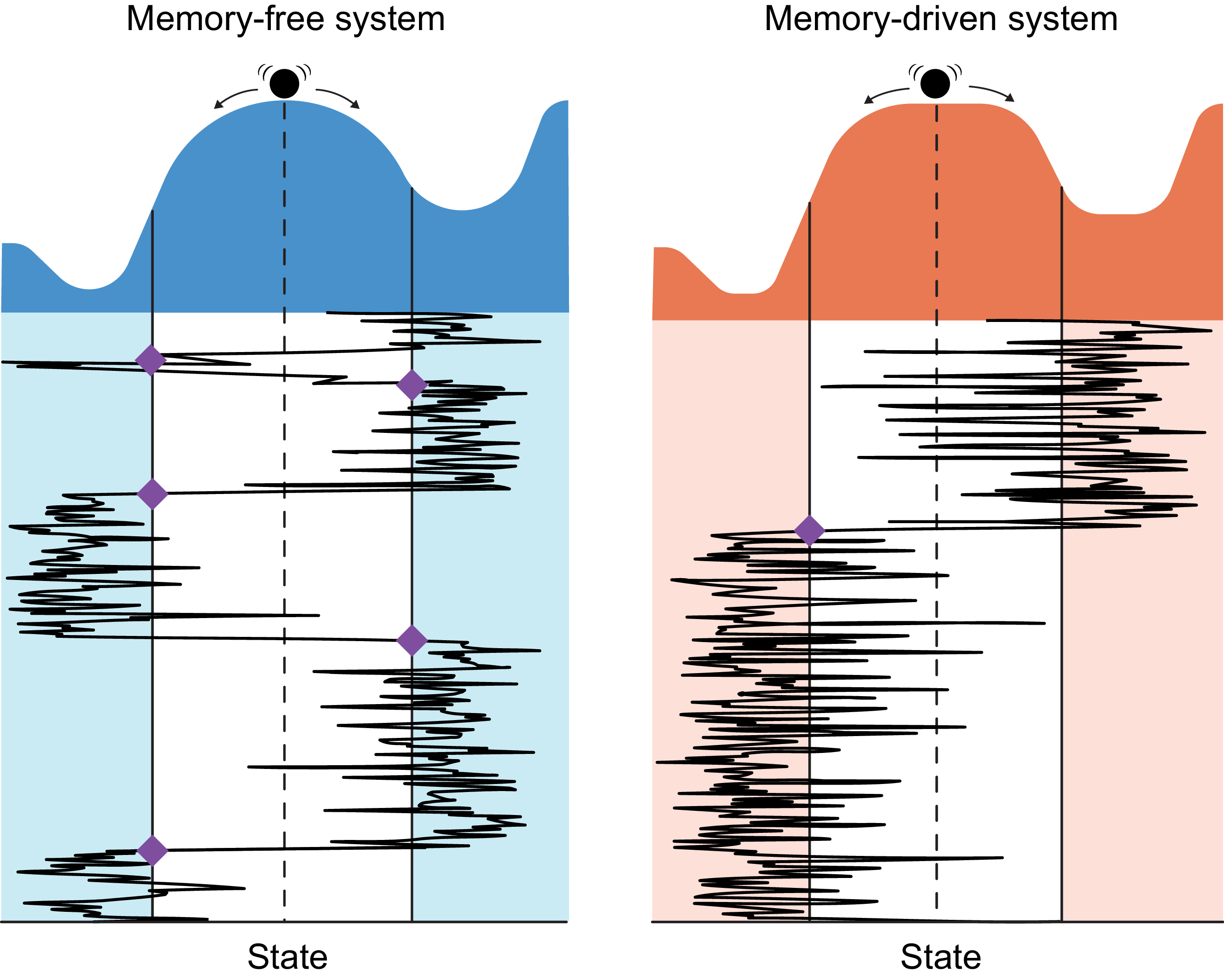}
      \caption{Schematic of committed transition detection for a noisy bistable system. Time increases vertically, and the horizontal axis is the state. Left: memory-free case with five committed transitions. Right: memory-driven case with one committed transition. Diamond mark detected onsets of alternating committed transitions.}
    \label{fig 3}
\end{figure}

We therefore distinguish raw crossings from committed transitions: a transition is recorded only when the state reaches and remains within the core of the opposite basin. This filters out near-threshold rattling and enables a fair comparison between memory-free and memory-driven dynamics. Figure~\ref{fig 3} illustrates the detection procedure and shows that committed transition counts can differ from raw crossings. 

Using committed transitions, memory broadens dwell-time distributions and makes them sensitive to the observation window. This appears as heavier tails in dwell times (Fig.~\ref{fig:residhist}) and slower decay in survival curves (Fig.~S38). Because the net switching frequency can increase or decrease depending on basin geometry and noise timing, a single escape rate is not an adequate summary; we therefore report dwell-time distributions rather than compressing switching into one rate parameter (Supplementary S4).

\subsubsection*{State dependent noise: multiplicative case}
If the noise amplitude itself depends on the state, the qualitative memory signatures remain, but the ordering of global sample variance can become model-dependent. Near a given stable state, the local noise level governs the immediate spread, while memory still lengthens correlation and transmits faster content (see Figs.~\ref{fig:multiseg}-\ref{fig:trans-mult}). However, because the trajectory with memory may spend more or less time in parts of the basin where the noise amplitude is larger or smaller, the overall variance need not be monotone in memory. This does not contradict the mechanism above; it reflects how state-dependent forcing and memory jointly shape the empirical spread. 

While these state-dependent effects offer a rich area for study, the remainder of our results focuses on additive noise to isolate the fundamental influence of memory on landscape stability. For this interpretation, we distinguish within-basin statistics from committed transitions and adopt a consistent stochastic framework for all following analyses.

\subsection*{Critical outcome scenarios under combined drivers}
When internal and external drivers act together, memory couples their effects in ways that matter for critical transition and recovery. A temporary parameter shift can continue to reshape the landscape even after the shift ends, while external perturbations persistently nudge the state. Depending on when the slow landscape motion completes relative to the state position, memory can delay or accelerate the apparent outcome. We identified three repeatable patterns in our models and examined them in our biological case study (see Models and methods).

\subsubsection*{Delayed or accelerated outcomes after stress removal}
Because the landscape continues to drift for a while, collapse or recovery can occur long after the parameter returns to baseline. External perturbations can either assist a late escape or a late return, depending on where the state sits relative to the slowly moving equilibrium boundary during this drift and on the observation window.

Yet for sufficiently brief pulses, memory can also facilitate a faster recovery than in a memory-free system. Consider two stable states with weak external noise. A short internal pulse may carry the memory-free trajectory past the unstable point into the other basin, where it quickly relaxes. With memory, the same pulse can leave the state near the moving ridge while the landscape continues to evolve; as the landscape's hill reforms behind the state after the pulse ends, the state slides back into the original valley. Because the memory trajectory never fully relaxes into the new basin, the return completes earlier than in the corresponding memory-free case that already switched. 
 
In our simulations, we illustrate these dynamics through animations: Video S4 shows a delayed collapse caused by memory, while Videos S5 and S6 demonstrate delayed and rapid recoveries following a regime collapse.

\subsubsection*{Rollback after a near crossing}
A short pulse may push the state just beyond the momentary unstable point. As the landscape readjusts, the landscape's hill can reform behind the state and, in the presence of memory, the state drifts back toward the original valley. External perturbations modulate timing, but the rollback is a deterministic consequence of landscape motion with memory (see Video S7). 

\begin{figure}[ht!]
    \centering
    \includegraphics[width=0.8\linewidth]{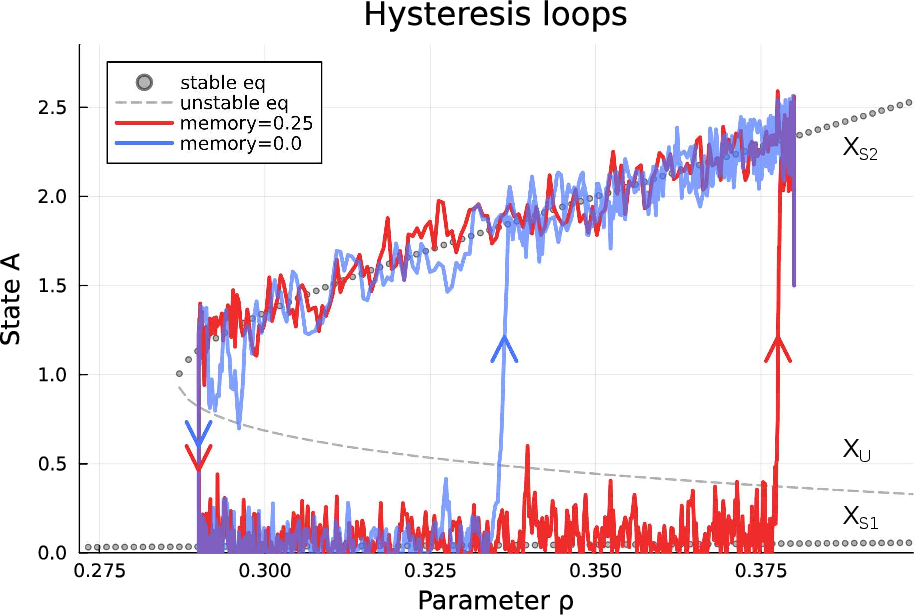}
  \caption{Memory broadens hysteresis in the quorum sensing model \eqref{eq: quorum1} (see Models and methods). It shows system trajectories on the bifurcation diagram with stable equilibrium branches ($X_{S1}$ and $X_{S2}$) shown in gray dots and the unstable branch ($X_U$) shown as a dashed curve. The system experiences both exogenous perturbations and gradual shifts of the endogenous control parameter $\rho$: it decreases from 0.38 to 0.29 to induce collapse, then returns to 0.38 for recovery. Red traces correspond to a system with memory equal to 0.25, and blue traces to a memory-free system. Circle markers indicate the collapse leg (as $\rho$ moves from 0.38 to 0.29), and triangle markers indicate the recovery leg (as $\rho$ returns from 0.29 to 0.38). With memory, the collapse and recovery paths separate more, so returning to the upper stable branch $X_{S2}$ requires stronger or longer favorable conditions.}
    \label{fig: widen hysteresis}
\end{figure}

\subsubsection*{Broadened hysteresis under gradual parameter shifts}
Under gradual parameter shifts, memory widens the range over which an apparent critical transition can occur. External fluctuations further spread the observed thresholds, so the forward and backward branches of the loop are more diffuse than in memory-free models (see Fig.~\ref{fig: widen hysteresis}; and Videos S8).

These scenarios demonstrate that memory can both impede and promote transitions, broaden hysteresis, and yield delayed collapse or faster and slower recovery, challenging memory-free interpretations and guiding monitoring and management.

\subsection*{Inference risk under memory mismatch}
A practical challenge is that memory may be present in real systems even when it is not explicitly modeled. We therefore ask what happens when memory-driven trajectories are fit with a memory-free model with adjusted parameters~\cite{ghosh2025fractional}. The key result is that trajectory matching does not guarantee correct inference of equilibria and tipping thresholds. 

We generate deterministic trajectories from a model with memory and then fit them using the same model structure without memory. The fitted model can reproduce the observed trajectories and qualitative convergence patterns, which might suggest that it also recovers the correct stability picture (see Supplementary S5; Fig.~\ref{fig: fit-data}). 

\begin{figure}[ht!]
    \centering
    \includegraphics[width=.75\linewidth]{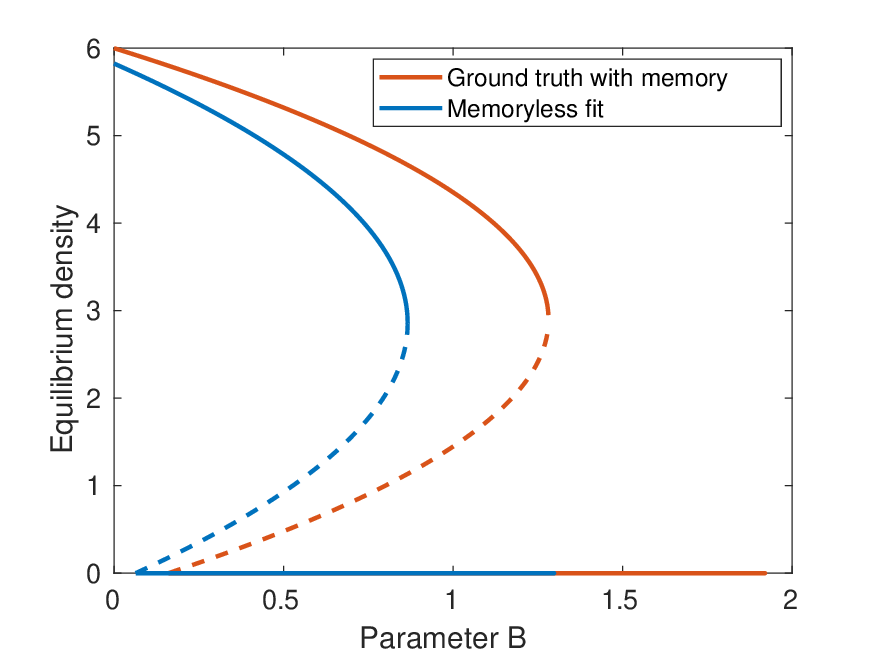}
    \caption{Misplaced bifurcation structure when neglecting memory. Red curves (baseline, memory-driven system that generated the data) show stable (solid) and unstable (dashed) equilibria versus the control parameter \(B\). Blue curves depict the bifurcation diagram of the best memory-free fit to those time series. Although the fit reproduces trajectories well, it shifts the equilibrium branches along \(B\), moving the tipping points to lower \(B\) and distorting (shrinking) the bistable interval. Thus, thresholds inferred from the fitted model are incorrect.}
    \label{fig:fit}
\end{figure}

Despite good trajectory fits, the memory-free calibration yields a different bifurcation diagram (Fig.~\ref{fig:fit}; also Figs.~\ref{fig:  land_fit_data}–\ref{fig:  land_fit_data2}). Equilibrium branches shift along the control parameter axis, tipping points move to lower parameter values, and the bistable interval shrinks. As a result, thresholds inferred from the fitted model are incorrect even though time-series agreement appears strong. 

This mismatch highlights an inference risk: ignoring memory can systematically bias conclusions about stability and tipping points. In applied settings, it is therefore important to treat landscape-based diagnoses as potentially moving and to test whether observed dynamics show signatures consistent with memory before trusting thresholds inferred from memory-free fits (see Supplementary S5 and Fig. \ref{fig:  fit data pulse}). 

\section*{Discussion}

\subsection*{Key contributions}

Long-term memory reshapes bistable stability landscapes. Adding memory (fractional derivatives) flattens basin floors (lower curvature) and can deepen or shallow basins depending on baseline geometry. These changes map to classical metrics: flatter floors slow recovery (reduced resilience), while deeper/steeper basins raise the disturbance needed to escape (greater resistance).

Because the landscape itself becomes history-dependent and time-evolving, dynamics can defy memory-free intuition. Short pulses may push the state past the unstable point; however, if the landscape continues to shift after the pulse, the barrier hill may reform ``behind'' the state, forcing the system to roll back. Conversely, collapse or recovery can occur long after stress is removed, reflecting delayed landscape reshaping rather than noise.

When endogenous (internal, parameter-driven) and exogenous (external, noise-driven) perturbations act together, memory couples their effects in complex ways. We identified three distinct critical outcome scenarios arising from this interplay: 1) Delayed or accelerated outcomes after stress removal. 2) Rollback after a near crossing. 3) Broadened hysteresis under gradual parameter shifts.

These novel behaviors do not appear in memory-free models and highlight how memory enriches the repertoire of system responses. Observers might otherwise attribute such phenomena to hidden stochasticity or external factors, but our study shows they can arise purely from internal memory. This motivates moving beyond the static ball-in-cup picture toward modeling and diagnostic frameworks that can distinguish memory-induced landscape evolution from external noise, enabling correct attribution of observed anomalies to their underlying mechanism.

\subsection*{Implications for modeling and inference}
A major applied risk is that memory-driven systems can be fit well by memory-free models while producing incorrect stability conclusions. Our fitting analysis shows that memory-free calibration can reproduce trajectories but still shift equilibrium branches, move tipping points along the control parameter axis, and distort the bistable interval. This means that even accurate time series fits can mislead landscape-based diagnosis and management decisions when memory is present but ignored.

\subsubsection*{Practical cues and cautious interpretation}
A moving landscape view suggests practical ways to reason about data. Sluggish recovery or persistent fluctuations near a threshold can reflect memory rather than noise, and interventions can have delayed effects because landscape reshaping continues after drivers return to baseline (see Supplementary S4). Diagnostics that explicitly compare time scales can help, including spectral signatures, scale-dependent variance, and dwell time analysis, together with explicit tests for memory.

\subsubsection*{Reinterpreting early warning signals}
Many frameworks view rising variance, increased autocorrelation, and spectral reddening as signatures of critical slowing down near a bifurcation~\cite{Boettiger2012}. Our results show that memory effects can generate similar patterns for entirely different reasons. By flattening basin floors, memory slows recovery even far from a tipping point; by altering noise transmission, it can simultaneously increase variance at native resolution while reducing it after smoothing. Furthermore, memory broadens dwell-time distributions (for additive noises) and couples internal and external drivers, allowing delayed collapses or recoveries to occur outside the parameter range predicted by standard theory. Consequently, classical indicators~\cite{laitinen2021} should be interpreted with caution: apparent signals may reflect evolving stability landscapes rather than proximity to a static critical point. Distinguishing true critical slowing down from memory-driven dynamics, therefore, requires combining time-scale explicit diagnostics (such as spectra and scale-dependent variance) with explicit tests for memory.

\subsection*{Scope, limitations, and future directions}

Our analysis focused on univariate gradient systems where a potential landscape exists. Real systems are often multivariate and non-gradient, where a single scalar potential may not apply~\cite{Rodriguez2020EscherStairs, Zhou2012Quasi}. We also considered one specific memory form: a power-law kernel~\cite{kilbas2006theory,podlubny1998fractional}. Other kernels, exponential, oscillatory, or discrete delay, may produce different behaviors~\cite{Yang2023}. The reconstructed potentials are path-specific diagnostic constructs, not true energy functions, meaning they depend on the system’s trajectory. 

Moreover, quantitative definitions of resilience and resistance depend on chosen thresholds and initial conditions (see supplementary S6). 
For instance, while we find that memory generally increases resistance, it can also have the opposite effect near an unstable equilibrium. In such cases, slow, memory-driven motion can keep the state lingering at the threshold, allowing even small perturbations to trigger a transition—effectively reducing local resistance (Fig.~\ref{fig: s38}). Finally, our study is theoretical and simulation-based; empirical detection of fractional memory remains challenging due to measurement noise, parameter drift, and the difficulty of estimating the memory order.

As next steps, we suggest validating these memory-induced dynamics in real systems—using targeted data analyses and controlled experiments to probe delayed collapses, accelerated recoveries, and related effects—while extending the framework to higher-dimensional and spatially coupled models to uncover network-level consequences. Parallel work should compare alternative memory forms~\cite{khalighi2025impact} (e.g., power-law, exponential, discrete delays~\cite{Yang2023}) to separate universal from kernel-specific behaviors, and develop data-driven tools to infer both the presence and the type of memory from observations, leveraging diagnostics like power spectra and dwell-time statistics alongside machine learning. Finally, revisiting stability metrics, mean exit times~\cite{Arani_exit_time_2021}, critical slowing down, and basin stability, under memory could help connect deterministic pictures with probabilistic resilience in a unified way.

In summary, incorporating historical memory effects into modeling reveals richer dynamical regimes and cautions against assuming fixed tipping points or rapid recoveries. The study offers both conceptual and methodological foundations for future interdisciplinary work on memory in complex systems.

\section*{Models and methods}

Shortly explain the models, population dynamics models, and then refer to the details in the supplementary.

\subsection*{Population dynamics models}
We use one-dimensional population dynamics models, where the population size $x(t)$ changes over time according to 
\[\frac{dx}{dt}=F(x,t),\]
where $x(t)$ is the population at time $t$, and $F(x,t)$ is the net growth rate (births, deaths, and other factors) as a function of population and time.
\subsection*{Polynomial models}
For this study, we consider $F(x)$ is a third-degree (cubic) polynomial 
\[F(x)=ax^3 + bx^2 +cx +d\]
to give us a bistable system, because it is cubic, it can have up to three real equilibria ($F(x)=0$). With suitable parameter values, the cubic curve can cross the x-axis three times, so the outer two equilibria can be stable and the middle equilibrium is typically unstable, perfect for our study to examine a bistable system behavior. If the population starts on one side, it moves toward the lower stable equilibrium; if it starts on the other side, it moves toward the higher stable equilibrium
A small perturbation near the middle unstable equilibrium can push the system to the other stable state.

We can define its potential function $V(x)$ such that~\cite{Rodriguez2020EscherStairs}:

\[\frac{dx}{dt}=F(x)=-\frac{dV}{dx}\]
then:
\[V(x)=-\int F(x)dx\]

for the cubic function:
\[V(x)=-\int (ax^3+ bx^2 +cx +d)dx = -\left(\frac{a}{4}x^4 +\frac{b}{3}x^3+\frac{c}{2}x^2+dx\right)\]

So, $V(x)$ is a fourth-degree (quartic) polynomial. Stable equilibria of 
$F(x)$ (where $F=0$ and $F'<0$) correspond to local minima of $V(x)$.
The unstable equilibrium corresponds to a local maximum of $V(x)$.
So, the double-well shape of the potential landscape represents the two alternative stable population states (bistability) separated by an energy barrier (the unstable point). Small noise can push the population over the barrier, causing it to switch between states.

\subsection*{Ecological and biological models}

Many ecological and biological models can be written as (Supplementary S1)
\[
\frac{dx}{dt}=\frac{p(x)}{q(x)},
\]
with $\deg p=3>\deg q$ and $q(x)>0$ on the domain of interest, so the qualitative dynamics are governed by the zeros and signs of $p(x)$ \cite{dakos2022MainRef,saade2023MainRef}. This form supports bistability and standard potential landscape analysis via $V(x)=-\int p(x)/q(x)\,dx$. In the supplementary, we classified typical ecological and biological models that follow this format. We employ the following two models in this form for our analyses.

\paragraph{Ecological case: herbivory model~\cite{dakos2022MainRef, saade2023MainRef}.}
A typical desertification model is
\begin{equation}
\frac{dx}{dt}=rx\!\left(1-\frac{x}{K}\right)-\frac{Bx}{A+x}
\end{equation}
where $x$ is vegetation biomass,
$r$ intrinsic vegetation growth rate,
$K$ carrying capacity,
$B$ herbivore pressure (density or grazing intensity),
$A$ half saturation constant of consumption.
It can be written as
\[
\dot x=\frac{p(x)}{q(x)},\quad
p(x)=a_3x^3+a_2x^2+a_1x+a_0,\quad
q(x)=A+x,
\]
with coefficients
\[
a_3=-\frac{r}{K},\qquad
a_2=r\!\left(1-\frac{A}{K}\right),\qquad
a_1=Ar-B,\qquad
a_0=0.
\]
Here $q(x)>0$ for $x\ge 0$ when $A>0$.

\paragraph{Biological case: quorum sensing toy model.}
Following Dockery and Keener’s reduced autoinducer model~\cite{DOCKERY200195},
\begin{equation}\label{eq: quorum1}
\frac{dx}{dt}=\frac{V x^2}{K+x^2}+x_0-\!d(\rho)\,x,
\end{equation}
where $x$ is the autoinducer concentration,
$V$ the maximal synthesis rate due to positive feedback, $K$ saturation constant in the term $x^2/(K+x^2)$, $x_0$ basal production rate, $\rho$ cell density level,
$d(\rho)$ effective loss rate (dilution or degradation) depending on $\rho$.
We obtain
\[
\dot x=\frac{p(x)}{q(x)},\quad
q(x)=K+x^2>0,\quad
p(x)=a_3x^3+a_2x^2+a_1x+a_0,
\]
with
\[
a_3=-d(\rho),\qquad
a_2=V+x_0,\qquad
a_1=-d(\rho)\,K,\qquad
a_0=x_0 K.
\]
This cubic numerator with a positive quadratic denominator captures the bistable switch as $\rho$ varies, while keeping analysis focused on the roots and slopes of $p$.

\subsection*{Stochastic forcing and stochastic differential equations}
To represent exogenous perturbations, we model the system state as a stochastic process governed by a stochastic differential equation (SDE). In one dimension, we write the dynamics in It\^{o} form as~\cite{KloedenPlaten1992, nolting2016balls}
\begin{equation}
  dX_t = f(t, X_t)\,dt + g(t, X_t)\,dW_t,
\end{equation}
where $X_t$ is the state, $f(t,X)$ is the drift term that captures the deterministic tendency of the system, $g(t,X)$ is the diffusion term that controls the magnitude of random fluctuations, and $W_t$ is a standard Wiener process. Additive noise corresponds to constant diffusion, $g(t,X)=\sigma$, while multiplicative (state-dependent) noise allows $g(t,X)$ to vary with the state, for example $g(t,X)=\sigma(X)$. More details are provided in Supplementary S4.

The drift governs the mean direction and rate of motion, whereas the diffusion governs variability, short-time roughness, and the probability of rare landscape boundary crossings. The probability density $p(x,t)$ associated with the SDE satisfies the Fokker-Planck equation~\cite{nolting2016balls}
\begin{equation}
  \frac{\partial p}{\partial t}
  = -\frac{\partial}{\partial x}\!\big(f(t,x)p(x,t)\big)
    + \frac{1}{2}\frac{\partial^2}{\partial x^2}\!\big(g^2(t,x)p(x,t)\big).
\end{equation}
In gradient like cases with constant diffusion, where $f(x)=-U'(x)$ for a potential function $U(x)$, the stationary density (when it exists) is proportional to $\exp\!\big(-2U(x)/\sigma^2\big)$, linking stochastic fluctuations to an effective stability landscape.

\subsection*{Quantifying memory with fractional derivatives}

By \textbf{memory} we mean the influence of past states, such as earlier population densities, on the current dynamics. We model this influence by replacing ordinary derivatives with fractional derivatives~\cite{Khalighi2022Ploscb, khalighi2025impact}. The order of the derivative, denoted $\alpha$ and constrained to $0<\alpha\le 1$, provides a measure of memory strength:
\[
\text{memory strength}=1-\alpha .
\]
A value $\alpha=1$ recovers the usual memory-free system, while a smaller $\alpha$ implies a stronger influence of history.

Mathematically, a fractional derivative is a convolution of the state variable with a memory kernel that decays following a power law, so the entire past contributes, but recent states carry more weight. Thus, a population equation that ordinarily reads
\[
\frac{dx}{dt}=F(x,t)
\]
becomes
\[
\mathcal{D}^{\alpha}x = F(x,t),
\]
where $\mathcal{D}^{\alpha}$ is Caputo fractional derivative of order $\alpha$~\cite{kilbas2006theory,podlubny1998fractional}. As $\alpha$ decreases, the kernel spreads over longer times, and memory increasingly shapes the evolution of the system. Details are explained in the Supplementary S2.

\subsection*{Evaluating the potential landscape with memory}

Our goal is to visualize a potential landscape that is consistent with models where memory matters. When the fractional order is $\alpha<1$, the drift at time $t$ depends on the past through a power-law kernel. We therefore evaluate the potential along system trajectories rather than from a closed form $V'(x)=-F(x)$ which only holds for $\alpha=1$.

Consider the fractional model $\mathcal{D}^{\alpha}x=F(x,t)$ with $0<\alpha\le 1$. For any trajectory $Z(t)$ of this model, the relative potential difference between two states $X_a$ and $X_b$ is defined by
\[
V\;=\;-\int_{X_a}^{X_b} \dot Z(t)\, dZ,
\]
where $\displaystyle Z(t)$ is generated by the fractional dynamics. For $\alpha=1$ this reduces to the classical potential because $\displaystyle \dot Z = F$, depending only on the current state. 

In practice, we need to find equilibrium points of $\displaystyle  \mathcal{D}^{\alpha}x=F(x,t)$ and label the two stable states $\displaystyle X_{S_1}, X_{S_2}$ and the unstable state $\displaystyle X_U$ with $\displaystyle  X_{S_1}<X_U<X_{S_2}$.
We integrate the fractional model from initial conditions just to the left and right of $X_U$ to obtain trajectories $Z(t)$ that sweep the valleys on each side.
Thus, we should estimate $\displaystyle dZ/dt$ along each trajectory and use it to evaluate the integral $\displaystyle -\!\int \dot Z(t)\, dZ$ numerically across the state grid. We set an arbitrary zero level, for example $V(X_U)=0$, and plot $V$ piecewise on $\displaystyle (X_{S_1}, X_U)$ and $(X_U, X_{S_2})$. Algorithmic details and accuracy choices for numerical differentiation and interpolation are given in the Supplementary S2.

\bibliographystyle{unsrt}
\bibliography{references}

\clearpage
\renewcommand{\thefigure}{S\arabic{figure}}
\setcounter{figure}{0}

\section*{\huge{Supplementary material}}

\section*{S1: Representation of ecological and biological models using polynomial quotients}\label{sec: sup samples}

Building upon the framework proposed by Dakos et al.~\cite{dakos2022MainRef} and Saade et al.~\cite{saade2023MainRef}, we aim to categorize ecological and biological models that can be expressed as quotients of polynomials, $\frac{p(x)}{q(x)}$, where $p(x)$ is a third-degree polynomial:
\begin{equation}\label{eq: poly}
p(x) = a_3 x^3 + a_2 x^2 + a_1 x + a_0,
\end{equation}
and the degree of $q(x)$ is lower than that of $p(x)$. This formulation simplifies analysis by allowing us to focus on the roots and signs of $p(x)$, which determine the qualitative dynamics of the original models. By ensuring $q(x)$ remains positive and monotonic—often a first- or second-degree polynomial—we can analyze system behavior by concentrating on $p(x)$.

\begin{figure}[ht!]
    \centering
    \includegraphics[width=1\textwidth]{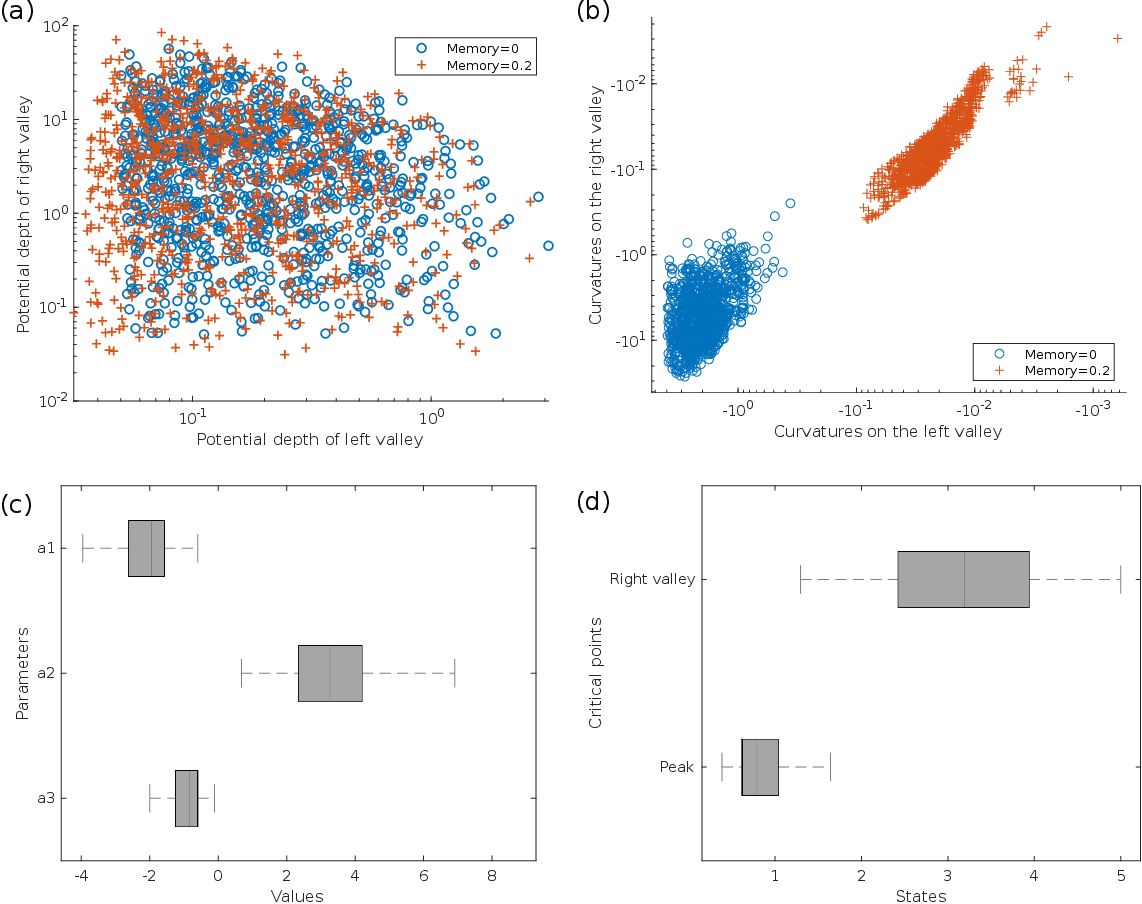}
    \caption{{(a)} The scatter plot of the potential depths generated by 1000 randomly parameterized polynomial models. The red crossed markers represent models with memory, while the blue circles represent models without memory. {(b)} Similarly, the scatter plot of the potential curvatures from the left and right valleys of the potential energy for the same set of polynomial models. {(c)} The distribution of values used for the coefficients a1, a2, and a3 in the sample models~\eqref{eq: poly}. {(d)} The distribution of the locations of the basin of attraction (the positive stable points) and the hilltop (unstable points) is depicted in this panel.}
    \label{fig: samples}
\end{figure}

We generated an ensemble of 1000 polynomial models with varying shapes, depths, and curvatures (see Figure~\ref{fig: samples}a–b). Initially, all models are simulated without memory, and we then introduce memory (by setting the memory strength to \(1 - \alpha = 0.2\); see Section S2).

The model coefficients are selected as follows (see Figure~\ref{fig: samples}c):

\begin{itemize}
    \item[-] \(a_0 = 0\), omitting the constant term;
    \item[-] \(a_1\) and \(a_3\) are negative values, each independently sampled from the uniform distribution from a specified interval;
    \item[-] \(a_2\) is a positive value, selected from the interval
    \(
        a_2 \in \left[ \sqrt{4a_1a_3},\ \sqrt{4a_1a_3} + 4 \mathcal{\eta} \right],
    \)
    where \(\mathcal{\eta}\) is a uniformly distributed random variable on the interval \([0, 1]\).
\end{itemize}

To clarify this $p/q$ representation, consider the classic herbivory model~\cite{dakos2022MainRef, saade2023MainRef}:
\begin{equation}\label{eq: herbivor}
\frac{dx}{dt} = rx\left(1 - \frac{x}{K}\right) - \frac{Bx}{A+x}.
\end{equation}
This equation can be rewritten in the $p/q$ format~\cite{saade2023MainRef}, where:
\[
p(x) = a_3 x^3 + a_2 x^2 + a_1 x + a_0,
\]
with $a_1 = Ar - B$, $a_2 = r(1 - A/K)$, $a_3 = -r/K$, and $q(x) = A + x$.

Likewise, Dockery and Keener’s autoinducer model~\cite{DOCKERY200195},
\begin{equation}\label{eq: quorum}
\frac{dx}{dt}=\frac{V x^2}{K+x^2}+x_0-\!d(\rho)\,x,
\end{equation}
can be written as:
\[
\dot x=\frac{p(x)}{q(x)},\quad
q(x)=K+x^2>0,\quad
p(x)=a_3x^3+a_2x^2+a_1x+a_0,
\]
with
\[
a_3=-d(\rho),\qquad
a_2=V+x_0,\qquad
a_1=-d(\rho)\,K,\qquad
a_0=x_0 K.
\]

This example demonstrates the broad applicability of the polynomial quotient framework to ecological and biological models. By expressing model dynamics in the form $\frac{p(x)}{q(x)}$, analysis is simplified: one can investigate critical points, stability, and dynamic behavior primarily by examining $p(x)$, provided that $q(x)$ is positive and monotonic. This unified representation not only facilitates comparison across different domains but also supports both theoretical understanding and the development of new models.

Below, we introduce the fundamental model components and present a classification scheme for ecological and biological models that can be reformulated in this polynomial quotient format. The main types of ecological and biological models that fit this polynomial quotient framework, together with their typical forms and representative examples, are summarized in Table~\ref{tab:poly-quotient-models}.

\subsection*{Fundamental model components}

The models under consideration are constructed using combinations of the following fundamental terms:

\begin{enumerate}
    \item \textbf{Linear Term (\( a x \))}: Represents density-independent processes such as constant birth or death rates, or exponential growth in population models.
    \item \textbf{Quadratic Term (e.g., \( r x \left(1 - \frac{x}{K} \right) \))}: Captures density-dependent growth limited by environmental carrying capacity, as in logistic growth.
    \item \textbf{Saturating Function \Big( \( \dfrac{a x^n}{H + x^n} \) \Big)}: Models saturation at high concentrations, as seen in enzyme kinetics (Hill function) or functional responses (Holling type II)~\cite{holling1959some}.
    \item \textbf{Nonlinear Inhibition Function \Big( \( \dfrac{a x^m + b x^{m-1} + \dots}{A x^n + B x^{n-1} + \dots} \), with \( n > m \) \Big)}: Represents inhibitory effects, such as repression in gene regulatory networks.
\end{enumerate}

\subsection*{Categorization of models}

We classify models into three main groups based on combinations of these fundamental terms:

\textbf{Group 1: Linear and Saturating Terms}

\textbf{General Form:}

\[
\frac{dx}{dt} = \frac{a x^n}{H + x^n} - c x,
\]
where, \( x \) is concentration or population size, \( a \) is maximum production rate, \( H \) is half-saturation constant, \( n \) is Hill coefficient (reflects cooperativity), and \( c \) is degradation rate constant.

\textbf{Applications or examples:}
Models in this category typically describe situations where a product enhances or catalyzes its own production, but this effect saturates at high concentrations. Examples include gene expression dynamics, where a gene product promotes its own synthesis while simultaneously being degraded in proportion to its concentration~\cite{alon2007introduction}, as well as enzyme kinetics involving autocatalytic reactions that display saturation behavior.

A further example in this class is the quorum-sensing toy model~\cite{DOCKERY200195}, which we also use as a biological case study in the main text (see~\eqref{eq: quorum}). In this model, the Hill type feedback term $Vx^2/(K+x^2)$ captures saturating self amplification of signal production, $x_0$ represents basal synthesis, and the loss term is approximately linear with an effective rate $d(\rho)$ that may depend on cell density $\rho$ (dilution or degradation). The model therefore has the same linear plus saturating structure as the general form above, corresponding to $n=2$ with an added constant source term.

A more complex example is the modeling of T-cell activation in adaptive immunity, where $x$ denotes the concentration of immune effector cells and $H$ represents the activation threshold~\cite{murray2002mathematical}. 

This class of models also describes Lake eutrophication, for instance, Carpenter's model (see Dakos et al.~\cite{dakos2022MainRef}) represents nutrient dynamics in lakes as follows:
\[
\frac{dx}{dt} = a - b x + r \frac{x^k}{x^k + h^k}
\]

Similarly, the Klausmeier model~\cite{Klausmeier1999} for vegetation patterns in semi-arid environments can be transformed into this format~\cite{saade2023MainRef}.

Another combination of these terms appears in models of bacterium-phagocyte dynamics~\cite{Malka_2010_BacteriumModel}, which may be represented in the general form:
\[
\frac{dx}{dt} = \frac{a x}{1 + b x} - \frac{c x}{1 + d x} - \beta x.
\]

\textbf{Group 2: Quadratic and Saturating Terms}

\textbf{General Form:}

\[
\frac{dx}{dt} = r x \left(1 - \frac{x}{K} \right) - \frac{a x^n}{x^n + H},
\]
where \( r \) is intrinsic growth rate, \( K \) is carrying capacity,  \( a \) and \( H \) are parameters of the saturating function, and \( n \) is Hill coefficient.

\textbf{Applications or examples:}

These models combine population growth constrained by environmental carrying capacity with mortality due to predation, which is often described using a saturating functional response~\cite{holling1959some, real1977kinetics}. A prominent example is the Noy-Meir model (herbivory model), which captures grazing systems where plant growth is offset by consumption from herbivores~\cite{noy1975stability}. Saade et al.~\cite{saade2023MainRef} demonstrate how this model can be reformulated in the polynomial quotients format, particularly when $n = 1$, resulting in a cubic polynomial $p(x)$.

This class also includes models for the dynamics of infections such as \textit{Listeria monocytogenes} in the gut, where the impact of the host’s immune response is incorporated into the system~\cite{RAHMAN2016101}.

\textbf{Group 3: Nonlinear inhibition with other terms}

These models combine nonlinear inhibition functions with linear or saturating terms and are suitable for complex phenomena involving inhibition.

\textbf{Applications or examples:}

\textit{Gene Repression with Negative Feedback:}  
These models describe scenarios where a gene product represses its own synthesis by inhibiting transcription~\cite{alon2007introduction}. As the concentration of the gene product ($x$) increases, the production rate decreases due to self-repression. For example, certain metabolites inhibit the transcription of enzymes required for their own synthesis by blocking the DNA-to-mRNA process~\cite{murray2002mathematical}. This dynamic can be represented as:
\[
\frac{dx}{dt} = \frac{\alpha}{1 + (x/K)^n} - \beta x,
\]
where $x$ is the concentration of a repressor protein, $\alpha$ is the maximal production rate, $K$ is the repression coefficient, $n$ is the Hill coefficient (indicating cooperativity), and $\beta$ is the degradation or dilution rate.

\textit{Goodwin Oscillator Model:}  
The above repression model can be extended to account for enzymatic degradation, as in the Goodwin oscillator~\cite{GOODWIN1965} and related cell-cycle models. Here, gene regulation is modeled with both production-inhibition and Michaelis–Menten-like degradation kinetics~\cite{Rombouts2023_EnzymaticDegradation}:
\[
\frac{dx}{dt} = a \frac{K^m}{K^m + x^m} - b \frac{x}{1 + c x} - \beta x.
\]
This form includes linear, saturating, and inhibiting terms. When $m=1$, it can be written as a polynomial quotient, $p(x)/q(x)$.

\textit{Nutrient Uptake with Substrate Inhibition:}  
Enzyme-catalyzed reactions that display substrate inhibition, such as those described by the Haldane-Andrews equation~\cite{cornish2013fundamentals}, are also included in this group:
\[
\frac{dx}{dt} = \frac{V x}{K_m + x + \frac{x^2}{K_s}} - \beta x,
\]
where $x$ is the substrate (nutrient) concentration, $V$ is the maximum uptake rate, $K_m$ is the half-saturation constant, $K_s$ is the inhibition constant, and $\beta$ is the microbial decay or loss rate. The first term describes substrate inhibition, while the second accounts for linear decay or maintenance losses.

\textit{Ecological Example – Desertification Model:}  
An ecological instance is the Rietkerk and Van de Koppel model for desertification, which can be reformulated in this framework~\cite{saade2023MainRef}:
\[
\frac{dx}{dt} = r x \left(1 - \frac{x}{K}\right) \frac{W}{W + B} - M x,
\]
with
\[
W = \frac{I(x + A q)}{(L + U x)(x + A)}.
\]
Further details on this transformation and model interpretation can be found in~\cite{saade2023MainRef, dakos2022MainRef}.

\begin{landscape}
\begin{table}[ht!]
\centering
\caption{Summary of ecological and biological model groups that can be represented as polynomial quotients,
illustrating general forms, key mechanisms, and representative applications.}
\begin{tabular}{p{7cm}|p{5cm}|p{10cm}}
\hline
\textbf{Group} & \textbf{General Form} & \textbf{Representative Examples/Applications} \\
\hline

1. Linear and Saturating Terms
& \[\frac{dx}{dt} = \frac{a x^n}{H + x^n} - c x\]
& \vspace{.1cm}
      
    Gene expression with autocatalysis and degradation \cite{alon2007introduction}\vspace{.5cm}

    T-cell activation in adaptive immunity \cite{murray2002mathematical}\vspace{.5cm}

    Lake eutrophication (Carpenter’s model) \cite{dakos2022MainRef}\vspace{.5cm}

    Klausmeier vegetation patterns \cite{Klausmeier1999, saade2023MainRef}\vspace{.5cm}

    Quorum Sensing in Pseudomonas aeruginosa
    \cite{DOCKERY200195}\vspace{.5cm}
    
    Bacterium-phagocyte dynamics \cite{Malka_2010_BacteriumModel}
\vspace{.1cm}
\\
\hline

2. Quadratic and Saturating Terms
& \[\frac{dx}{dt} = r x \left(1 - \frac{x}{K} \right) - \frac{a x^n}{x^n + H}\] 
& \vspace{.1cm}
 Noy-Meir herbivory/grazing model~\cite{noy1975stability, saade2023MainRef}
\vspace{.5cm}
    
 Infection dynamics with immune response (e.g., \textit{Listeria} in gut) \cite{RAHMAN2016101}

 \vspace{.1cm}
\\
\hline

3. Nonlinear Inhibition with Other Terms
&  
  \[\frac{dx}{dt} = \frac{\alpha}{1 + (x/K)^n} - \beta x\]  

\[\frac{dx}{dt} = a \frac{K^m}{K^m + x^m} - b\frac{x}{1 + c x} - \beta x\] 

\[ \frac{dx}{dt} = \frac{V x}{K_m + x + \frac{x^2}{K_s}} - \beta x\]

& 
\vspace{.1cm}
     Gene repression with negative feedback~\cite{alon2007introduction, murray2002mathematical}
     \vspace{.5cm} 
     
     Goodwin oscillator (gene regulation)~\cite{GOODWIN1965, Rombouts2023_EnzymaticDegradation} \vspace{.5cm} 
     
     Enzyme kinetics with substrate inhibition (Haldane-Andrews)~\cite{cornish2013fundamentals} \vspace{.5cm} 
     
     Desertification model (Rietkerk and Van de Koppel)~\cite{saade2023MainRef, dakos2022MainRef}

\end{tabular}
\label{tab:poly-quotient-models}
\end{table}
\end{landscape}

\section*{S2: Concept of memory and potential energy}\label{Sec: Sup. memory}
\begin{figure}[ht!]
    \centering
    \includegraphics[width=\textwidth]{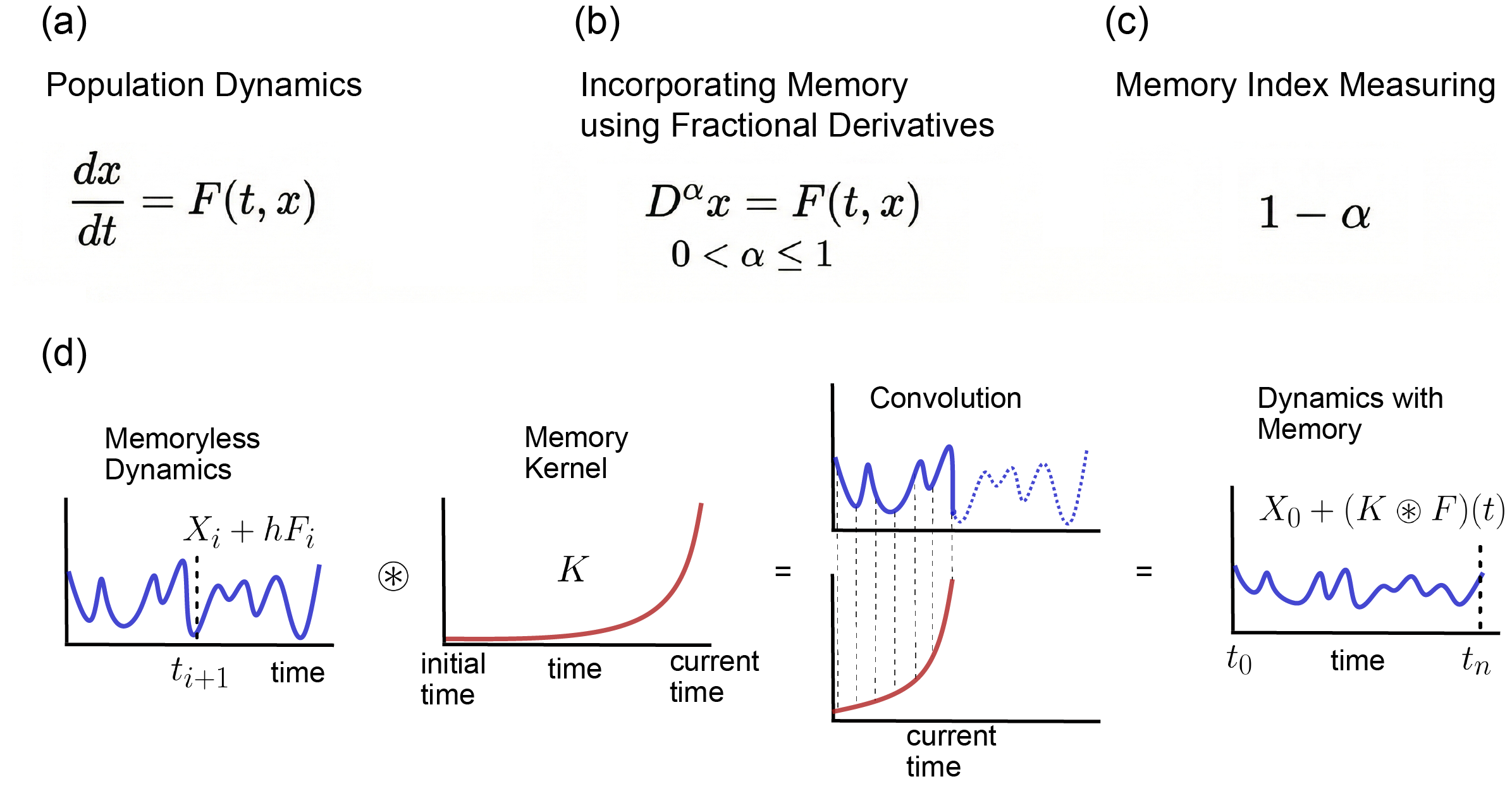}
    \caption{(a) Standard population dynamics described by a differential equation. (b) Extension to a fractional differential equation with order $0<\alpha\leq 1$, incorporating memory effects. (c) The memory index, defined as $1-\alpha$, quantifies the degree of memory present. (d) Memory is described by the convolution of the system function $F$ with a power-law kernel $K$, illustrating the integral kernel modifies the population dynamics $X$.} 
    \label{fig: memory_modeling}
\end{figure}

This section introduces the concept of memory employed in our study, along with the relevant theoretical and mathematical foundations. Figure~\ref{fig: memory_modeling} illustrates the overview of the concept that we used.

\subsection*{Incorporating memory into models using fractional derivatives}
Memory effects in dynamic systems can be mathematically represented using fractional derivatives of order $\alpha$~\cite{Khalighi2022Ploscb, khalighi2025impact}:

\begin{equation}\label{eq:fracModel}
\mathcal{D}^{\alpha}X(t)=F(t,X(t)),
\end{equation}

where $\mathcal{D}^{\alpha}$ denotes Caputo fractional derivative operator~\cite{kilbas2006theory, podlubny1998fractional}:
\begin{equation}\label{eq:caputo}
\mathcal{D}^{\alpha} X(t) = I^{1-\alpha}_{t_0}X'(t) = \frac{1}{\Gamma(1-\alpha)}\int_{t_0}^{t} \frac{X'(\tau) \, d\tau}{(t-\tau)^{\alpha}}.
\end{equation}

This fractional derivative involves the Riemann–Liouville fractional integral, defined by:

\begin{equation}\label{eq:RLInt}
I^{\alpha}_{t_0}X(t) = \frac{1}{\Gamma(\alpha)}\int_{t_0}^{t} \frac{X(\tau)}{(t-\tau)^{1-\alpha}}d\tau, \quad \alpha \in \mathbb{R}^+.
\end{equation}

In Equation \eqref{eq:fracModel}, the fractional derivative integrates contributions from past states of the derivative $X'(\tau)$, weighted by the power-law term $(t - \tau)^{-\alpha}$. The gamma function $\Gamma(\alpha)$ serves as a normalization factor. This structure demonstrates the non-local characteristics of fractional operators, explicitly incorporating historical memory by summing contributions from all previous times $\tau$.

Equivalently, this relationship can be represented using a Volterra integral equation \cite{kilbas2006theory}:

\begin{equation}\label{eq:transformation}
X(t_n) = X(t_0) + \int_{t_0}^{t_n} K_{\alpha}(t-\tau) F(\tau, X(\tau)) , d\tau,
\end{equation}
where the memory kernel $K_{\alpha}(t-\tau)$ is defined as:
\[
K_{\alpha}(t-\tau) = \frac{(t-\tau)^{\alpha - 1}}{\Gamma(\alpha)},
\]
and $t_0$ is the initial time and $t_n$ is any later time. In practice, we often evaluate this relation on a discrete grid $\{t_i\}_{i=0}^n$, so $t_n$ can be taken as the $n$th sampled time point.

Unlike classical derivatives, Caputo fractional derivatives incorporate the historical states of a function. The convolution with a power-law kernel, $K_{\alpha}(t - \tau)$, ensures that every previous state $\tau$ influences the current derivative at time $t$. According to Ogle’s classification \cite{ogle_Eco_memory2015}, such kernels represent long-term memory: they lack explicit time lags, yet exhibit slowly decaying influence over time. In this formulation, the singular kernel models systems where memory effects gradually diminish. Specifically, there exists a parameter $\alpha > 0$ such that $\lim_{t \to \infty} t^{-\alpha} K_{\alpha}(t-\tau)$ remains finite for fixed $\tau$. This property makes the kernel particularly suitable for capturing real-world dynamics with memory effects~\cite{Tarasova2018, METZLER20001}.

The fractional order $\alpha$ adjusts the intensity of memory: smaller $\alpha$ values correspond to stronger and longer-lasting memory, while $\alpha = 1$ reduces the derivative to the traditional, memory-free (local) derivative. Thus, we quantify memory strength by $1 - \alpha$, where 0 indicates no memory and values closer to 1 indicate stronger memory. The derivative order $\alpha$ controls both the memory duration and its decay rate.

\subsection*{Sensitivity to initial conditions in memory-driven dynamics}
A key challenge when incorporating memory into dynamical systems is the heightened sensitivity to initial conditions. In a memory-free model, the reconstructed stability landscape is independent of where the system begins; trajectories and the resulting landscape overlay are identical no matter the starting point. By contrast, memory ``remembers'' any past trajectory, so the same model parameters can yield different landscapes depending on the initial state.

\begin{figure}[ht!]
    \centering
    \includegraphics[width=\textwidth]{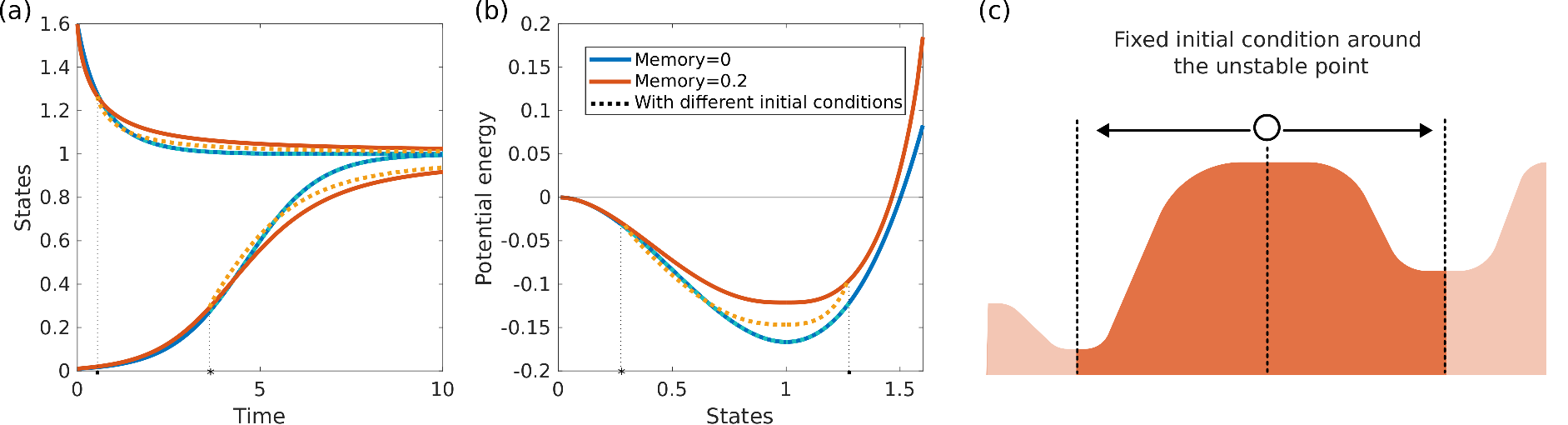}
\caption{
(a) Dynamics of a simple logistic model, \(\displaystyle\frac{dX}{dt} = X(1-X)\), with and without memory, shown for various initial conditions.
(b) The corresponding potential energy landscapes from panel (a), highlighting how the presence of memory leads to distinct patterns depending on the initial state.
(c) Illustration of how initial conditions are addressed in our study. For potential energy reconstruction, we use fixed initial values near the unstable point and limit our analysis to the region between the two stable states, focusing on the interval surrounding the unstable point up to stable points.
}
    \label{fig: initial conditions}
\end{figure}

We illustrate this in Figure~\ref{fig: initial conditions} with a logistic growth example model: without memory, every run collapses onto a single trajectory and landscape; with memory, different starting values produce distinct trajectory shapes and slightly altered landscapes. To control for this in our study, we always initialize trajectories just around the unstable equilibrium (the hilltop) when reconstructing the basins (Figure~\ref{fig: initial conditions}c). This consistent choice ensures that all landscape features—and comparisons across models—are based on the same portion of state space. 

\subsection*{Potential energy for systems with memory effects}

To evaluate the stability landscape of a system incorporating memory effects, we use the model defined in equation~\eqref{eq:fracModel}. The potential energy can be determined by solving the following integral equation:
\begin{equation}\label{eq: model with kernel}
V = -\int_{X_a}^{X_b} \dot Z(t)\, dZ,
\end{equation}
where $Z(t)$ is a time-dependent function that incorporates the influence of past states through the memory kernel $K$. In a discretized form~\eqref{eq:transformation}, using the index $i$ to denote the time step, $Z(t)$ can be written as:
\begin{equation}
Z_i = Z_0 + \int_{t_0}^{t_i} K(t_i - \tau) F(\tau, Z(\tau))\, d\tau, \quad X_a \leq Z_i \leq X_b.
\end{equation}

The parameters $X_a$ and $X_b$ define the boundaries of the state space considered within the stability landscape, chosen specifically to encompass both valleys and the peak of a bistable system.

\subsection*{Potential energy after perturbation and dynamics stability landscapes}

In this study, we demonstrate that stability landscapes change in response to pulse perturbations. Our main focus is on abrupt changes—such as sudden shifts in light, temperature, chemical concentrations, or biotic factors—rather than gradual transitions. To model such events, we introduce a sudden stress into the simulation. In memory-free models, the response to perturbation is immediate; the stability landscape transforms abruptly if we assume that the system's conditions remain unchanged before and after the disturbance.

However, in systems with memory, the situation is fundamentally different. After a perturbation, the stability landscape does not shift instantaneously but instead recovers gradually, often retaining traces of the disturbance for some time. This lingering effect leads to slower recovery or even irreversible changes in the stability landscape. As a result, the post-perturbation landscape becomes qualitatively distinct from the original configuration. Specifically, the valleys and peaks of the landscape evolve after a perturbation, and time is required for the system to approach its previous state. Thus, the depths and curvatures of potential wells may also assume new shapes.

\begin{figure}[ht!]
    \centering
    \includegraphics[width=1\textwidth]{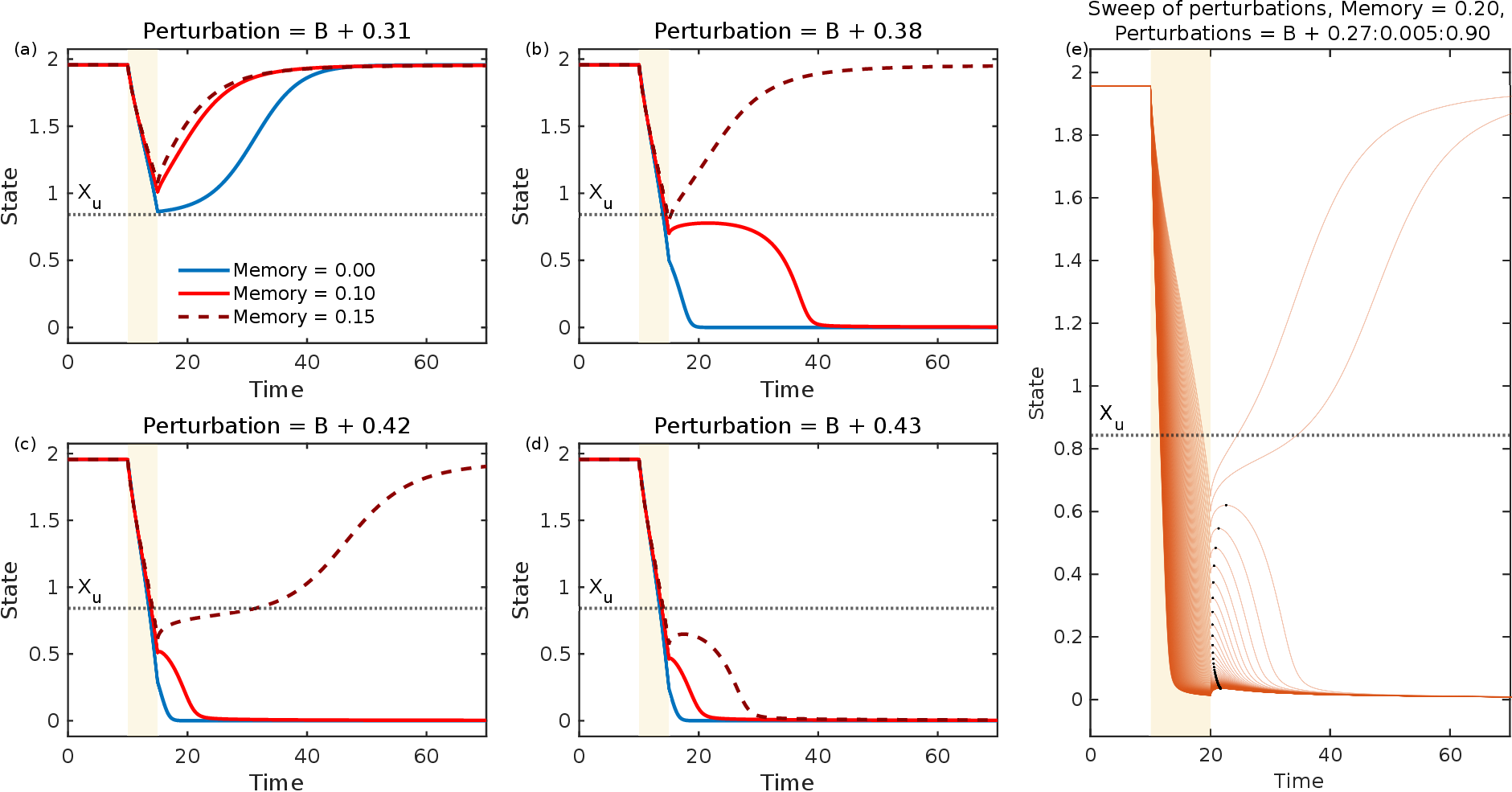}
    \caption{ (a)-(d) Dynamics for different perturbation strengths and memory strength levels. (e)  Dynamics of the system under varying perturbation strengths ($0.27$ to $0.9$) with a fixed memory strength of $0.2$. Black dots indicate the points where the system's state changes direction and begins converging toward the zero stable state.}
    \label{fig: dyn_pert1}
\end{figure}

\begin{figure}[ht!]
    \centering
    \includegraphics[width=1\textwidth]{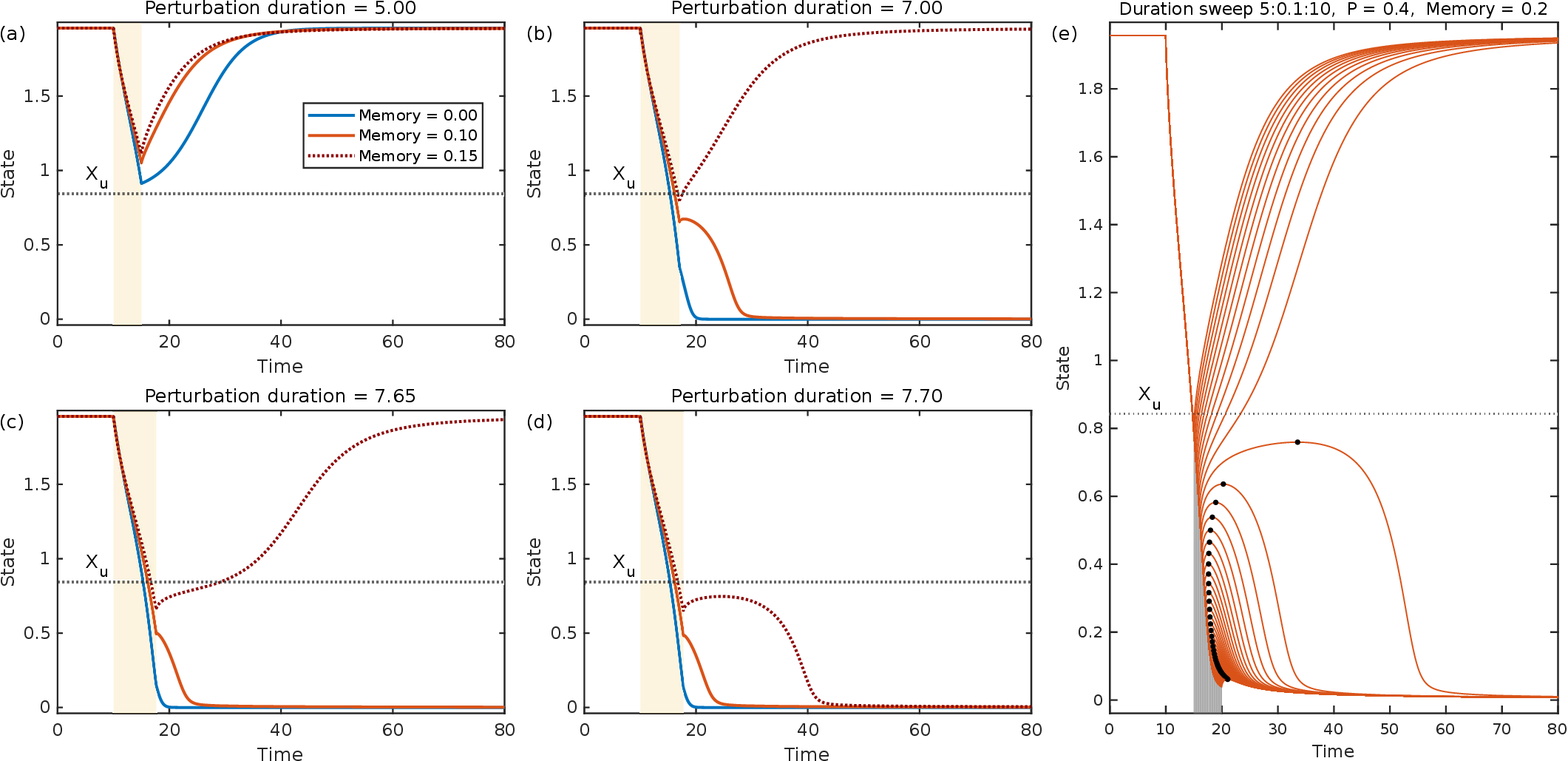}
    \caption{(a)-(d) Dynamics for different perturbation lengths and memory strength levels. (e) Dynamics of the system under varying perturbation periods ($5$ to $10$) with a fixed memory strength of $0.2$. Black dots indicate the points where the system's state changes direction and begins converging toward the zero stable state.}
    \label{fig: dyn_pert2}
\end{figure}

To see what happens in simulations, consider a system, e.g., herbivory model~\eqref{eq: herbivor}, initially at a positive stable state subjected to a pulse perturbation, here, by changing parameter $B$, that displaces the system state. Figures~\ref{fig: dyn_pert1}a-d and \ref{fig: dyn_pert2}a-d illustrate various scenarios involving different perturbation strengths, durations, and memory intensities. Notably, these figures reveal situations where, during perturbations, states with memory cross the unstable point of the pre-perturbation landscape; $X_U$. After perturbation, the states sometimes return to the unstable points—occasionally crossing them, or at other times approaching but ultimately retreating. Such behavior would be impossible without memory.


Figures~\ref{fig: dyn_pert1}e and \ref{fig: dyn_pert2}e further consolidate these findings by comparing different perturbation strengths and levels for a fixed memory setting. These plots demonstrate cases where the system's state crosses the unstable point during perturbation. Remarkably, the system with memory sometimes manages to return to the original basin of attraction despite initially crossing the unstable point. This ability to recover is especially intriguing because, in memory-free systems, once conditions revert post-perturbation, states that surpass the unstable point are typically drawn towards zero.

An additional intriguing observation arises when considering the recovery process in systems with memory. As the state attempts to return to the original basin bottom, the evolving landscape occasionally changes more rapidly than the state's movement towards stability. Consequently, the state trajectory can reverse direction midway through recovery and move toward zero instead. This change of direction is marked explicitly with black dots on the curves in Figures~\ref{fig: dyn_pert1}e and \ref{fig: dyn_pert2}e, clearly illustrating the dynamic nature of the post-perturbation stability landscape, gradually returning to its original configuration rather than instantaneously resetting.

To further clarify this phenomenon by illustrating the potential landscapes, we present examples such as those in Figure~\ref{fig: land_herb_RS}, which depict trajectories on stability landscapes before, during, and after perturbations under the conditions used in Figure~\ref{fig: mem_herb_RS}e-h. In systems with memory, the trajectory of the ball appears to exhibit unusual behaviors. For instance, the ball seems to surpass the hill (unstable point) during the perturbation and then, following the perturbation, appears to move upward towards the hill again. Such a scenario would be impossible in memory-free systems. The explanation lies in the dynamic nature of the landscape itself; rather than the ball moving upward, the hill is shifting dynamically, gradually returning to its original configuration. The interplay between the velocity of the ball's movement and the evolution speed of the stability landscape (basin of attraction) dictates the trajectory and final resting point of the state. Supplementary Videos~S2 and~S3 show these scenarios.

\subsubsection*{Mathematical characterization of post-perturbation landscape dynamics}
Consider a system initially at the bottom of the valley, indicating positive equilibrium, denoted as $X_0 = X_{S_2}$. Suppose the system has the potential to recover after a disturbance. During the perturbation, a perturbed function $F_p$ causes the system to lose stability at $X_{S_2}$ and transition toward $X_{U'}$, near the unstable state $X_U$. The potential energy under this perturbation is given by:
\begin{equation}\label{Eq: Vp}
V_p = -\int_{X_{U'}}^{X_{S_2}} \dot Z(t)\, dZ,
\end{equation}
with each $Z_i \in [X_{U'}, X_{S_2}]$ defined as:
\begin{equation}\label{Eq: Vp1}
Z_i = Z_0 + \int_{t_{S_2}}^{t_i} K(t_i-\tau) F_p(\tau, Z(\tau))\, d\tau.
\end{equation}

Here, $t_{S_2}$ is the time when the system is at state $X_{S_2}$, and $t_i$ is the current time as the system transitions from $X_{S_2}$ to $X_{U'}$.

After the perturbation, assume the system settles at a new state $X_{S_2'}$, close to the original stable state $X_{S_2}$. The potential energy in this new state is given by:
\begin{equation}\label{Eq: V2'}
V = -\int_{X_{U'}}^{X_{S_2'}} \dot Z(t)\, dZ,
\end{equation}
where each $Z_j \in [X_{U'}, X_{S_2'}]$ is defined as:
\begin{equation}\label{eq: funPert}
Z_j = Z_0 + \int_{t_{S_2}}^{t_j} K(t_j-\tau) G(\tau, Z(\tau))\, d\tau,
\end{equation}
with $t_j$ as the current time during the transition from $X_{U'}$ to $X_{S_2'}$, and $G$ as a piecewise function:
\begin{equation}
G(\tau, Z(\tau)) = 
\begin{cases}
    F_p(\tau, Z(\tau)), & \text{if } \tau \in [t_{S_2}, t_{U'}], \\
    F(\tau, Z(\tau)),   & \text{if } \tau \in [t_{U'}, t_{S_2'}],
\end{cases}
\end{equation}
where $t_{U'}$ and $t_{S_2'}$ are the times at which the system reaches $X_{U'}$ and $X_{S_2'}$, respectively.

Equation~\eqref{Eq: V2'} shows that, after a perturbation, the system's potential valley cannot fully return to its original shape due to the influence of the memory kernel $K$ and the piecewise function $G$.

In contrast, for a memory-free system, the post-perturbation potential energy~\eqref{Eq: V2'} depends only on $F(\tau, Z(\tau))$:
\begin{equation}\label{eq: poten_noM}
V = -\int_{X_{U'}}^{X_{S_2'}} \dot Z(t)\, dZ,
\end{equation}
with
\begin{equation}
Z_j = Z_0 + \int_{t_{U'}}^{t_j} F(\tau, Z(\tau))\, d\tau,
\end{equation}
where the system’s states are independent of the past, and the integral starts at $t_{U'}$, the time at which the system reaches $X_{U'}$. Thus, only the immediate state following the perturbation depends on the perturbed dynamics ($\frac{dX}{dt} = F_p$), while subsequent states evolve according to the original function $F$. As the potential values derived from Equation~\eqref{eq: poten_noM} match those before the perturbation, the potential landscape remains unchanged and unaffected by the pulse disturbance.

These comparisons highlight that memory not only shapes transient responses but can also lead to qualitatively different stability landscapes following perturbations.

\subsection*{Numerical implementation}

Here, we describe how to numerically evaluate the potential energy of systems with memory. Instead of seeking an exact analytical solution, we compute the potential energy numerically by tracing the path of the potential valleys as the ball (system states) moves through them.

The first step is to determine the equilibrium points of the model~\eqref{eq:fracModel}. For clarity, let us denote $X_{S_1}$ and $X_{S_2}$ as the stable states and $X_U$ as the unstable state, with $X_{S_1} < X_U < X_{S_2}$. Our goal is to construct a potential energy profile across the state space, denoted as $X_t = (X_{S_1}, X_{S_2})$.

To achieve this, we numerically solve Equation~\eqref{eq:fracModel} using initial conditions slightly greater or less than $X_U$. This allows us to compute the potential energy along both sides of the valleys, specifically within intervals such as $(X_U, X_{S_1})$ and $(X_U, X_{S_2})$, which correspond to regions near the unstable point. For $X > X_{S_2}$ (the region beyond the right valley), the potential energy may differ in memory systems depending on the initial condition $X_a > X_{S_2}$. While we do not analyze this case in detail here, we include results for these regions to visualize the full potential landscape.

For evaluating Equation~\eqref{eq: model with kernel}, instead of directly integrating $F$, we numerically solve the fractional equation~\eqref{eq:fracModel} (using tools such as \texttt{FdeSolver}~\cite{khalighi2024algorithm} in Julia or equivalent code in MATLAB) to obtain the states $X$ over the desired domain. Once we have $X$, we compute its derivatives. Thus, in Equation~\eqref{eq: model with kernel}, we substitute $F$ with the numerical derivative $dX/dt$, and then perform the integration numerically.

To improve accuracy, we use a fourth-order central difference scheme ($O(h^4)$) to approximate the first derivative at each point:
\begin{equation}
    \frac{dX_t}{dt} = \frac{X_t(n-2) - 8X_t(n-1) + 8X_t(n+1) - X_t(n+2)}{12h}
\end{equation}
After computation, the results are sorted in ascending order of the state values. These sorted values are partitioned into a predefined grid of $2(X_{S_2}-X_{S_1})/h$ points. To ensure uniformity and improve computational precision, we use Cubic Hermite interpolation to generate evenly spaced states and their corresponding derivatives~\cite{cubic}. Finally, we apply numerical integration methods such as Simpson’s rule or the trapezoidal rule, using an optimized step size, to integrate the derivatives with respect to the states.

\subsection*{Anomalous diffusion in population dynamics motivates fractional order models}

This section motivates fractional calculus as a principled framework for modeling memory by linking environmental heterogeneity, trapping, and subdiffusion to fractional order population dynamics.

In ecological and biological systems, population dynamics and spatial spread depend significantly on both intrinsic demographic processes and environmental spatial structure. Classical population models, such as logistic growth, typically assume homogeneous environments and instantaneous responses to local population densities. However, real-world environments often feature physical heterogeneity, including obstacles or barriers that influence movement patterns and delay habitat colonization.

For instance, in soils and porous geological structures~\cite{rajyaguru2024diffusion, datta2025bacterial}, complex pore networks trap microorganisms or particles, resulting in sublinear growth of the mean squared displacement over time. Similarly, within biological cells, the densely crowded cytoplasm restricts molecular movement, causing subdiffusive behavior~\cite{WEISS20043518, yu2018subdiffusion}. Comparable effects occur in polymeric and viscoelastic materials, where internal architectures hinder particle dynamics, resulting in memory effects~\cite{lim2025anomalous, joo2020anomalous, kim2025short}. Collectively, these phenomena lead to \textit{anomalous diffusion}~\cite{METZLER20001}, wherein the movements of individuals deviate from classical Brownian motion, reflecting underlying memory-driven dynamics.

To illustrate anomalous diffusion, we implemented an agent-based, spatially explicit simulation of logistic population growth on a two-dimensional lattice (see Supplementary video S1 and Fig.~\ref{fig:compare_obs_noobs}a). Agents occupy discrete grid cells and can move to adjacent empty cells (up, down, left, right), mimicking local dispersal. Reproduction occurs probabilistically according to the logistic model:
\begin{equation}
p\_{\text{birth}} = r \left(1 - \frac{X}{K}\right),
\end{equation}
where $X$ represents the current population size, $K$ is the carrying capacity, and $r$ is the intrinsic growth rate. Reproduction only occurs if adjacent space is available, thus incorporating local resource competition.

We analyzed two scenarios: (1) a homogeneous environment without obstacles, and (2) a heterogeneous environment with obstacles resembling a porous medium. In the second scenario, agents near obstacles experience temporary immobility, leading to \textit{trapping} and prolonged waiting times before movement or reproduction can occur. These spatial constraints slow down population expansion, delaying the attainment of carrying capacity compared to obstacle-free environments.

Such delays exemplify \textit{subdiffusive behavior}, characterized by the mean squared displacement increasing more slowly than linearly with time. The mean squared displacement (MSD), denoted $\langle x^2(t) \rangle$, quantifies how far, on average, individuals or particles have moved from their initial positions after a given time $t$. Specifically, it is calculated by following each agent, squaring the distance it has traveled from its starting point, and then averaging over all agents or simulation runs. In classical diffusion (such as Brownian motion), the MSD grows proportionally with time, but in subdiffusive systems, environmental obstacles and crowding slow down movement, resulting in slower MSD growth. Mathematically, this behavior follows a power-law relationship~\cite{METZLER20001}:
\begin{equation}
\langle x^2(t) \rangle \propto t^{\alpha}, \quad 0 < \alpha < 1,
\end{equation}
where the exponent $\alpha$ describes how strongly the environment impedes movement; lower values of $\alpha$ correspond to more pronounced slowing. This contrasts with the linear relationship ($\alpha = 1$) typical of classical diffusion. Standard integer-order differential equations fail to capture these dynamics because they inherently lack memory.

To describe the observed memory effects and anomalous diffusion, we can utilize \textit{fractional-order differential equations}, generalizing the time derivative to a fractional order ($\alpha$). The fractional logistic equation:
\begin{equation}
\mathcal{D}^{\alpha}X = rX\left(1 - \frac{X}{K}\right),
\end{equation}
naturally integrates historical system states. Thus, population growth rates at a given time depend not only on the current state but also on the cumulative history of the system, weighted by a power-law memory kernel.

Simulation results confirm that environments with spatial obstacles generate population dynamics qualitatively matching fractional models. Specifically, growth rates are slower, and equilibrium states are approached gradually. This underscores fractional calculus as a biologically meaningful framework for modeling population dynamics in heterogeneous environments. Memory, emerging as an environmental property, acts as a hidden yet influential parameter shaping population growth, observable indirectly through memory-induced dynamics.
To connect this coarse-grained fractional description to a direct diagnostic of history dependence, we next test whether the same population size $N$ can lead to different futures depending on hidden spatial configuration and trapping history, using a distributional restart experiment.

\subsubsection*{Restart test exposes history dependence from spatial trapping}

We asked whether the forward population dynamics are determined only by the coarse state \(X\) (here, the population size \(N\)) or whether they also depend on hidden microstructure in the agent configuration. We performed a distributional restart test at step \(s=35\), launching two ensembles with fresh random seeds: (i) a \emph{snapshot ensemble} that continues from the exact microstate at step \(35\) (agent positions preserved and, in the obstacle case, trap timers and trapped status preserved), and (ii) a \emph{macro matched ensemble} that preserves only the population size \(N\) while randomizing agent positions (and, in the obstacle case, resetting trap timers) under the same obstacle mask. For each ensemble we ran \(R=300\) replicates, summarized the forward trajectories by the mean and the central \(95\%\) band, and compared the full distributions at selected horizons using the two-sample Kolmogorov-Smirnov test \cite{massey1951kolmogorov}.

\begin{figure}[ht!]
  \centering  \includegraphics[width=\linewidth]{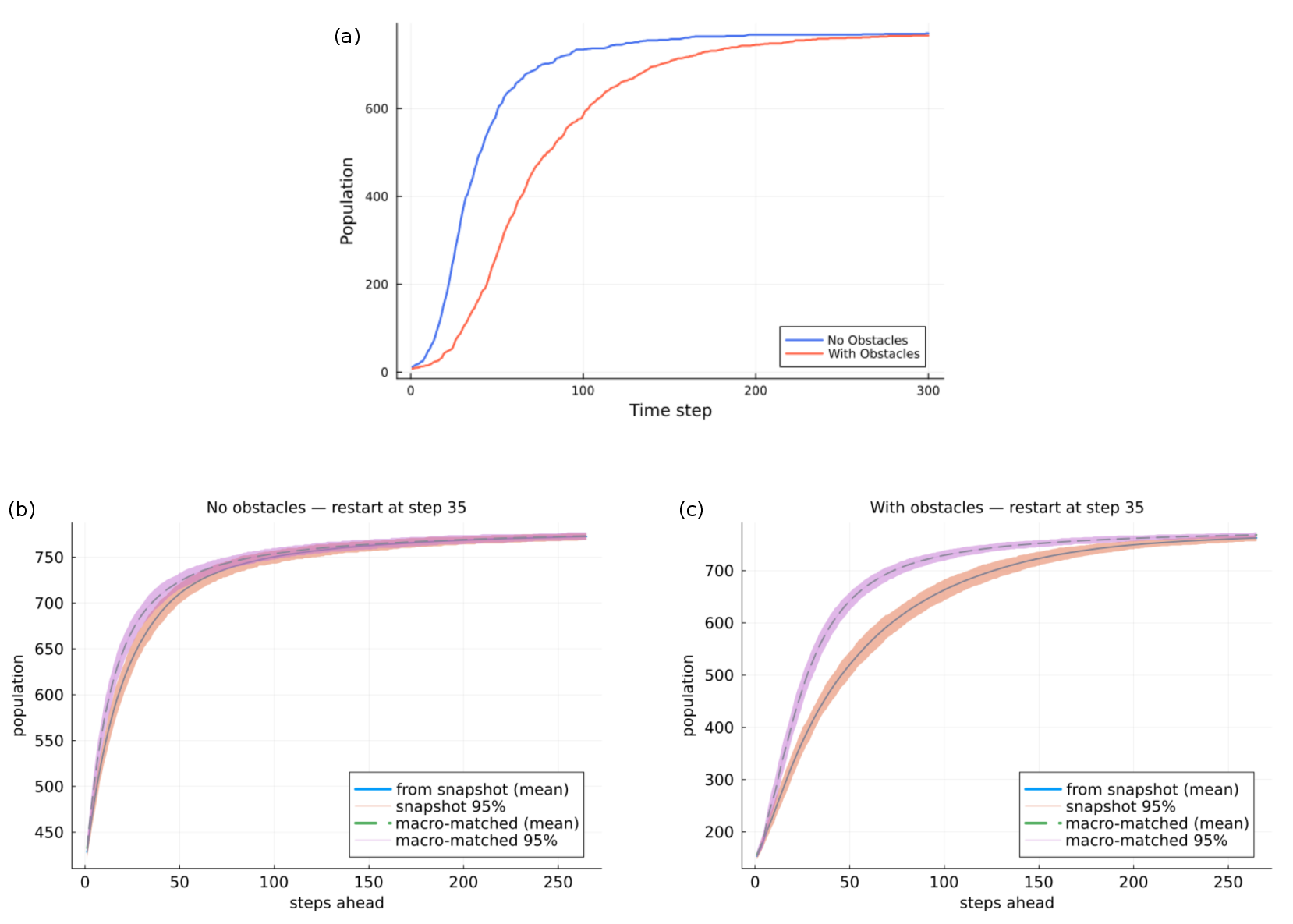}
  \caption{Population dynamics and evidence for history dependence using a distributional restart at step \(35\).
(a) Single realization comparing no obstacles and obstacles with trapping, showing that obstacles slow the approach to carrying capacity.
(b) No obstacle restart at step \(35\): snapshot and macro matched ensembles converge as the horizon increases, consistent with rapid mixing and near sufficiency of population state \(N\) at long horizons.
(c) Obstacle restart at step \(35\): snapshot and macro matched ensembles remain separated, indicating persistent dependence on microstate and trap history beyond \(N\) alone. Shaded regions show the central \(95\%\) band across \(R=300\) replicates.}
  \label{fig:compare_obs_noobs}
\end{figure}

In the baseline trajectories, obstacles mainly delay the transient growth while the long-run population approaches a similar level, motivating a distributional test that separates transient microstate effects from coarse state dependence. We repeated the restart experiment for multiple restart times (e.g., \(s \in \{25,30,35,50,70\}\)) and obtained qualitatively consistent results; therefore, we report only the representative case \(s=35\) in Fig.~\ref{fig:compare_obs_noobs}.

At a forecast horizon \(h\), let \(X_h\) and \(Y_h\) be the samples of population size from the snapshot and macro matched ensembles. The Kolmogorov-Smirnov statistic \(D \in [0,1]\) is the maximum gap between the two empirical cumulative distribution functions. Larger \(D\) indicates stronger separation, with \(D \approx 1\) meaning that the samples are almost completely ordered relative to one another. The \(p\) value is the probability of observing a gap at least as large as \(D\) if both samples were drawn from the same distribution.

\paragraph{No obstacles:}
The restart test indicates that microstructure effects are transient and wash out over time. At short horizons the distributions are strongly different (for \(s=35\): \(D=0.93\) at \(h=10\), \(D=0.96\) at \(h=25\), \(D=0.82\) at \(h=50\), all with \(p \approx 0\)). However, by the final horizon, the ensembles nearly coincide (\(D=0.11\), \(p=0.053\) at \(h=265\)). This is consistent with rapid mixing in a homogeneous environment: once enough time passes, conditioning on \(N\) becomes almost sufficient for predicting the future population distribution, and the detailed spatial arrangement retains little predictive power. This behavior is visible in Fig.~\ref{fig:compare_obs_noobs}b, where the two ribbons approach one another as the forecast horizon increases.

\paragraph{With obstacles:}
In contrast, obstacles induce strong and persistent microstate dependence. At short horizons we observe near complete separation (\(D=0.86\) at \(h=10\), \(D=1.0\) at \(h=25\), \(D=1.0\) at \(h=50\), all with \(p \approx 0\)), and the difference remains large even at the final horizon (\(D=0.7\), \(p \approx 0\) at \(h=265\)). The ribbon plot in Fig.~\ref{fig:compare_obs_noobs}c shows the same effect: macro matched restarts typically recover faster and reach higher population levels earlier than snapshot restarts. This indicates that, in heterogeneous media, the forward distribution depends on microstate variables not captured by \(N\), such as spatial clustering relative to obstacles and the internal trap history. In other words, the environment stores effective memory through immobilization and spatial structure, so the process is not effectively Markov in \(N\) alone. Because the macro matched restarts also reset trap timers and trapped status, part of the faster recovery reflects the removal of internal trapping memory in addition to spatial randomization.

\section*{S3: Impact of memory on dynamical systems under endogenous perturbations}\label{sec: sup3}

In the following sections, we present an analysis of how memory influences the geometry of stability landscapes and key ecological properties, resilience, and resistance, in one-dimensional bistable population models. We first examine how incorporating memory reshapes potential landscapes and alters the system’s approach to equilibrium. We then quantify the consequences of these geometric changes for resilience and resistance, both in representative case studies and across large ensembles of randomly parameterized models. By introducing a novel basin sharpness index, we further dissect how memory’s effects vary depending on the underlying landscape structure. Finally, we synthesize the interrelationships among all major variables to provide an integrated view of the mechanisms through which memory modulates stability in dynamical systems.

\subsection*{Quantitative metrics}\label{sec: quantitative metrics}

\subsubsection*{Depth}
Potential depth is defined as the difference in potential energy between the bottom of the stable valley (the positive stable state) and the peak of the hill separating the two stable states (the unstable equilibrium). This measure reflects how deeply the system is trapped in a stable state and thus how resistant it is to transitions.

\subsubsection*{Flatness and curvature}

To characterize the shape of the potential landscape around stable equilibria, we consider curvature and its inverse, which we define as \textit{flatness}. While curvature is mathematically meaningful, it takes negative values in our context, making comparisons less intuitive. Therefore, we use flatness as a positive, more interpretable alternative.

Curvature quantifies how sharply or gently the potential landscape bends around an equilibrium. It is defined as the second derivative of the potential function and, for one-dimensional systems, is given by~\cite{dakos2022MainRef}:

\[
\text{Curvature} = \left.\frac{d^2V}{dx^2}\right|_{x = x^*} = -F'(x^*) = -\lambda,
\]

where \( x^* \) is the equilibrium point and \( \lambda \) is the eigenvalue of the linearized system at that point. A larger magnitude of curvature implies steeper potential wells (see Figure~\ref{fig S1}d).

To facilitate comparison across models, we define flatness as:

\begin{equation}
\text{Flatness}(i) = \lceil \text{max}(|\text{curvature}(i)|) \rceil - |\text{curvature}(i)|,
\end{equation}

where \( \text{curvature}(i) \) is the curvature of the \( i \)th model, and \( \text{max}(|{curvarture}|) \) is the maximum curvature magnitude among all models. This results in a positive flatness value that increases with the flatness of the potential well.

The index \( i \) can refer to either a model with a specific memory level or a sample model in a memory-free setting. Flatness offers a clear, interpretable metric for comparing potential shapes across varying memory intensities or structural parameters.

\subsubsection*{Basin sharpness index}
In the following, we show that, across our case studies, potential depth and basin flatness (curvature magnitude) are often strongly correlated. We therefore introduce a compact metric that summarizes both features in a single number. The basin sharpness index combines normalized depth and normalized curvature magnitude as:

\begin{equation}
\text{Basin sharpness index} = \frac{\text{normalized(depth)} + \text{normalized}(|\text{curvature}|)}{2}.
\end{equation}

This index ranges from 0 to 1. A value close to 0 indicates a shallow and flat basin, while a value near 1 reflects a deep and steep basin.

\subsection*{Impact of memory on stability landscape}
We focus on one‐dimensional, bistable population dynamics models, such as the herbivory example. Each model has two stable equilibria: one at zero (extinction) and one at a positive, biologically meaningful density. Our analysis centers on the positive basin (the right‐hand valley) and the intervening hill separating it from the extinction state.

\begin{figure}[ht!]
    \centering
    \includegraphics[width=1\textwidth]{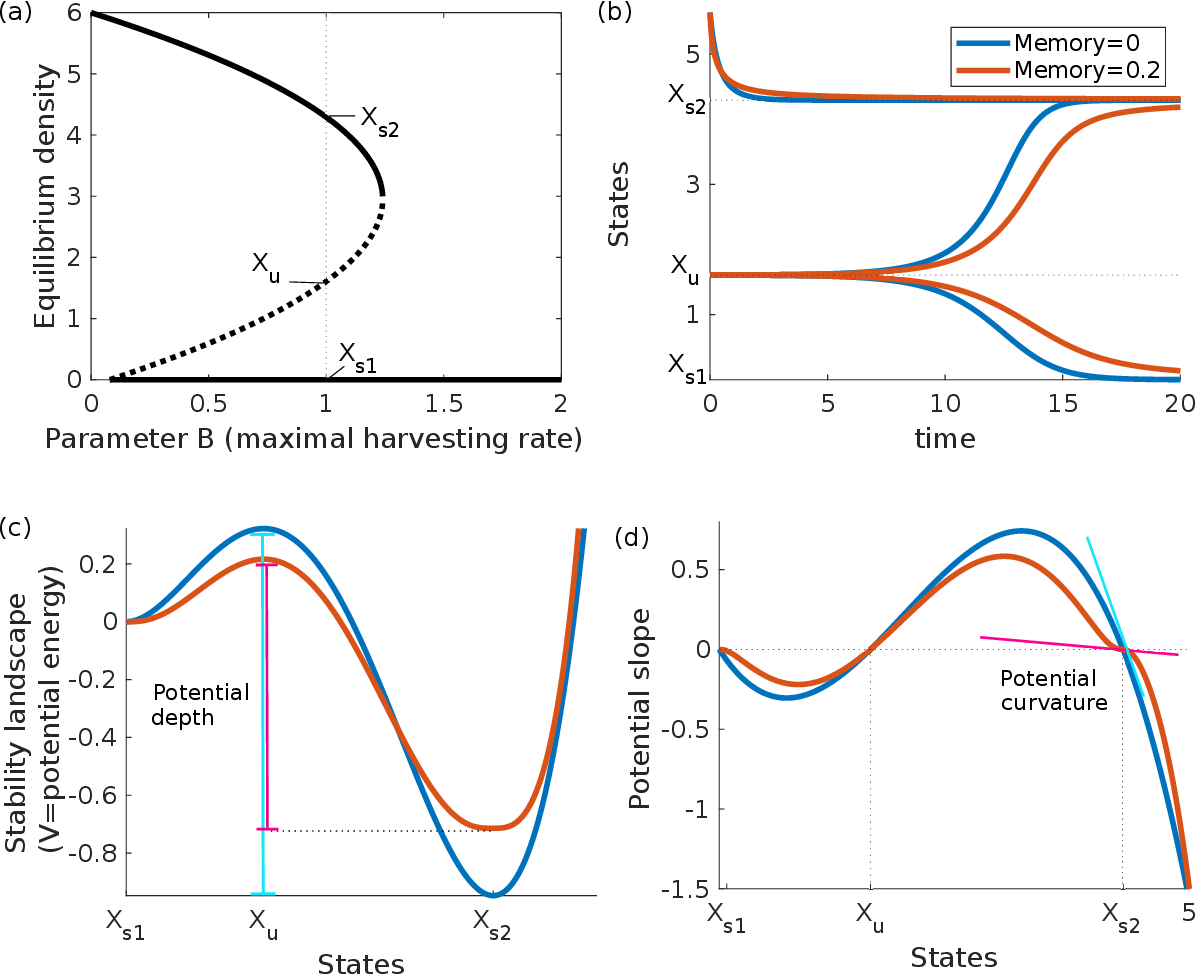}
    \caption{
    (a) The bifurcation diagram shows the equilibria of the herbivory model as a function of parameter \(B\). At \(B=1\), two stable states (\(X_{S1}\) and \(X_{S2}\)) and one unstable state (\(X_U\)) coexist. 
    (b) The dynamics of the bistable model, with (red curves) and without (blue curves) memory, are illustrated for three different initial abundances, with \(B=1\). 
    (c) The potential energy landscape, reconstructed from the trajectories in panel (b). Here, potential depth is defined as the distance from the bottom of the valley (positive stable state $X_{S2}$) to the top of the hill (unstable state).
    (d) The potential slope represents the growth rate of the system, indicating how quickly the dynamics evolve within the state space. Stable (\(X_{S1}\), \(X_{S2}\)) and unstable (\(X_U\)) points correspond to the roots of the potential slope curve. The potential curvature (red and blue lines) shows the steepness of the slope at the stable points.
    }
    \label{fig S1}
\end{figure}

When memory is introduced into a model, it acts like a drag on system motion: trajectories take longer to approach equilibrium. Because we reconstruct potential landscapes directly from those trajectories, slower convergence appears as gentler slopes, both at the basin floor and on the surrounding hill. Figure~\ref{fig S1} shows how, under memory, population trajectories of the herbivory model converge more gradually toward the positive attractor, and the corresponding landscape exhibits an altered hill and a noticeably flatter basin bottom.

To generalize beyond a single example, we generated 1000 ensembles of randomly parameterized, bistable polynomial models. For each positive basin, we measured two geometric features: 1) {Depth}, and 2) {Flatness}.

\begin{figure}[ht!]
    \centering
    \includegraphics[width=1\textwidth]{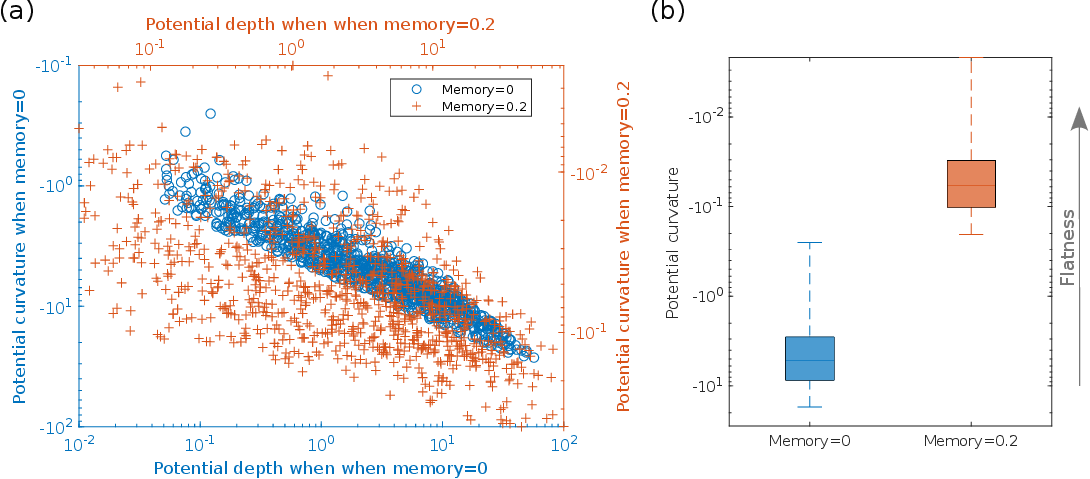}
    \caption{
    Illustration of how memory affects the potential depth and curvature of the right valleys in 1000 randomly parameterized polynomial models. 
    (a) Scatter plots compare the potential curvature and potential depth in the positive basin of attraction for models with memory values of 0 and 0.2. 
    (b) Distributions of potential curvature for models with and without memory. A lower magnitude of curvature indicates a flatter basin bottom. Thus, systems with memory tend to have flatter basins of attraction, characterized by less sharply curved contours.
    }
    \label{fig: Dp_Curv}
\end{figure}

\begin{figure}[ht!]
    \centering
    \includegraphics[width=1\textwidth]{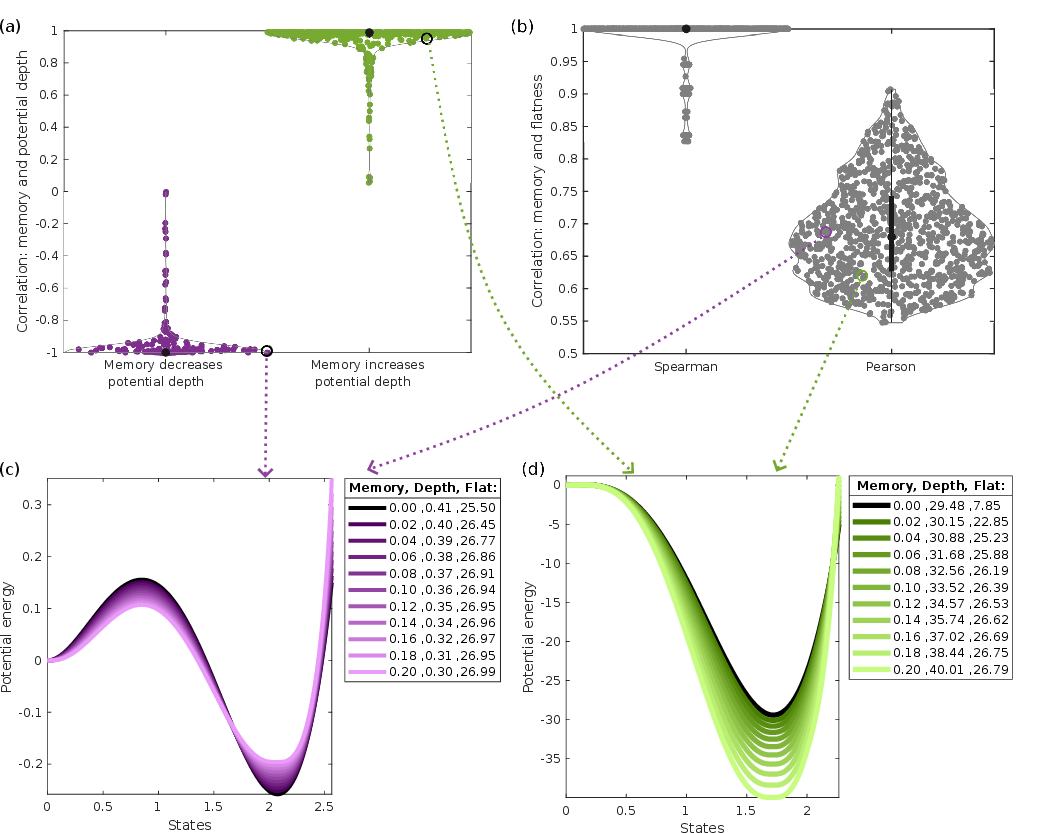}
    \caption{
    Violin plots illustrating the correlation between memory and two stability landscape metrics: (a) potential depth and (b) potential curvature of the right valley (corresponding to the positive stable state) across 1000 randomly parameterized polynomial models. In panel (a), both positive and negative correlations are shown: green dots indicate models where memory increases potential depth, while purple dots indicate models where memory decreases it. Panel (b) shows that the correlation is always positive, meaning memory consistently leads to a flatter basin of attraction (i.e., reduced curvature). 
    Panels (c) and (d) present representative models: (c) shows a case where increasing memory reduces potential depth, and (d) shows a case where memory increases potential depth. For each, the effect of different memory strengths on the stability landscape, potential depth, and flatness is demonstrated.
}

    \label{fig: Cor_Mem_dp_Fl}
\end{figure}

As memory strength increases, flatness grows monotonically at the basin bottom, indicating that memory consistently buffers system motion. In contrast, the effect on depth is heterogeneous: some basins become deeper (enhanced resilience), while others become shallower (reduced resilience). These divergent trends are summarized in Figures \ref{fig: Dp_Curv} and Figure~\ref{fig: Cor_Mem_dp_Fl} quantifies the correlations between memory strength and each landscape metric. Collectively, our results demonstrate that memory universally softens the slopes of stability landscapes, even as its impact on basin depth depends on model specifics.

\FloatBarrier

\subsection*{Impact of memory on resilience and resistance}

\begin{figure}[ht!]
    \centering
    \includegraphics[width=.5\textwidth]{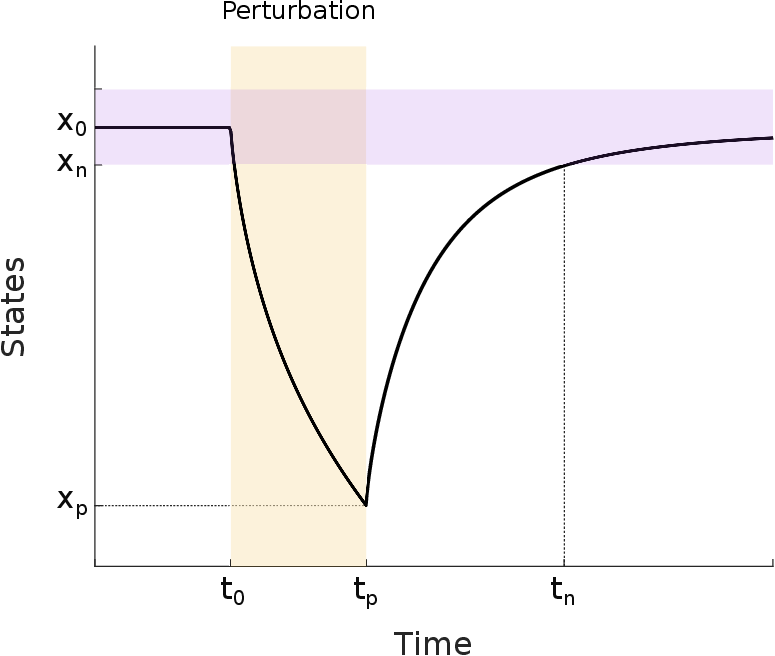}
    \caption{Illustration of the resilience with a quantitative approach. The system's original stable state is denoted as $X_0$. At the time $t_0$, a perturbation occurs, lasting until $t_p$, resulting in a change in the system state of $|X_0-X_p|$ during the time interval $(t_0, t_p)$. Resilience is quantified as the recovery rate index, which measures how quickly the system returns to a vicinity of $X_0$, specifically denoted as $X_n$, after the perturbation. The threshold for this vicinity is defined as $(X_{S2} - X_U)\times10^{-4}$, where $X_{S2}$ is the positive stable state (same as $X_0$ here) and $X_U$ is the unstable point between the two stable states. 
    The recovery rate is calculated based on how long it takes for the system to return to the original stable state, $X_n$, compared to the magnitude of the state change caused by the disturbance, represented by $|X_0-X_p|$. The mathematical definition of the resilience (or recovery rate) is as follows~\cite{Shade_fmicb2012}:
    \( \displaystyle
    \text{Resilience}=\frac{\left(\frac{|X_0-X_p|-|X_0-X_n|}{|X_0-X_p|+|X_0-X_n|} \right)}{t_n-t_p},
    \) 
    where $t_n$ represents the time at which the system state reaches the determined vicinity of the original stable state, indicating the completion of recovery. 
    }
    \label{fig: RL}
\end{figure}

Having shown that memory reshapes stability landscapes, we next assess how these changes translate into resilience and resistance. Following the definitions given in the main manuscript, we quantify:
{Resilience} is the rate of recovery back to the positive equilibrium after a brief, pulse‐type perturbation.
{Resistance} is the maximum perturbation amplitude the system can absorb without tipping into the extinction state (Figure~\ref{fig: RL}).

\begin{figure}[ht!]
    \centering
    \includegraphics[width=1\textwidth]{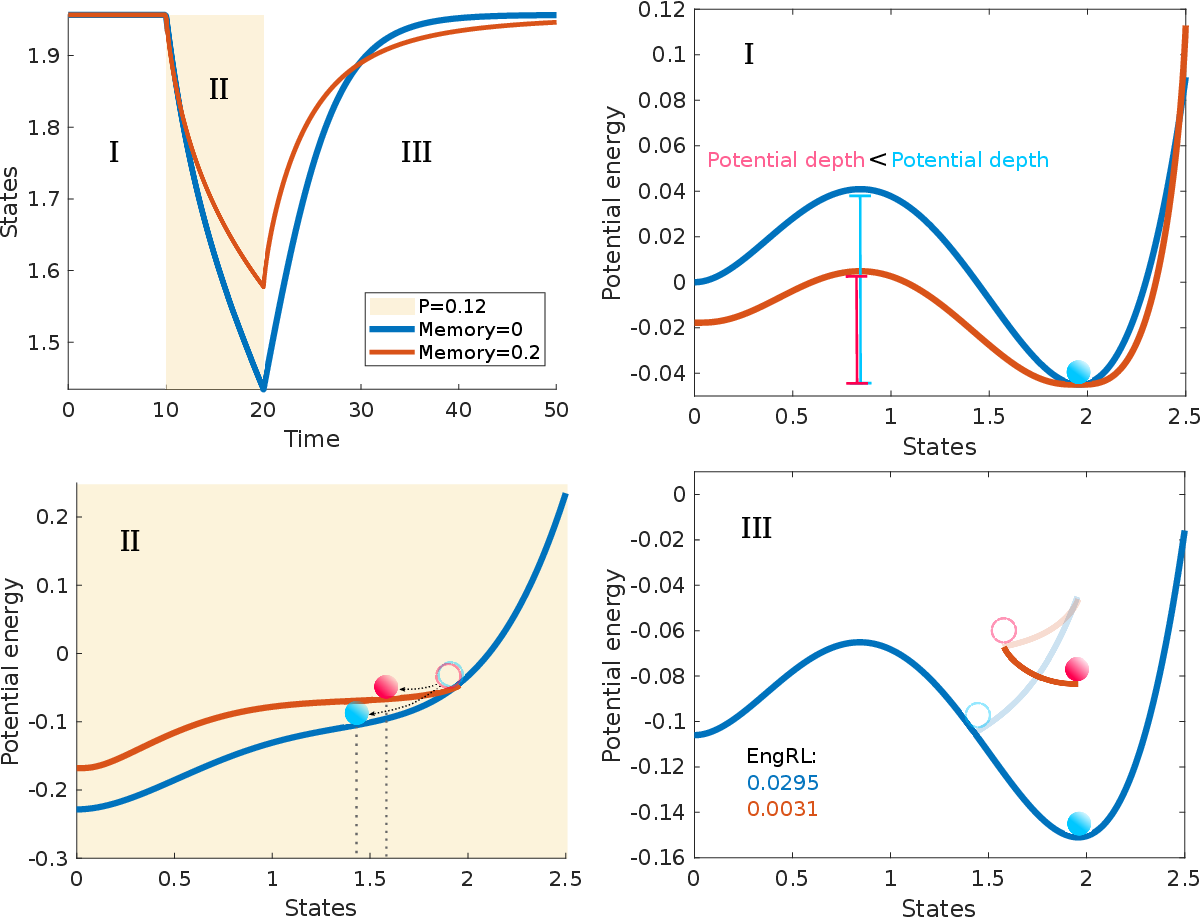}
    \caption{
    Impact of memory on the resilience of the herbivory model when memory reduces potential depth.
    (Top left) Model dynamics with and without memory are shown under three scenarios: (I) the stable state, (II) response to a perturbation (0.12 added to parameter $B$; parameters: $r=0.8$, $K=3$, $A=0.2$, $B=0.6$), and (III) recovery to the original stable state.
    \textbf{(I)} Comparison of potential energy in the stable state reveals that the potential depth is smaller in the presence of memory.
    \textbf{(II)} Following the perturbation, the system exhibits a single basin of attraction at 0.
    \textbf{(III)} During recovery, the model without memory has a recovery rate of 0.0295, while the model with memory recovers more slowly, at a rate of 0.0031. The relative effect of memory on resilience is quantified as the fraction derived from subtracting the recovery rates and dividing by their sum, yielding a value of -0.8088. 
    }
    \label{fig: mem_herb_RL1}
\end{figure}

\begin{figure}[ht!]
    \centering
    \includegraphics[width=1\textwidth]{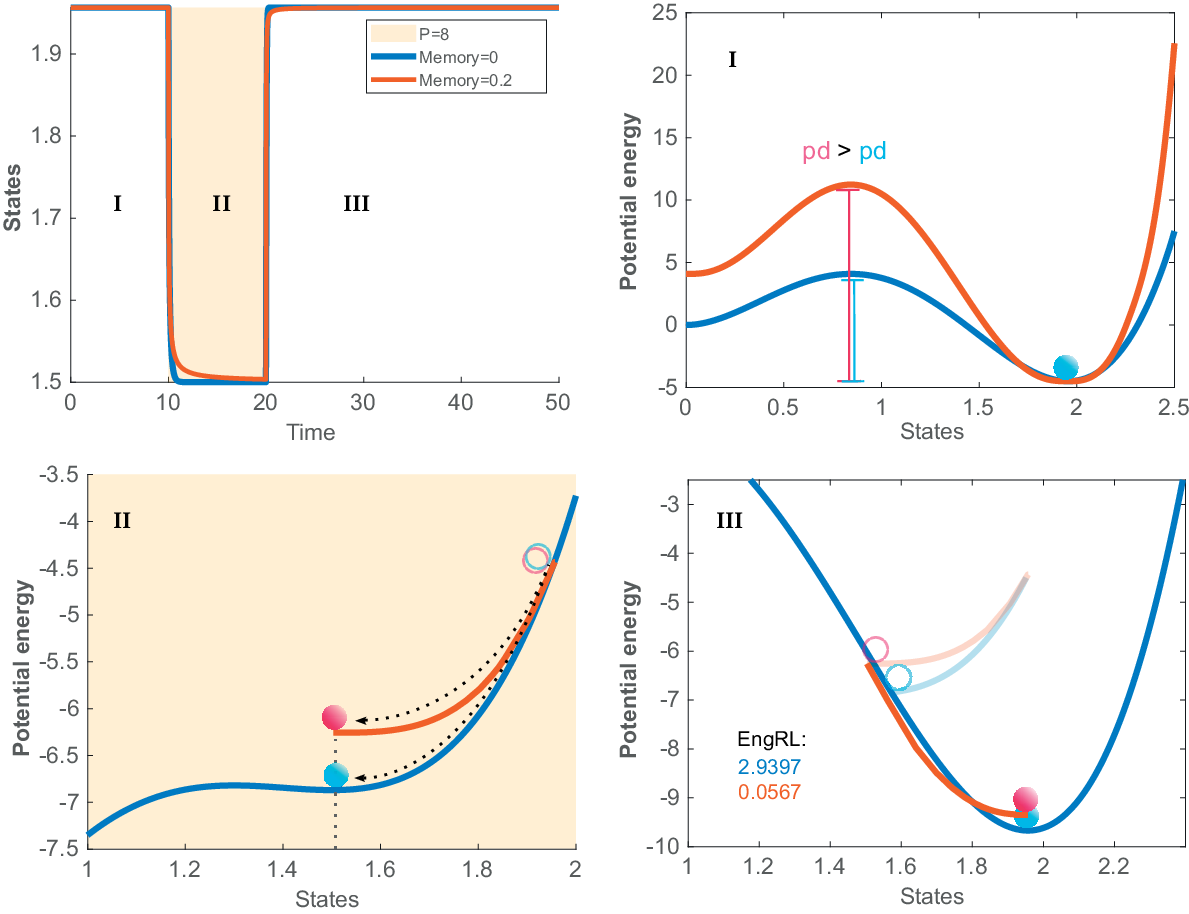}
    \caption{
    Impact of memory on the resilience of the herbivory model when memory increases potential depth.
    (Top left) Dynamics of the model with and without memory are shown under three scenarios: (I) stable state, (II) response to a perturbation (8 added to parameter $B$; parameters: $r=80$, $K=3$, $A=0.2$, $B=60$), and (III) recovery to the original stable state. 
    \textbf{(I)} In the stable state, the potential energy comparison shows that memory leads to a greater potential depth than the case without memory. 
    \textbf{(II)} After the perturbation, the system remains bistable, but a new basin of attraction emerges around 1.5. 
    \textbf{(III)} During recovery, the model without memory exhibits a recovery rate of 2.9397, while the model with memory recovers more slowly at 0.0567. The relative effect of memory on resilience is calculated as the fraction of the difference in recovery rates divided by their sum, resulting in -0.9622.
    }   
    \label{fig: mem_herb_RL2}
\end{figure}

We begin with the herbivory model under two contrasting memory scenarios, one in which memory shallows the basin (Figure~\ref{fig: mem_herb_RL1}) and one in which memory deepens it (Figure~\ref{fig: mem_herb_RL2}). In each case, we displace the system from its positive stable state via a short pulse and then track both the maximum excursion away from equilibrium and the subsequent return time. Although memory increases resistance, producing smaller deviations from equilibrium, it invariably slows the recovery process. This pattern of increased resistance accompanied by reduced recovery speed is evident in both shallow- and deep-basin scenarios, illustrated by the altered trajectory paths on the reconstructed stability landscapes (Figures~\ref{fig: mem_herb_RL1}-\ref{fig: mem_herb_RL2}).

\begin{figure}[ht!]
    \centering
    \includegraphics[width=1\textwidth]{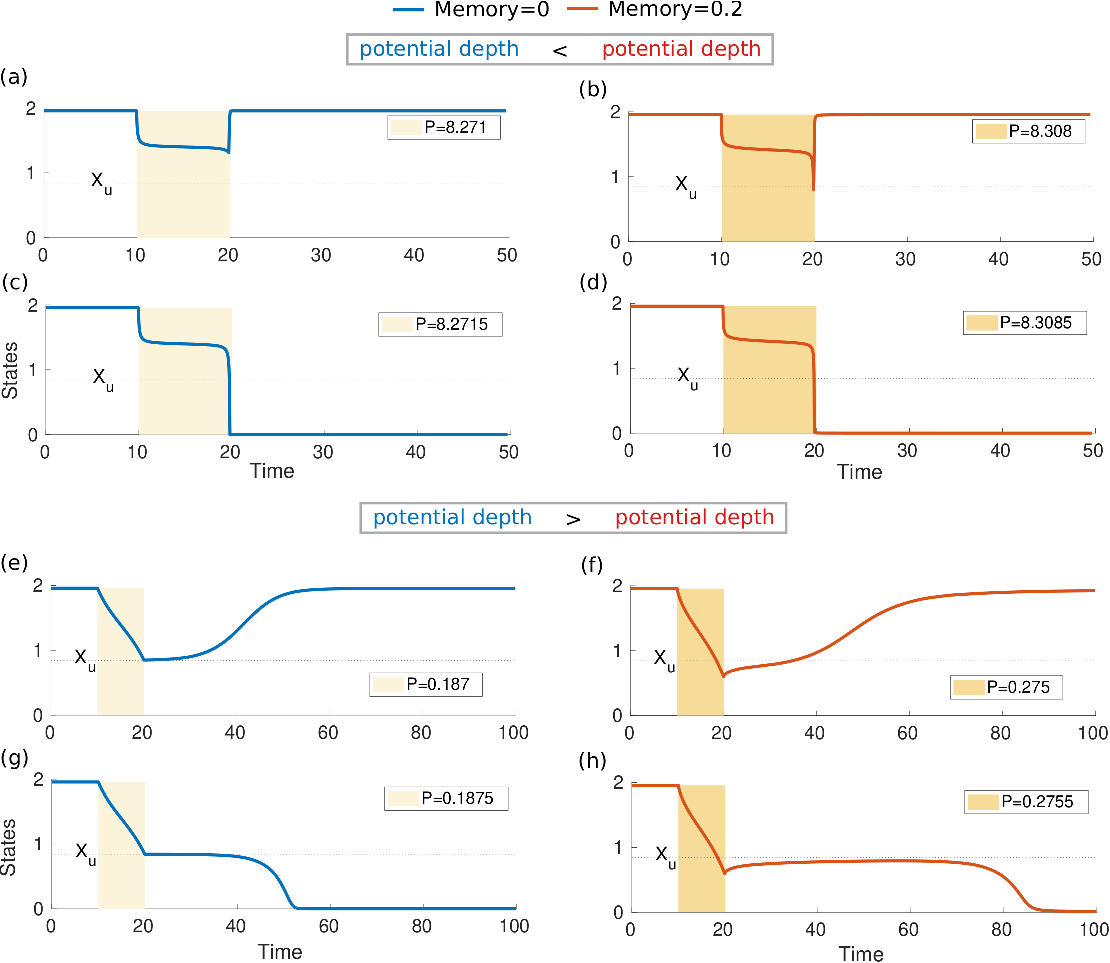}
    \caption{
    Memory increases resistance in the herbivory model, as shown here. Panels (a), (b), (e), and (f) depict the magnitude of perturbation that the system can tolerate in two scenarios: (a) and (e) without memory, and (b) and (f) with memory. Panels (a–d) correspond to the case where memory increases potential depth, while panels (e–h) represent the scenario where memory decreases potential depth. A small perturbation can shift the system to an alternative stable state: (c) and (g) show this transition without memory, and (d) and (h) with memory. The perturbation $P$ quantifies the minimum change in herbivore population density ($B$) required to trigger a state shift, measured with a precision of $5 \times 10^{-4}$. Comparing (a) and (b), the difference in resistance (the value of $P$) is 0.0370, with a relative effect of memory of 0.0022. For (e) and (f), the resistance difference is 0.0880, and the relative effect is 0.1905. The parameter values used are identical to those in Figures~\ref{fig: mem_herb_RL1} and~\ref{fig: mem_herb_RL2}.
    }
    \label{fig: mem_herb_RS}
\end{figure}

\begin{figure}[ht!]
    \centering
    \includegraphics[width=1\textwidth]{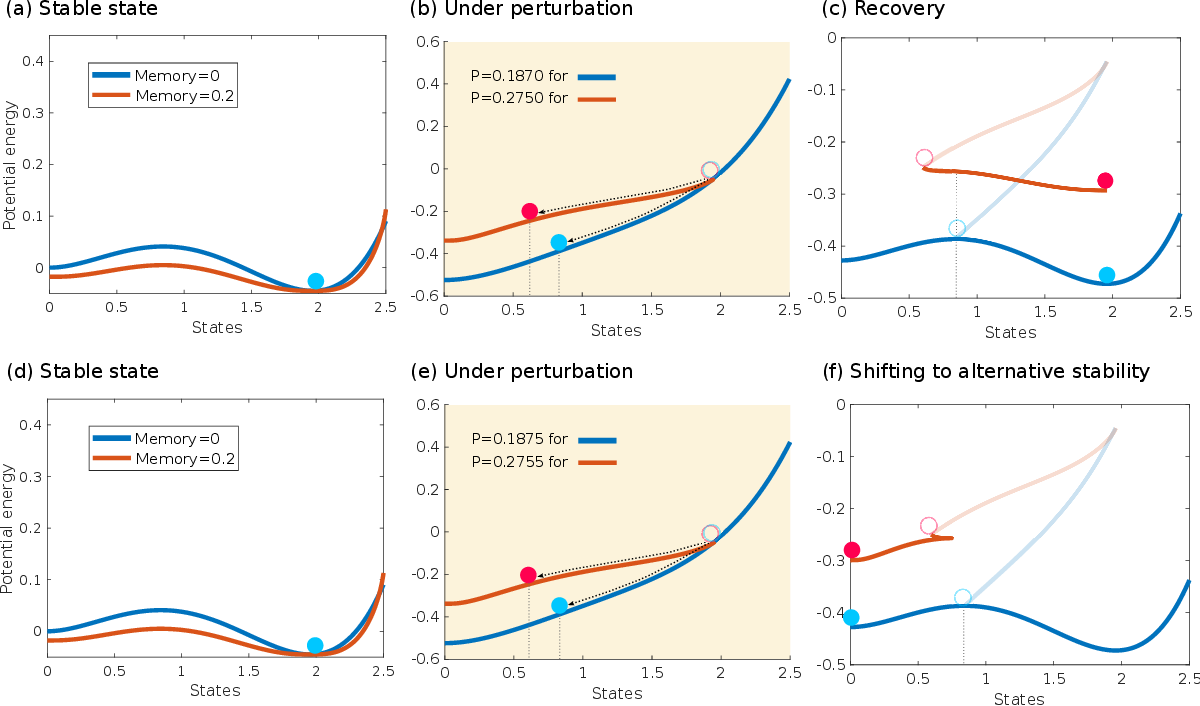}
    \caption{
    Illustration of potential energy for the model dynamics presented in Figure~\ref{fig: mem_herb_RS}(e–h), comparing cases with and without memory under three scenarios: 
    (a, d) stable state, 
    (b, e) response to a perturbation, and 
    (c) recovery to the original stable state or 
    (f) transition to an alternative stable state.
    }
    \label{fig: land_herb_RS}
\end{figure}

To further quantify resistance, we determine, for each model, the minimum perturbation strength required to push the system across the hill into the zero state. As shown in Figure~\ref{fig: mem_herb_RS}, introducing memory always raises this critical perturbation threshold, confirming its uniformly positive effect on resistance regardless of whether it increases or decreases basin depth.

\begin{figure}[ht!]
    \centering
    \includegraphics[width=1\textwidth]{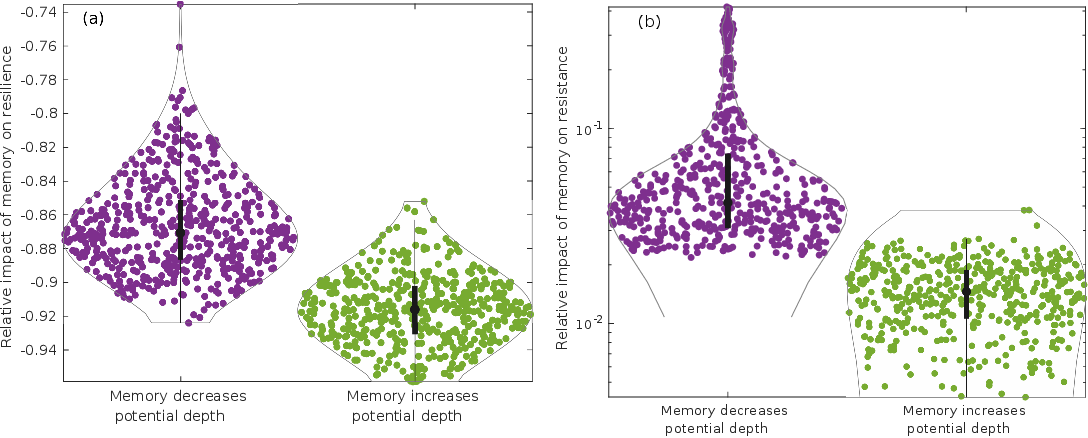}
    \caption{
    Illustration of how memory affects the resilience and resistance of 1000 randomly parameterized polynomial models. Models are categorized into two groups based on the effect of memory on potential depth: those where memory increases potential depth (green dots) and those where it decreases potential depth (purple dots).
    {(a)} Violin plots show the relative effect of memory on resilience for both groups. Memory generally reduces resilience, with a stronger negative impact observed when memory increases potential depth.
    {(b)} Violin plots display the relative effect of memory on resistance. Here, memory typically increases resistance, with the effect being more pronounced in models where memory decreases potential depth.
    }
      \label{fig: mem_RL_RS}
\end{figure}

Finally, we extend these analyses to our full ensemble of one‐dimensional, bistable models. For each sample, we compute the relative change in recovery rate and in critical perturbation amplitude (i.e., the difference divided by the sum, to normalize across scales). Figure~\ref{fig: mem_RL_RS} demonstrates that memory consistently diminishes resilience and enhances resistance across all models. Moreover, the magnitude of these effects correlates with the sign of the depth change: resilience loss is most pronounced when memory deepens the basin, whereas resistance gain is greatest when memory shallows it. Together, these results reveal a universal inverse relationship between resistance and resilience induced by memory effects.

\FloatBarrier
\subsection*{Basin sharpness index and relative memory effects}

To explore how memory’s influence varies across different basin geometries, we first introduce a quantitative measure of basin shape. The chosen models allow the valleys to stretch and alter their gradient metaphorically, showcasing a dynamic interplay of steep and flat basins. Empirically, we observe a negative correlation between depth and flatness in our positive basins: deeper wells tend to be steeper, while shallower wells are flatter. 
We therefore define a basin sharpness index as the average of the normalized depth and the normalized absolute curvature at the basin minimum, where each metric is scaled to lie between 0 and 1. An index near 1 corresponds to a deep, steep basin; an index near 0 indicates a shallow, flat basin.

\begin{figure}[ht!]
    \centering
    \includegraphics[width=.6\linewidth]{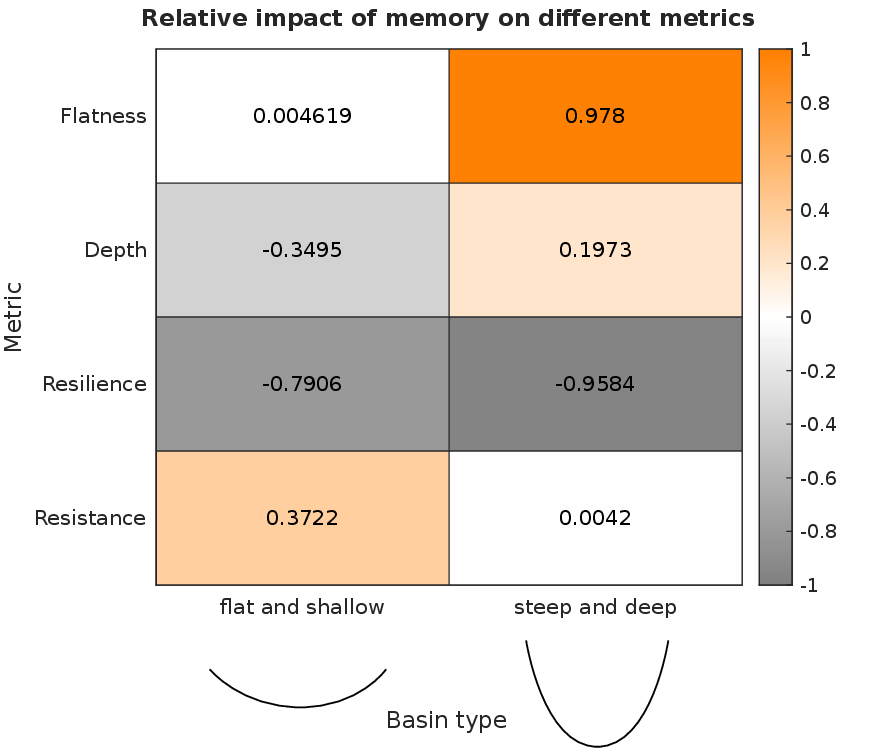}
    \caption{
    Relative impact of memory on flatness, depth, resilience, and resistance as a function of the basin sharpness index, shown for two types of landscapes: flat and shallow versus steep and deep. Results are based on 1000 randomly parameterized polynomial models. The relative impact of memory on flatness is pronounced in models with steep, deep basins of attraction, but close to zero (less than 0.005) in flat, shallow basins. For potential depth, memory has a positive effect on models with steep and deep valleys, but a negative impact on flat and shallow ones. Memory’s influence on resilience is generally more significant than on resistance, and its relative impact is not much different based on the basin types. Memory has a greater impact on resistance in models with flat, shallow basins, while its effect in steep and deep landscapes is minor (less than 0.005).
}
    \label{fig: detailed results}
\end{figure}

Using this index, we partition our model ensemble into \textit{steep–deep} and \textit{flat–shallow} basins and then compute, for each group, the relative change induced by memory in four key quantities: flatness, depth, resilience, and resistance. Figure~\ref{fig: detailed results} summarizes these group‐level comparisons. We find that:

1) Flatness is more strongly influenced by memory in steep–deep basins than in flat–shallow ones.

2) Depth is more sensitive to memory in flat–shallow basins.

3) Resilience generally declines under memory, with the largest relative decrease occurring in basins where memory deepens.

4) Resistance increases under memory most noticeably in flat–shallow basins, and is least affected in steep–deep basins.

5) Overall, the proportional decrease in resilience exceeds the proportional increase in resistance, indicating that memory’s slowing of recovery is a stronger effect than its bolstering of resistance.

\begin{figure}[ht!]
    \centering
    \includegraphics[width=1\textwidth]{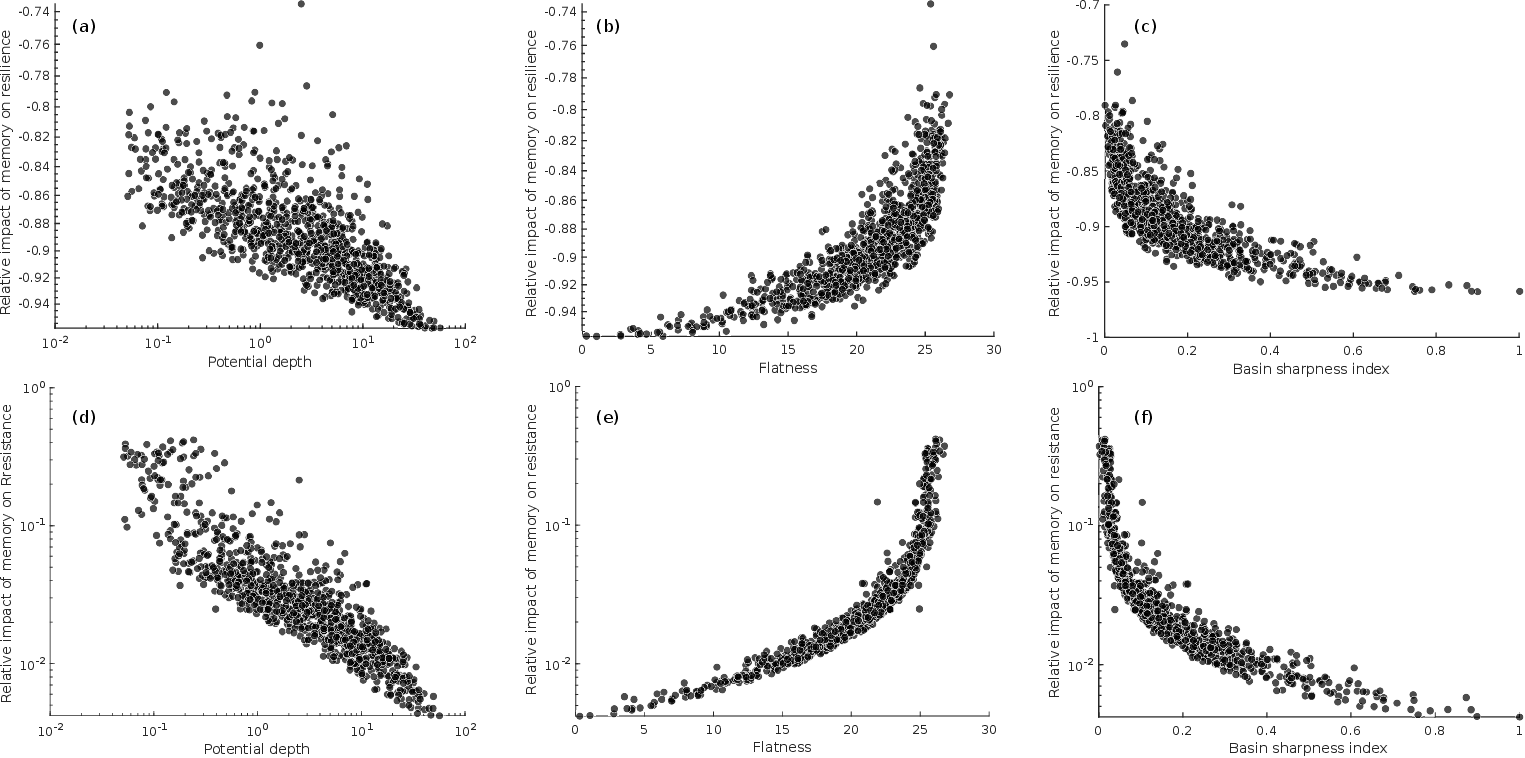}
    \caption{
    Illustration of how memory affects resilience and resistance in relation to potential depth and flatness in 1000 randomly parameterized polynomial models.
    Scatter plots showing the relationship between: (a) the relative impact of memory on resilience and the potential depth of models without memory,
    (b) the relative impact of memory on resilience and the flatness of models without memory,
    (c) the relative impact of memory on resilience and the basin sharpness index of models without memory,
    (d) the relative impact of memory on resistance and the potential depth of models without memory,
    (e) the relative impact of memory on resistance and the flatness of models without memory, and 
    (f) the relative impact of memory on resistance and the basin sharpness index of models without memory.
    }
    \label{fig: Mem_RL_RS_Basin}
\end{figure}

\begin{figure}[ht!]
    \centering
    \includegraphics[width=1\textwidth]{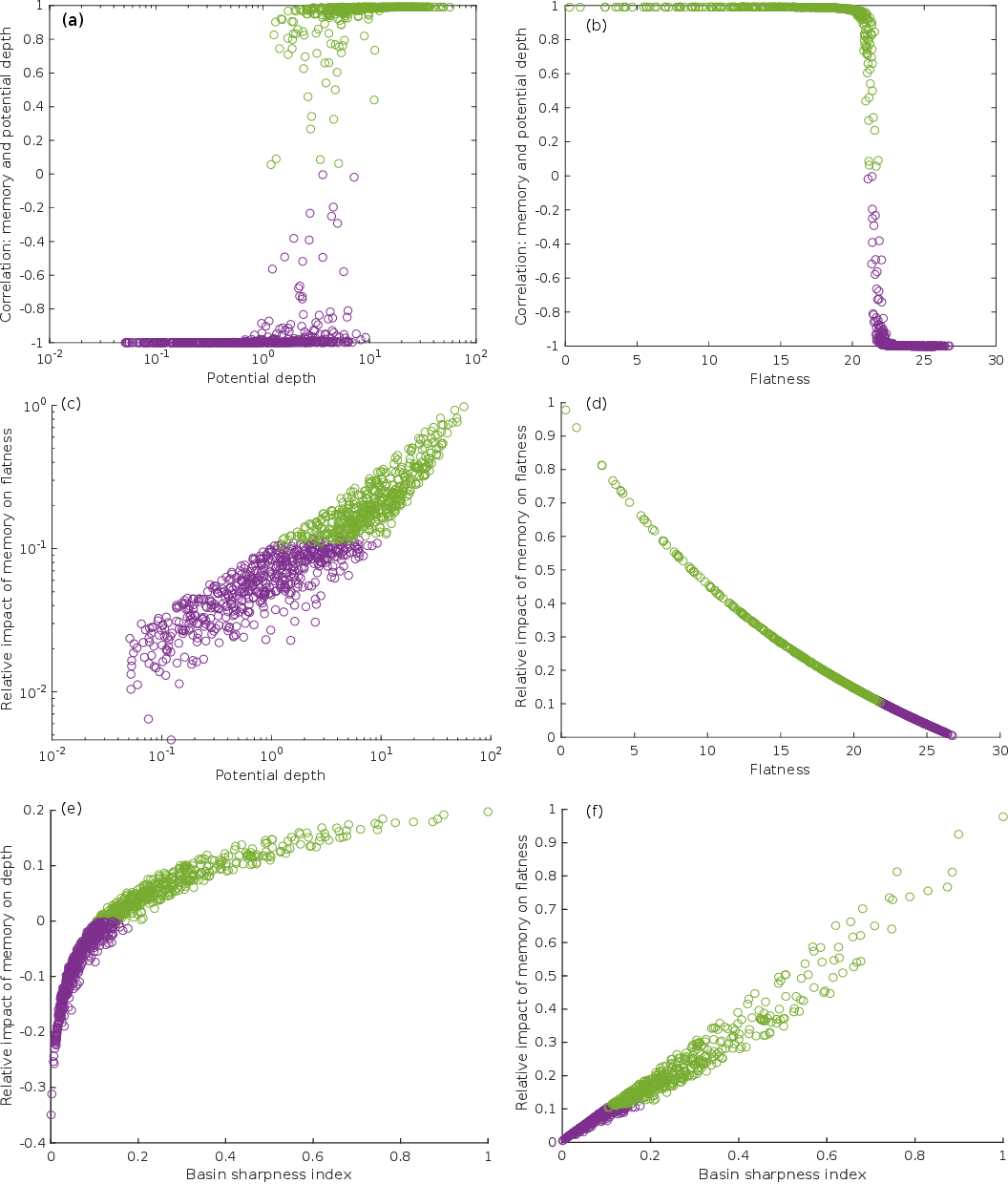}
    \caption{
    Illustration of how memory influences potential depth and flatness in the right valleys of 1000 randomly parameterized polynomial models. Models are grouped by the effect of memory on potential depth: green circles indicate models where memory increases potential depth, and purple circles indicate models where memory decreases it. 
    (a) Scatter plot showing the relationship between memory-induced changes in potential depth and the baseline potential depth of models without memory (x-axis).
    {(b)} Scatter plot showing memory-induced changes in potential depth versus the baseline flatness of models without memory (x-axis).
    {(c)} Scatter plot of the relative impact of memory on flatness as a function of baseline potential depth (x-axis).
    {(d)} Scatter plot of the relative impact of memory on flatness as a function of baseline flatness (x-axis).
    {(e)} Scatter plot of the relative impact of memory on potential depth as a function of the basin sharpness index without memory (x-axis).
    {(f)} Scatter plot of the relative impact of memory on flatness as a function of the basin sharpness index without memory (x-axis).
    }
    \label{fig: Mem_fl_dp}
\end{figure}

To visualize these continuous trends, Figures~\ref{fig: Mem_RL_RS_Basin} and \ref{fig: Mem_fl_dp} present scatter plots of memory strength versus landscape metrics and versus ecological responses. These plots confirm that memory’s buffering effect on motion (flatness) and the inverse relationship between resilience and resistance depend systematically on the underlying basin geometry.

\FloatBarrier
\subsection*{Linking stability landscape geometry to the magnitude of memory effects}

\begin{figure}[ht!]
    \centering
    \includegraphics[width=\linewidth]{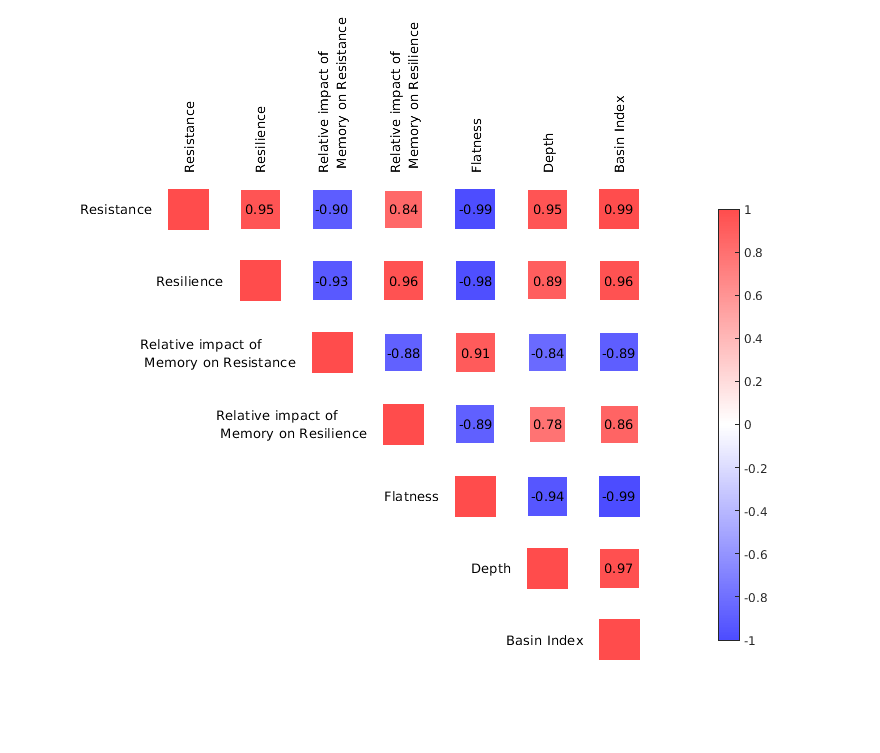}
    \caption{
    Summary of the detailed findings for 1000 randomly parameterized polynomial models. 
    Spearman correlations were calculated to examine relationships among key resilience metrics: potential depth, flatness, basin sharpness index, resilience, resistance, and the relative impact of memory on resilience and resistance. A negative correlation between depth and flatness indicates that shallower valleys are associated with flatter basins of attraction. Both resilience and resistance are negatively correlated with flatness and positively correlated with depth. The relative impact of memory on resistance exhibits an opposite trend, being negatively correlated with depth and positively correlated with flatness, whereas the relative impact of memory on resilience is positively correlated with depth and negatively correlated with flatness.
    Models with steeper and deeper landscapes, characterized by greater potential depth and flatter attraction basins, exhibit higher resilience and resistance, supporting the expectation that more resilient systems are associated with deeper and steeper stability landscapes.
    }
    \label{fig: summary results}
\end{figure}

To consolidate all our results, we calculated the pairwise correlations among the core variables—depth, flatness, basin sharpness index, resilience, resistance, and their relative changes under memory—and display them in Figure~\ref{fig: summary results}.

First, the geometric trade‐off in one‐dimensional bistable models is immediately apparent: basin depth and flatness are strongly negatively correlated. Depth correlates positively with the basin sharpness index, while flatness correlates negatively with it. These relationships reaffirm that deeper wells tend to be steeper, whereas shallower basins are flatter.

\begin{figure}
    \centering
    \includegraphics[width=1\textwidth]{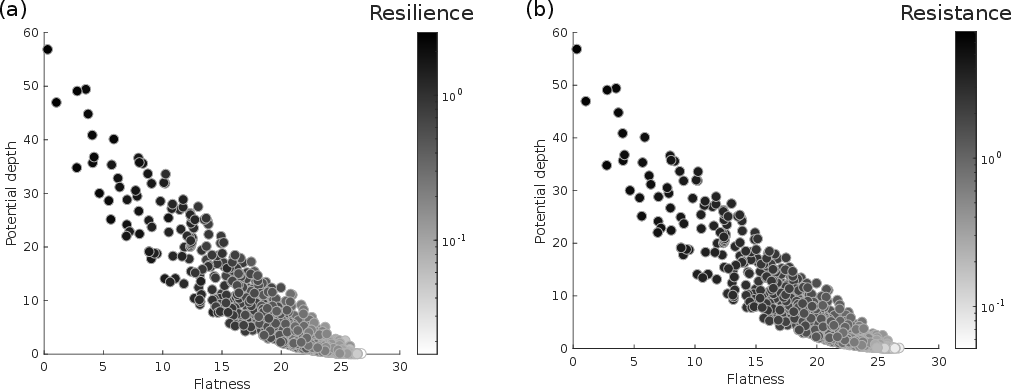}
    \caption{
    Scatter plots of potential depth versus flatness for 1000 randomly parameterized polynomial models without memory, where a grayscale gradient indicates the values of (a) resilience and (b) resistance.
    }
    \label{fig: dp_fl_RL_RS}
\end{figure}

Next, mapping ecological function onto landscape shape, we observed that both resilience and resistance increase with basin depth but decrease with flatness. In other words, deeper, steeper basins not only recover more quickly but also tolerate larger perturbations. Figure~\ref{fig: dp_fl_RL_RS} illustrates this visually: each point represents a model’s flatness versus depth, with the gray-scale shading indicating its resilience or resistance level.
Together, these patterns confirm our intuition—and validate our metrics—that deeper, steeper basins indeed confer both greater resistance and faster recovery, showing that our measures align with biological expectations.

Finally, we examined how memory modulates these patterns, which is not a trivial fact. The proportional increase in resistance under memory correlates positively with flatness and negatively with depth, indicating that memory most enhances resistance in flat, shallow basins. In contrast, the proportional decline in resilience under memory correlates negatively with flatness and positively with depth, showing that memory most undermines resilience in deep, steep basins. These contrasting correlations illustrate that memory enhances resistance primarily in shallow, flat landscapes, while more markedly diminishing resilience in deeper, steeper ones.

\FloatBarrier
\section*{S4: Impact of memory on dynamical systems under exogenous perturbations}\label{sec: sup4}
\subsection*{Within basin variability under additive noise}

\subsubsection*{Model and local linearization}
We study a one-dimensional double well with fractional dynamics
\[
\mathcal{D}^{\alpha}x(t)=f(x(t))+\sigma\,\xi(t), \qquad f(x)=x-x^3, \qquad 0<\alpha\le 1,
\]
where $\mathcal{D}^{\alpha}$ is the Caputo derivative, $\xi(t)$ is Gaussian white noise with unit spectral density, and $\sigma>0$ is kept fixed when $\alpha$ changes. The deterministic term $f(x)$ is the \emph{drift} (the systematic restoring force that pushes trajectories toward stable wells). The stochastic term $\sigma\,\xi(t)$ represents \emph{diffusion} in the sense of additive random fluctuations, with $\sigma$ controlling the noise intensity. In the memory free limit $\alpha=1$, this reduces to the standard Langevin form $\dot x=f(x)+\sigma\xi(t)$, equivalently $dx=f(x)\,dt+\sigma\,dW_t$, where $f$ is the drift and $\sigma$ is the diffusion amplitude (so the diffusion strength is proportional to $\sigma^2$).

We need to do linearization at the stable equilibrium inside a single basin. Thus, we define $y(t)=x(t)-x^\star$ as the deviation of the state from the stable equilibrium $x^\star\in\{\pm1\}$. Inside one basin, typical fluctuations are small, and the absolute position is less informative than the deviation from equilibrium. Studying $y$ focuses the analysis on within basin variability and yields a linear model that keeps the leading effects of stability strength and memory while allowing closed form characterizations of variance and correlation.

Because the Caputo derivative is linear and $x^\star$ is constant,
$\mathcal{D}^{\alpha}y(t)=\mathcal{D}^{\alpha}x(t)$~\cite{kilbas2006theory}. Expand the drift around $x^\star$:
\[
f(x^\star+y)=f(x^\star)+f'(x^\star)\,y+\mathcal O(y^2).
\]
Since $x^\star$ is an equilibrium, $f(x^\star)=0$, and for a stable equilibrium $f'(x^\star)<0$. Keeping only the linear term and substituting $x(t)=x^\star+y(t)$ into the stochastic system gives
\[
\mathcal{D}^{\alpha}y(t)=f'(x^\star)\,y(t)+\sigma\,\xi(t).
\]
Writing $\lambda:=-f'(x^\star)>0$ and moving terms to the left yields the
fractional Ornstein–Uhlenbeck (OU) equation~\cite{rybakov2020spectral}
\begin{equation}
D_t^{\alpha}y(t)+\lambda\,y(t)=\sigma\,\xi(t).
    \label{eq:frac_ou}
\end{equation}
For the double well $f(x)=x-x^3$ we have $f'(x)=1-3x^2$, hence
$f'(\pm1)=-2$ and $\lambda=2$.

The parameter $\lambda$ equals the local curvature of the potential
$V(x)=\tfrac14 x^4-\tfrac12 x^2$ at $x^\star$ because $f(x)=-V'(x)$, so larger $\lambda$ means stronger local restoring force. The order $\alpha\in(0,1]$ controls memory: $\alpha=1$ recovers the classical OU process with exponential relaxation, while $\alpha<1$ produces slower, Mittag–Leffler type relaxation and longer correlation~\cite{MAINARDI19961461}. The linear approximation is valid when typical deviations remain small compared to the basin width, which holds for sufficiently small noise $\sigma$ and time windows that do not include barrier crossings. Within this regime, the fractional OU model captures the leading behavior of within-basin variance and temporal correlation under additive noise while isolating the roles of $\lambda$ and $\alpha$.

\begin{figure}[t]
  \centering
  \includegraphics[width=0.78\linewidth]{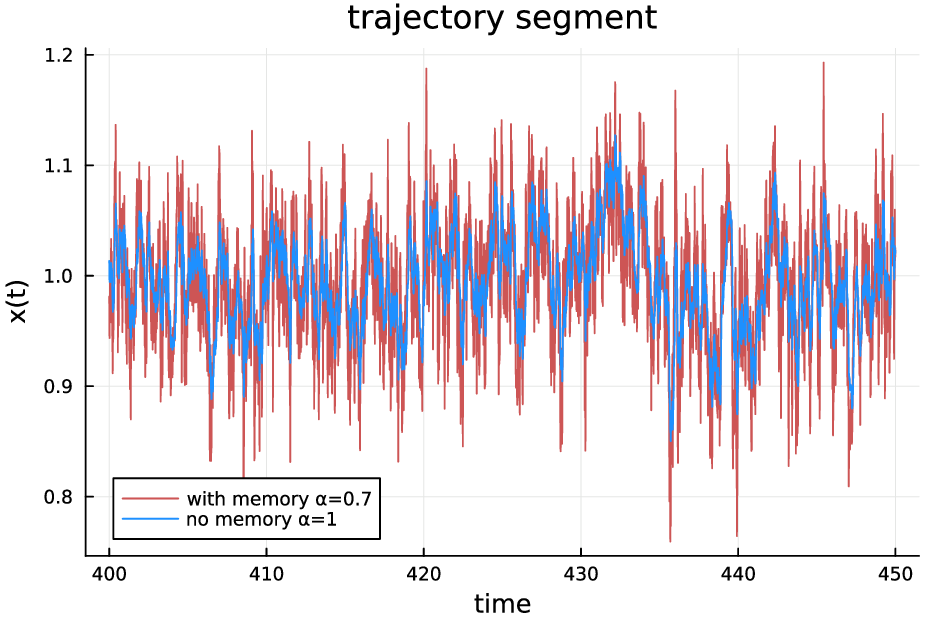}
  \caption{{Trajectory segment at native resolution.}
  The memory trajectory ($\alpha=0.7$) is visibly rougher than the classical case ($\alpha=1$), even though both remain confined to the same basin over this window.}
  \label{fig:traj}
\end{figure}

We first build intuition directly in the time series. Figure~\ref{fig:traj} displays a representative window from both simulations at native resolution. The dynamics with memory ($\alpha=0.7$) shows more rapid excursions riding on the same slow drift inside the well. This is the visual manifestation of weaker attenuation of fast components by the fractional operator, a point we quantify next in the frequency domain.

\subsubsection*{Frequency domain analysis: transfer function and spectrum}

We use Fourier transforms for linear time invariant models, which, moving to the frequency domain, turns differentiation and convolution into simple multiplication~\cite{stinga2023fractional, Welch1967}. This lets us treat the system as a filter driven by noise and read off how each frequency is amplified or suppressed.

Consider the linear fractional relaxation model driven by noise, schematically $\mathcal{D}^{\alpha} y(t) + \lambda\,y(t) = \xi(t)$, with $0<\alpha\le 1$, relaxation rate $\lambda>0$, and zero past $y(t)=0$ for $t<0$. Taking Fourier transforms with the convention $\mathcal{F}\{\mathcal{D}^{\alpha}y\}=(i\omega)^{\alpha}Y(\omega)$ gives the input output relation~\cite{stinga2023fractional}
\[
\bigl(\lambda + (i\omega)^{\alpha}\bigr) Y(\omega)=\Xi(\omega).
\]
The \emph{transfer function} is the frequency response of the system~\cite{hartley1998solution},
\[
H(\omega)=\frac{Y(\omega)}{\Xi(\omega)}=\frac{1}{\lambda+(i\omega)^{\alpha}},
\]
which tells how much a sinusoid at frequency $\omega$ is scaled and phase shifted by the dynamics.

For a stationary real signal $y(t)$, the power spectral density (PSD) $S_y(\omega)$ distributes the variance across frequencies and obeys~\cite{DROZDOV2007237}
\[
\mathrm{Var}(y)=\int_{0}^{\infty} S_y(\omega)\,d\omega
\]
for the one-sided convention (nonnegative frequencies only). A one-sided spectrum is standard when working with real data; it folds the negative frequencies onto the positive axis. For linear time invariant systems,
\[
S_y(\omega)=|H(\omega)|^2\,S_{\xi}(\omega),
\]
so the output spectrum equals the input spectrum times the gain squared.

With white input noise of level $\sigma^2$ under the one sided convention, $S_{\xi}(\omega)=\sigma^2$ for $\omega\ge 0$. Using $|(i\omega)^{\alpha}|^2=\omega^{2\alpha}$ and $\mathrm{Re}\{(i\omega)^{\alpha}\}=\omega^{\alpha}\cos(\tfrac{\pi\alpha}{2})$ \cite{DROZDOV2007237, lovejoy2019fractional},
\begin{equation}
S_y(\omega)=|H(\omega)|^2 S_\xi(\omega)
=\frac{\sigma^2}{\lambda^2+2\lambda\,\omega^{\alpha}\cos(\tfrac{\pi\alpha}{2})+\omega^{2\alpha}}.
\label{eq:psd}
\end{equation}
This shows explicitly that $S_y$ is the input level $\sigma^2$ modulated by the squared magnitude of the transfer function.

Here, two limits set the shape:
\[
\omega\to 0:\quad S_y(\omega)\to \frac{\sigma^2}{\lambda^2}\quad\text{(a flat low frequency plateau),}
\]
\[
\omega\to \infty:\quad S_y(\omega)\sim \sigma^2\,\omega^{-2\alpha}\quad\text{(a power law tail with slope }-2\alpha\text{ in log log).}
\]
The plateau reflects that very slow fluctuations are strongly damped by $\lambda$, so power saturates to a constant set by $\sigma$ and $\lambda$. The tail shows how quickly power decays at fast frequencies. For the classical case $\alpha=1$ the slope is near $-2$; for $\alpha<1$ it is flatter $-2\alpha$, which shifts more variance to higher frequencies and makes sample paths look rougher. Because variance is the area under the PSD, a flatter tail places more area in the fast band, matching Fig.~\ref{fig:traj}.

The height of the low frequency plateau equals $\sigma^2/\lambda^2$, so it provides a direct handle to estimate the ratio of noise level to relaxation rate from data, and it explains why different models with the same $(\sigma,\lambda)$ meet at low frequencies. Model differences then show up in the high frequency slope $-2\alpha$.

\begin{figure}[ht!]
  \centering
  \includegraphics[width=0.78\linewidth]{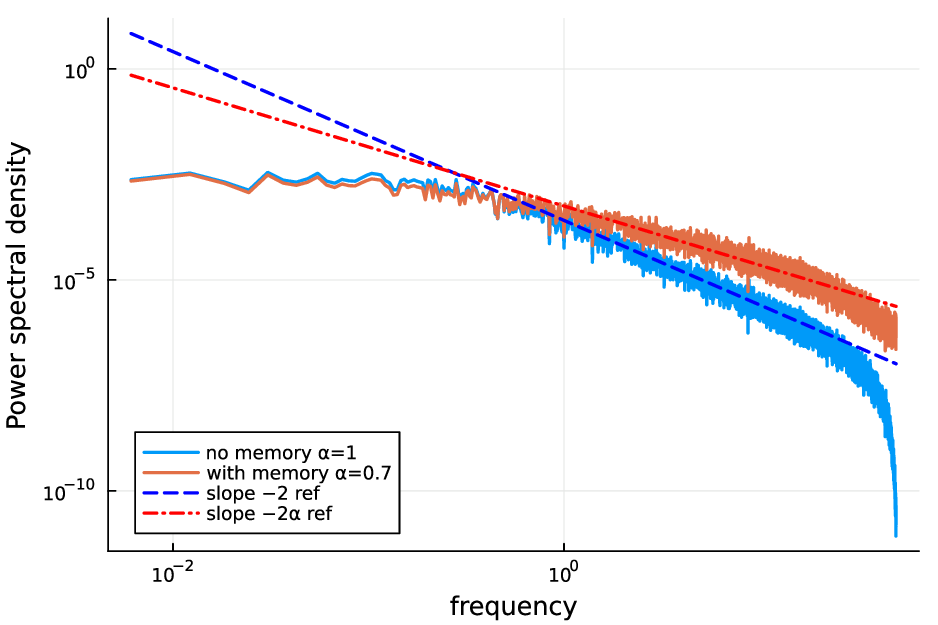}
  \caption{{Power spectral density.}
  Both spectra share the same low frequency plateau. In the tail, $\alpha=1$ decays with slope near $-2$, whereas $\alpha=0.7$ follows the flatter slope $-2\alpha\approx -1.4$ and lies above, confirming weaker attenuation of fast components. The roll off near the right edge is produced by the Nyquist limit and windowing in the Welch estimator~\cite{Welch1967}.}
  \label{fig:psd}
\end{figure}

With sampling step $h$, the highest resolvable frequency is the Nyquist $f_N=1/(2h)$, or in angular units $\Omega_c=\pi/h$. The bend toward the right edge of Fig.~\ref{fig:psd} is the finite sampling roll off near this cutoff.

\subsubsection*{Time domain analysis: impulse response and correlation}

For a linear time invariant system, the \emph{impulse response} $g(t)$ is the output produced by a unit impulse input $\xi(t)=\delta(t)$ under zero past. It is also called the Green function~\cite{lovejoy2019fractional}. For any input $\xi$, the output is the convolution
\[
y(t)=(g*\xi)(t)=\int_{0}^{t} g(t-u)\,\xi(u)\,du,
\]
thus, $g$ fully characterizes how the present output depends on past input.

To find the impulse response of \eqref{eq:frac_ou}, take the Laplace transform of $D_t^{\alpha}y+\lambda y=\xi$ with zero initial conditions. Using $\mathcal{L}\{D_t^{\alpha}y\}(s)=s^{\alpha}Y(s)$~\cite{kilbas2006theory},
\[
\bigl(s^{\alpha}+\lambda\bigr)Y(s)=\Xi(s)
\quad\Rightarrow\quad
H(s)=\frac{Y(s)}{\Xi(s)}=\frac{1}{s^{\alpha}+\lambda}.
\]
For an impulse input, $\Xi(s)=1$, so $Y(s)=H(s)$ and $g$ is the inverse Laplace transform of $H$:
\[
g(t)=\mathcal{L}^{-1}\!\left\{\frac{1}{s^{\alpha}+\lambda}\right\}(t)
=t^{\alpha-1}E_{\alpha,\alpha}\!\big(-\lambda t^{\alpha}\big), \qquad t>0.
\]
This identity follows from the general formula
\[
\mathcal{L}^{-1}\!\left\{\frac{s^{\alpha-\beta}}{s^{\alpha}+a}\right\}(t)=t^{\beta-1}E_{\alpha,\beta}\!\big(-a\,t^{\alpha}\big),
\]
Where is $E_{\alpha,\alpha}$, is the two parameter Mittag Leffler function, defined as~\cite{kilbas2006theory}
\[
E_{\alpha,\beta}(z)=\sum_{k=0}^{\infty} \frac{z^{k}}{\Gamma(\alpha k+\beta)},
\]
so here
\[
E_{\alpha,\alpha}(z)=\sum_{k=0}^{\infty} \frac{z^{k}}{\Gamma(\alpha k+\alpha)}.
\]
In particular $E_{\alpha,\alpha}(0)=1/\Gamma(\alpha)$.

For white forcing with one sided level $\sigma^{2}$, the variance equals the noise level times the $L^{2}$ energy of $g$,
\[
\mathrm{Var}[y]=\sigma^{2}\int_{0}^{\infty} g(t)^{2}\,dt.
\]
Using $E_{\alpha,\alpha}(-\lambda t^{\alpha})=1/\Gamma(\alpha)+\mathcal{O}(t^{\alpha})$ as $t\downarrow 0$,
\[
g(t)\sim \frac{t^{\alpha-1}}{\Gamma(\alpha)},\qquad
g(t)^{2}\sim \frac{t^{2\alpha-2}}{\Gamma(\alpha)^{2}}.
\]
The integral $\int_{0}^{\varepsilon} t^{2\alpha-2}\,dt$ is finite if and only if $2\alpha-2>-1$, that is $\alpha>\tfrac{1}{2}$. Hence the short time contribution is integrable for $\alpha>\tfrac{1}{2}$ and diverges for $\alpha\le\tfrac{1}{2}$. With sampled data of step $h$ the lower limit is effectively $t\approx h$, so the divergence is cut off and the measured variance increases as $\alpha$ decreases.

\begin{figure}[ht!]
  \centering
  \includegraphics[width=0.78\linewidth]{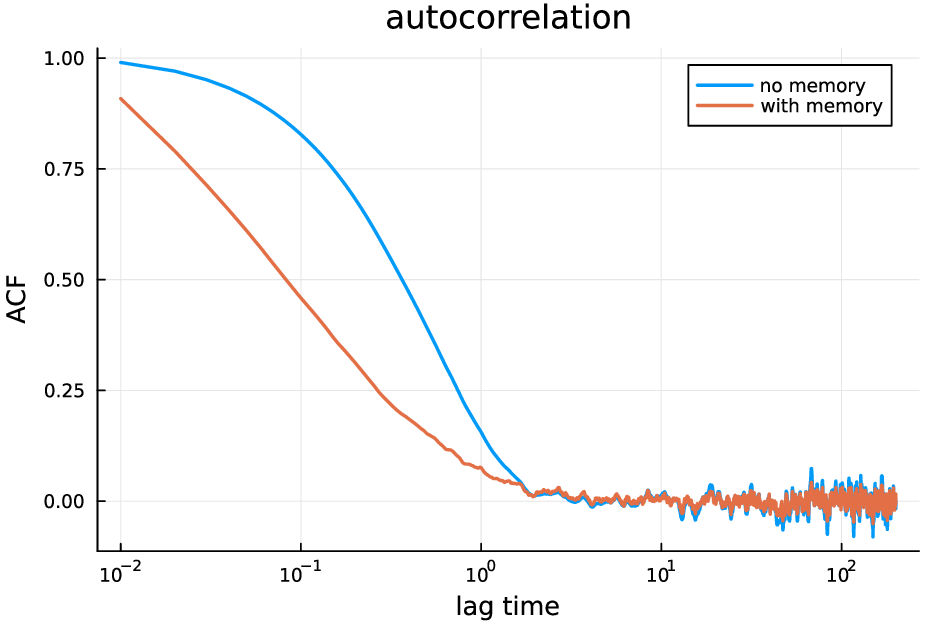}
  \caption{{Autocorrelation of the raw series.}
  The fractional series ($\alpha=0.7$) decorrelates more quickly at short lags than the classical series ($\alpha=1$). A slow tail, not fully visible at this scale, remains and will drive differences after coarse graining.}
  \label{fig:acf}
\end{figure}

For white forcing, the autocovariance is the self convolution of $g$,
\[
R_y(\tau)=\sigma^{2}\int_{0}^{\infty} g(u)\,g(u+\tau)\,du,
\qquad
\rho(\tau)=\frac{R_y(\tau)}{R_y(0)}.
\]
This produces a fast initial drop in $\rho(\tau)$ at very short lags (reflecting the strong weight of $g$ near zero time) followed by a slow tail from the long memory of the fractional kernel $t^{\alpha-1}E_{\alpha,\alpha}(-\lambda t^{\alpha})$.
Figure~\ref{fig:acf} shows the fast drop clearly: at the very shortest lags the fractional series stays below the classical one, which matches the spectral flattening in Fig.~\ref{fig:psd}. The slow tail persists beyond the plotted window and will become important once we change resolution.

\subsubsection*{Scale dependent statistics and coarse graining}
To summarize persistence at a chosen scale we use the windowed integrated autocorrelation time~\cite{geyer1992practical, Welch1967, flyvbjerg1989error}
\[
\tau_{\mathrm{int}}(W)=\int_0^{W}\rho(\tau)\,d\tau,
\]
estimated in discrete time with step $h$ by $h\big(1+2\sum_{k=1}^{\lfloor W/h\rfloor}\widehat{\rho}_k\big)$.

To study how resolution changes statistics we apply block averaging over $m$ consecutive samples~\cite{flyvbjerg1989error},
\[
\bar x_j=\frac{1}{m}\sum_{i=1}^{m}x_{(j-1)m+i},
\]
which suppresses fluctuations faster than the block length $mh$ and changes the effective sampling step to $mh$. A standard calculation gives~\cite{Welch1967}
\[
\mathrm{Var}(\bar x_j)=\frac{\gamma_0}{m}\left[1+2\sum_{k=1}^{m-1}\Big(1-\frac{k}{m}\Big)\rho_k\right],
\]
so the block variance is a weighted sum of short to moderate lags.

\begin{figure}[th!]
  \centering
  \includegraphics[width=0.78\linewidth]{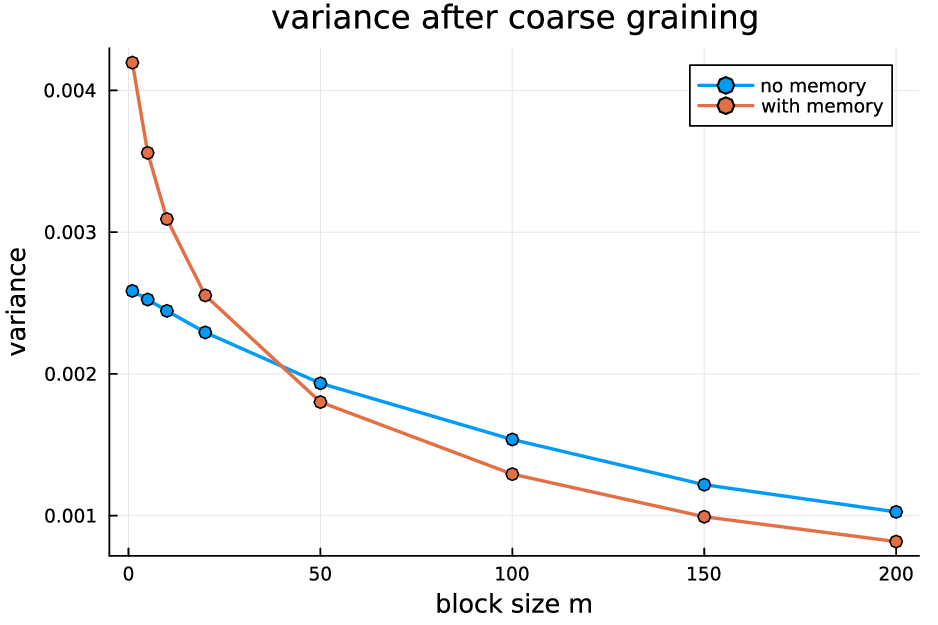}
  \caption{{Variance after coarse graining.}
  The memory case begins with a larger variance but declines more rapidly with block size and crosses below the classical curve after modest averaging, consistent with the weighted variance formula and the short lag behavior in Fig.~\ref{fig:acf}.}
  \label{fig:varcg}
\end{figure}

Figure~\ref{fig:varcg} walks through the consequence. At native resolution ($m=1$) the memory case has the larger variance because its high frequency content passes the filter. As $m$ grows, averaging is initially more effective for the fractional series since its very short lag correlations are smaller (Fig.~\ref{fig:acf}); the variance for $\alpha=0.7$ drops faster and crosses below the $\alpha=1$ curve around the reported crossover scale.

The slow tail that remains after short lags have been averaged out is best captured by $\tau_{\mathrm{int}}$. Figure~\ref{fig:persist} shows that $\tau_{\mathrm{int}}(W)$ increases with block size in both models, but once coarse averaging removes the fastest components, the fractional curve becomes larger. At coarse scales, the process with memory is therefore the more persistent one.

\subsubsection*{Numerical implementation}
We simulate $\mathcal{D}^{\alpha}x=f(x)+\sigma\,\xi(t)$ with $x(0)=1$, total time $T=1000$, step $h=0.01$, and $\sigma=0.1$. The same noise realization is used for $\alpha=1$ and $\alpha=0.7$ to isolate dynamical filtering. After discarding the first $100$ time units, we compute variances, autocorrelations, and $\tau_{\mathrm{int}}(W)$. Coarse graining uses $m\in\{1,5,10,20,50,100,150,200\}$. Power spectra are estimated by a Welch periodogram with a Hann window~\cite{Welch1967}, segment length $65536$, and $75\%$ overlap. One-sided frequency plots stop at the Nyquist limit $f_N=1/(2h)$.

\begin{figure}[ht!]
  \centering
  \includegraphics[width=0.78\linewidth]{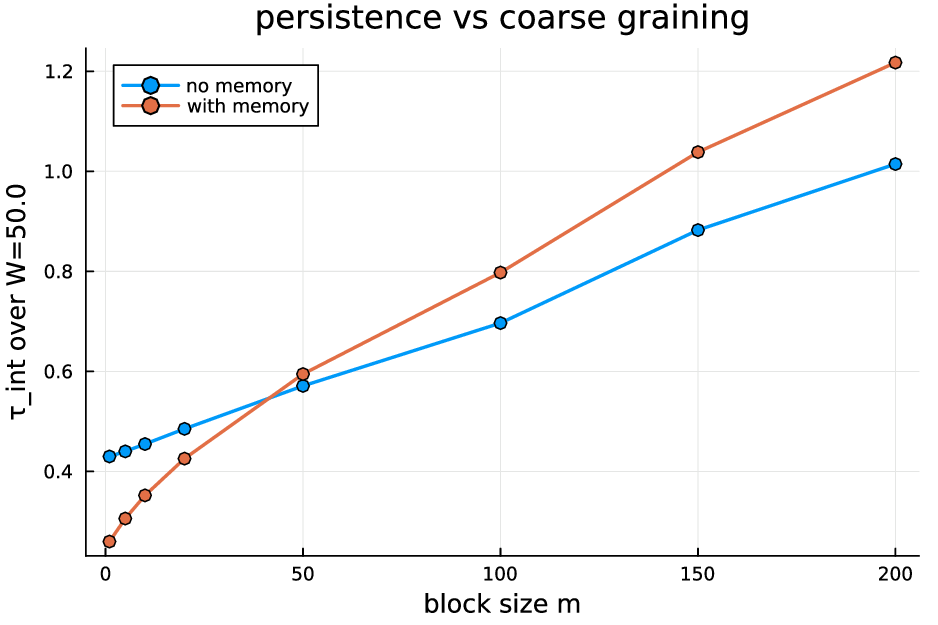}
  \caption{{Persistence versus coarse graining.}
  The windowed integrated autocorrelation time grows with block size for both models. After modest averaging the fractional curve exceeds the classical one, revealing greater persistence at coarse scales.}
  \label{fig:persist}
\end{figure}

With the same external noise level, a fractional response with $\alpha<1$ attenuates fast components less. The immediate outcome is a larger raw variance and lower short delay similarity at native resolution (Figs.~\ref{fig:traj} and \ref{fig:acf}). After modest averaging that removes the fastest content, the slow correlation tail of the fractional model dominates, so the memory case becomes more persistent at coarse scales (Figs.~\ref{fig:varcg} and \ref{fig:persist}). These observations are exactly consistent with the spectral picture in Fig.~\ref{fig:psd} and support the statement in the main text that the effect of memory depends on the observational scale.

\FloatBarrier
\subsection*{Switching and dwell times}

\subsubsection*{Definitions and detection of \emph{committed} switches}
Let the unstable saddle be at $x=0$ and the two basin cores be
\[
\mathcal{C}_+ := \{x:\; x \ge x_{\mathrm{core}}\}, 
\qquad 
\mathcal{C}_- := \{x:\; x \le -x_{\mathrm{core}}\},
\]
with $x_{\mathrm{core}}\in(0,1)$. 
Define the core label $\ell(t)\in\{-1,0,+1\}$ by 
$\ell(t)=+1$ if $x(t)\in\mathcal{C}_+$, $\ell(t)=-1$ if $x(t)\in\mathcal{C}_-$, and $\ell(t)=0$ otherwise.
With sampling step $h$, a \emph{committed core entry} at time $t_k$ occurs when there exists $s\in\{+1,-1\}$ such that
\[
\ell(t_k)=s\quad\text{and}\quad \ell(t_k+jh)=s\ \ \text{for all}\ \ j=1,\dots,\texttt{hold\_steps},
\]
i.e., the trajectory has entered a basin core and remained there for a minimum hold time 
$\tau_{\mathrm{hold}}=\texttt{hold\_steps}\cdot h$.  
To avoid counting re-entries to the \emph{same} core, we keep only \emph{alternating} committed entries $t_1<t_2<\cdots<t_m$ whose labels flip sign ($+1\to-1$ or $-1\to+1$). 
The \emph{dwell times} (residency times) are then the inter-event intervals
\[
T_i := t_{i+1}-t_i,\qquad i=1,\dots,m-1.
\]
This definition suppresses the “micro–switches” created by fast excursions around the separatrix, which are especially frequent in the fractional case. 
In all figures below, we use $(x_{\mathrm{core}},\tau_{\mathrm{hold}})=(0.9,\,0.2)$ unless stated otherwise.

We simulate the double-well with additive noise
\[
\mathcal{D}^{\alpha}x(t)=x(t)-x(t)^3+\sigma\,\xi(t),\qquad 0<\alpha\le 1,
\]
over a long horizon $T=5000$ with step $h=0.01$ and noise level $\sigma=0.34$. 
The \emph{same} discrete noise path is used for $\alpha=1$ (“no memory”) and $\alpha=0.7$ (“with memory”) to isolate dynamical filtering. 
We mark alternating committed entries and compute the dwell sequence $\{T_i\}$ for each model.

\subsubsection*{Time–series view with detected switches}
\begin{figure}[ht!]
  \centering
  \includegraphics[width=\linewidth]{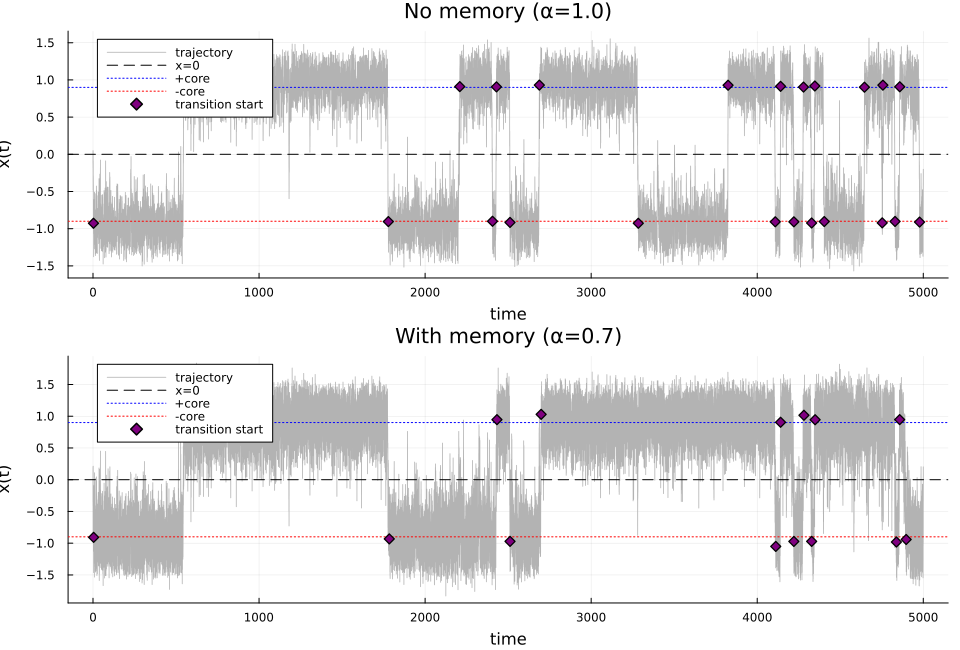}
  \caption{{Trajectories with committed switch detection.}
  Gray: $x(t)$; dotted: $\pm x_{\mathrm{core}}$; diamonds: detected starts of \emph{alternating committed} transitions. 
  Top: $\alpha=1$, with 23 committed transitions; bottom: $\alpha=0.7$, with 15 committed transitions.}
  \label{fig:trajswitch}
\end{figure}

Before turning to distributions, we show how the rule above behaves on the raw trajectories. 
In Fig.~\ref{fig:trajswitch} (two panels) the gray trace is $x(t)$, dotted lines mark $\pm x_{\mathrm{core}}$, and purple diamonds mark the detected \emph{starts} of committed transitions.

For this shared-noise run, the memory-free process records more committed switches (23 vs.\ 15), while the fractional process shows extended residence episodes interspersed with clusters of short jitters near the saddle that do \emph{not} pass the core/hold filter. 
This illustrates why raw zero-crossing counts would overstate switching under memory: many near-threshold rattles are not true regime changes.

\subsubsection*{From dwell lists to distributions}
Given the dwell sample $\{T_i\}_{i=1}^n$, the empirical cumulative distribution function (CDF) and survival (complementary CDF) are
\[
F_n(t)=\frac{1}{n}\sum_{i=1}^n \mathbf{1}\{T_i\le t\}, 
\qquad 
S_n(t)=1-F_n(t)=\frac{1}{n}\sum_{i=1}^n \mathbf{1}\{T_i>t\}.
\]
The (instantaneous) \emph{hazard} is $h_n(t)=f_n(t)/S_n(t)$ with $f_n=F_n'$, estimated nonparametrically by smoothing $F_n$ or via kernel density divided by $S_n$. 
For a memory-free (Markov) escape with constant rate $r$, one expects $S(t)\approx e^{-rt}$ (straight line on a semi-log plot) and $h(t)\equiv r$. 
In contrast, fractional dynamics create \emph{age-dependent} escape: $S(t)$ deviates from a simple exponential and $h(t)$ typically decreases with $t$ (e.g., stretched-exponential or Mittag–Leffler-like survival). 
Because such tails are sensitive to the observation window and censoring, we report survival curves alongside histograms instead of a single “escape rate.”

\begin{figure}[th!]
  \centering
  \includegraphics[width=0.78\linewidth]{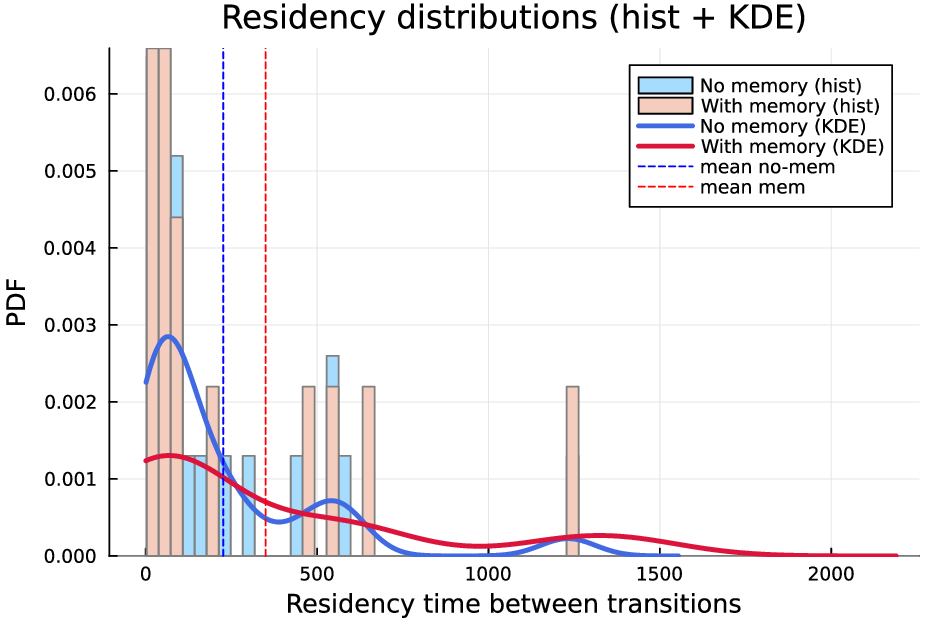}
  \caption{{Dwell-time distributions (histogram + KDE).}
  Transparent bars: equal-width histograms (normalized to a probability density function). 
  Solid lines: KDEs truncated at $t\ge 0$. 
  Vertical dashed lines mark sample means. 
  The memory case shows a heavier tail.}
  \label{fig:residhist}
\end{figure}

Figure~\ref{fig:residhist} illustrates histograms of dwell time distribution with nonparametric kernel density estimations (KDE).  
Both models produce many short dwells (left modes), but the memory case exhibits a heavier right tail: the mean dwell time is larger despite a smaller median (representative run: median $79$ vs.\ $94$, mean $350$ vs.\ $226$ for memory vs.\ no-memory). 
This “shorter typical, longer occasional” pattern is characteristic of broad, skewed dwell distributions.

\begin{figure}[th!]
  \centering
  \includegraphics[width=0.7\linewidth]{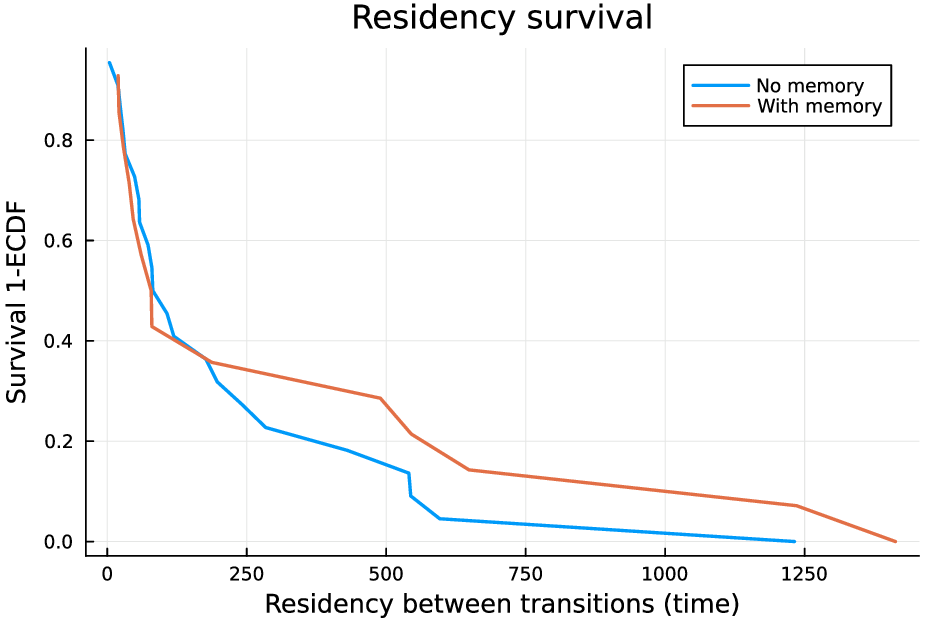}
  \caption{{Dwell-time survival.}
  The memory curve decays more slowly, indicating a broader distribution with more probability mass at long times.}
  \label{fig:survival}
\end{figure}

Figure~\ref{fig:survival} plots $S_n(t)=1-F_n(t)$. 
Over a wide time range the survival curve for the memory model lies above the classical curve, and the decay slows perceptibly in the tail (compare the slope beyond $t\approx 300$). 
This is the distributional signature of age-dependent escape: the longer the system has stayed in a basin, the less likely it is to leave in the next instant.

\subsubsection*{Mathematical support and expectations}

To picture the local barrier-crossing,
near the stable point $x^\star\in\{\pm1\}$, linearization yields the fractional OU surrogate 
$\mathcal{D}^{\alpha} y + \lambda y = \sigma\xi(t)$ with $\lambda=2$. 
For $\alpha=1$ and small $\sigma$, Kramers theory~\cite{hanggi1990reaction,berglund2009e} gives an approximately constant escape rate $\displaystyle{r\propto \sqrt{|f'(0)|\,f'(x^\star)}\,e^{-\Delta V/\sigma^2}}$, hence exponential dwell times. 
When $0<\alpha<1$, the nonlocal operator introduces memory in the relaxation trajectory; the exit problem becomes non-Markovian. 
Fractional escape problems (and, more generally, subordinated diffusions) generically produce survival laws between stretched-exponential and Mittag–Leffler forms~\cite{METZLER2000238}, 
\[
S(t)\approx E_\alpha\!\left[-\Big(\tfrac{t}{\tau}\Big)^\alpha\right]
\quad\text{or}\quad 
S(t)\approx \exp\!\left[-(t/\tau)^\beta\right],\ \ 0<\beta<1,
\]
whose hazards decrease with $t$. 
Thus one expects a broadened dwell distribution, stronger sensitivity to the observation window, and the possibility that means diverge in the idealized limit (in practice bounded by finite $T$). 
Our empirical survival curves (Fig.~\ref{fig:survival}) exhibit a slowing tail.

We found that using “committed” detection is essential here. The fractional model admits rapid near-threshold \emph{rattling} caused by its weaker high-frequency attenuation. 
If one were to count every sign change of $x(t)$, these micro-excursions would inflate the apparent switching rate for $\alpha<1$. 
Requiring entry into $\mathcal{C}_\pm$ \emph{and} a hold time $\tau_{\mathrm{hold}}$ creates a coarse-grained, behaviorally meaningful notion of regime change that is robust to such jitter. 
Figures~\ref{fig:trajswitch}–\ref{fig:survival} are computed with this committed rule.

\subsubsection*{Censoring, window effects, and robustness}
Dwell-time data from a finite window $[0,T]$ are subject to left/right censoring: the run that straddles $t=0$ and the terminal partial run at $T$ are excluded when we use differences of detected \emph{starts} ($T_i=t_{i+1}-t_i$). 
This aligns the two models and avoids biasing the survival tail with an artificially censored last interval. 
We verified that the qualitative comparisons are unchanged across a range of $(x_{\mathrm{core}},\tau_{\mathrm{hold}})$ and under alternative binning or kernel bandwidths (not shown). 
Because the two simulations share the same noise path, differences in Figs.~\ref{fig:residhist}–\ref{fig:survival} are attributable to the dynamical filtering by $\alpha$.

\subsubsection*{Numerical implementation}
For the run displayed in Fig.~\ref{fig:trajswitch} (shared noise, $\sigma=0.34$, $T=5000$, $h=0.01$, $x_{\mathrm{core}}=0.9$, $\tau_{\mathrm{hold}}=0.2$), we observed 23 transitions in the no-memory case compared with 15 in the memory case. The corresponding dwell counts were 22 and 14, with median dwell times of 94.36 and 79.35, respectively. The mean dwell times were 226.12 for the no-memory system and 349.52 for the system with memory.
These numbers will vary with the seed and parameter choices, but the qualitative pattern is consistent: memory broadens the dwell distribution and slows the survival tail.

Memory changes switching in two complementary ways. 
At fine resolution, it generates more near-threshold activity, which inflates na\"ive zero-crossing counts; the committed rule avoids this artefact. 
At the level of true regime changes, memory broadens the dwell-time distribution and makes the survival tail decay more slowly, so that “typical” dwells (medians) can be shorter while rare, very long dwells pull the mean upward. 
Because these effects depend on basin curvature and on the timing of noise bursts relative to the memory kernel, the net switching \emph{frequency} can increase or decrease across parameter regimes—hence our focus on full dwell distributions and survival curves rather than a single escape rate.

\FloatBarrier
\subsection*{State dependent noise: multiplicative case}

We consider the double well drift $f(x)=x-x^{3}$ with state dependent noise amplitude $g(x)$ and memory parameter $\alpha\in(0,1]$.
When the noise amplitude depends on the state ($g=g(X_t)$), the stochastic integral $\int g(X_t)\,\mathrm{d}W_t$ is ambiguous unless we choose a convention. Two standard choices are:

- \textbf{Itô} (default here)~\cite{KloedenPlaten1992}: evaluate $g$ at the \emph{left} endpoint of each time step. This is the convention we use unless stated otherwise.

- \textbf{Stratonovich} (denoted $\circ$)~\cite{KloedenPlaten1992}: evaluate $g$ at the \emph{midpoint} in time (equivalently, the limit of smooth noise). 

These two give different drifts but the same diffusion. Specifically,
\[
\mathrm{d}X_t = f(X_t)\,\mathrm{d}t + g(X_t)\,\circ\,\mathrm{d}W_t
\quad\Longleftrightarrow\quad
\mathrm{d}X_t = \big[f(X_t)+\tfrac12 g(X_t)g'(X_t)\big]\mathrm{d}t + g(X_t)\,\mathrm{d}W_t.
\]
Thus the Stratonovich model is identical to an Itô model with an \emph{extra deterministic drift} $\tfrac12 gg'$.

For the classical case $\alpha=1$, we write
\begin{equation}
\mathrm{d}X_t = f(X_t)\,\mathrm{d}t + g(X_t)\,\mathrm{d}W_t,
\label{eq:M1}
\end{equation}
interpreted as Itô unless stated otherwise. In discrete time with step $h$, the Itô Euler–Maruyama update~\cite{KloedenPlaten1992} is
\[
X_{k+1}=X_k+ f(X_k)\,h + g(X_k)\,\sqrt{h}\,\eta_k,\qquad \eta_k\sim\mathcal{N}(0,1).
\]
Here, left endpoint evaluation $\Rightarrow$ Itô type discretization, that is we plug in $X_k$ (the start of the step) into $f$ and $g$. 
If instead we want the Stratonovich model, a midpoint or Heun scheme~\cite{KloedenPlaten1992} is appropriate, and in Itô form it adds the well known correction $\tfrac12 gg'$ to the drift.

For the fractional case $\alpha<1$, we use a Caputo fractional drift driven by white in time forcing,
\begin{equation}
{}^{\mathrm{C}}D_t^{\alpha} X_t = f(X_t) + g(X_t)\,\xi(t),
\label{eq:M2}
\end{equation}
with $0<\alpha<1$. Numerically, we represent the white process as piecewise constant on $[t_k,t_{k+1})$,
\[
\xi(t)=\frac{\eta_k}{\sqrt{h}}\quad \text{for } t\in[t_k,t_{k+1}), \qquad \eta_k\sim\mathcal{N}(0,1),
\]
and evaluate $g$ at the left endpoint $X_k$. This produces the same Itô type treatment of the multiplicative term as in the classical case; the memory enters through the fractional drift (via the Caputo history sum), not through the noise interpretation.

Inside a potential well (near a stable point), the deterministic pull is set by the local slope of $f$. In Stratonovich form that pull is effectively modified by $\tfrac12 gg'$: it can shift the location of equilibria and change the apparent curvature of the well, which in turn affects residence times and barrier crossing rates. This effect only appears when $g$ depends on $x$.

In examples, we take $g(x)=g_0\,\big(1-\beta x^2\big)_+$, which concentrates noise near the saddle $x\approx 0$ and reduces it inside the wells. Here $(z)_+=\max\{z,0\}$, so $g(x)=0$ once $|x|>1/\sqrt{\beta}$.

\subsubsection*{Local linear analysis and fine-scale signature}
Let \(x_\star\in\{\pm1\}\) with \(f'(x_\star)=-2\). For small excursions \(X_t=x_\star+Y_t\),
\begin{equation}
\mathcal{D}^{\alpha}Y_t \approx -2Y_t + g_\star\,\xi(t),\qquad g_\star:=g(x_\star).
\label{eq:M4}
\end{equation}
Thus, locally we obtain a fractional OU surrogate: decreasing \(\alpha\) lengthens correlation (larger windowed \(\tau_{\mathrm{int}}\)) and attenuates fast frequencies less, exactly as in the additive case, with \(\sigma\) replaced by \(g_\star\).

\begin{figure}[th!]
  \centering
  \includegraphics[width=0.78\linewidth]{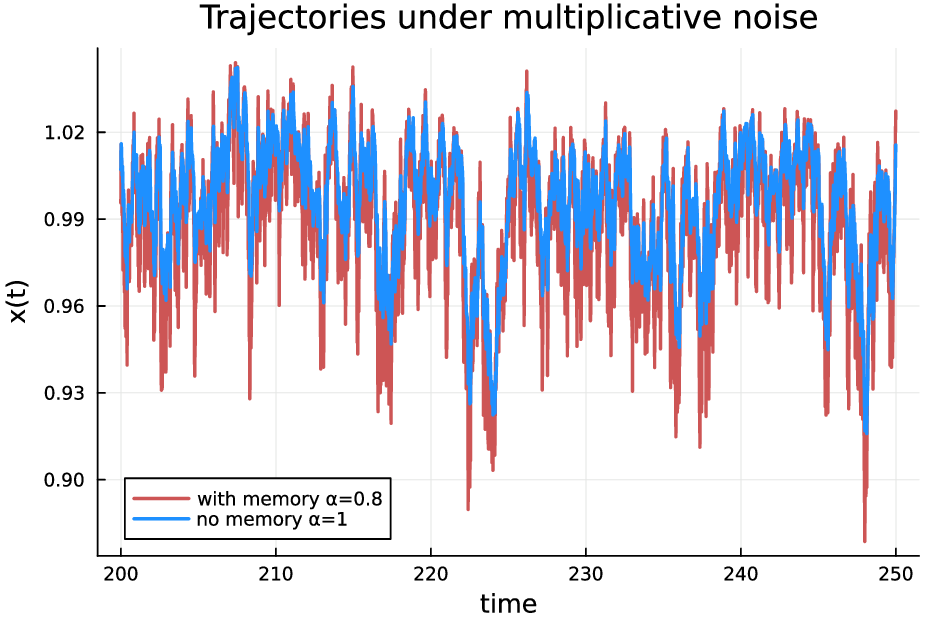}
  \caption{{Trajectories under multiplicative noise (small segment).}
  With the same noise path and \(g(x)\), the memory case (\(\alpha=0.8\)) is rougher at fine scales than \(\alpha=1\).}
  \label{fig:multiseg}
\end{figure}

A short window at native resolution under the same multiplicative-noise path shows that the memory trajectory (\(\alpha=0.8\)) exhibits visibly sharper, faster wiggles than \(\alpha=1\) (Figure \ref{fig:multiseg}).  
This is the time-domain footprint of the spectral flattening implied by \eqref{eq:M4}: memory still transmits more fast content even when the amplitude depends on state.

\subsubsection*{Global variance and occupancy}

When the noise amplitude $g(x)$ varies across the state space, the overall variance of $X_t$ mixes contributions from regions with different noise levels. To make this explicit, partition the state space into disjoint sets $B_1,\dots,B_J$ (for example, the two basins and a small saddle neighborhood, or bins in $|x|$). Define a coarse state indicator
\[
Z(t)=j \quad \text{iff} \quad X_t \in B_j.
\]
The occupancy of region $B_j$ is the probability $\mathbb{P}(Z=j)$, i.e. the long run fraction of time the trajectory spends in $B_j$.
For a time series sampled at $t_n=nh$, the empirical occupancy is
\[
\widehat{\mathbb{P}}(Z=j)=\frac{1}{N}\sum_{n=1}^{N}\mathbf{1}\{X_{t_n}\in B_j\},
\]
and the corresponding conditional moments are
\[
\widehat{\mathbb{E}}[X\mid Z=j]=\frac{1}{N_j}\sum_{n:\,X_{t_n}\in B_j}X_{t_n}, 
\qquad
\widehat{\operatorname{Var}}(X\mid Z=j)=\frac{1}{N_j-1}\sum_{n:\,X_{t_n}\in B_j}\bigl(X_{t_n}-\widehat{\mathbb{E}}[X\mid Z=j]\bigr)^2,
\]
with $N_j=\sum_n \mathbf{1}\{X_{t_n}\in B_j\}$.

By the law of total variance,
\begin{equation}
\operatorname{Var}(X)=\mathbb{E}\!\big[\operatorname{Var}(X\mid Z)\big]+\operatorname{Var}\!\big(\mathbb{E}[X\mid Z]\big)
=\sum_{j=1}^{J}\mathbb{P}(Z=j)\,\operatorname{Var}(X\mid Z=j)\;+\;\sum_{j=1}^{J}\mathbb{P}(Z=j)\,\bigl(\mu_j-\mu\bigr)^2,
\label{eq:M5}
\end{equation}
where $\mu_j=\mathbb{E}[X\mid Z=j]$ and $\mu=\mathbb{E}[X]$. The first term is the average of the \emph{within region} variances; the second is the \emph{between regions} variance due to different local means.

Memory modifies two ingredients:
\begin{itemize}
\renewcommand\labelitemi{--}
\item \emph{Local variability} $\operatorname{Var}(X\mid Z=j)$, through the short time roughness discussed earlier. For multiplicative noise, regions with larger $g(x)$ tend to show larger micro fluctuations.
\item \emph{Occupancy} $\mathbb{P}(Z=j)$, by changing residence times and transition rates between regions. If memory makes the process linger where $g$ is large, the weights on high variance regions increase.
\end{itemize}
Therefore, the global sample variance does not need to be monotone in $\alpha$: it can rise or fall depending on how the time spent across regions reweights the local contributions, without contradicting the local mechanism.

Regarding coarse graining and block size $m$,
let $X_k^{(m)}=\frac{1}{m}\sum_{i=1}^{m}X_{(k-1)m+i}$ be block averages over $m$ samples. The same decomposition holds for $X^{(m)}$ if we define $Z^{(m)}$ by the region of the block (for example, the modal region within the block). Coarse graining suppresses high frequency variance inside each region but leaves the occupancy weights $\mathbb{P}(Z)$ unchanged, so differences driven by time spent in high $g(x)$ zones can persist across $m$.

\begin{figure}[th!]
  \centering
  \includegraphics[width=0.78\linewidth]{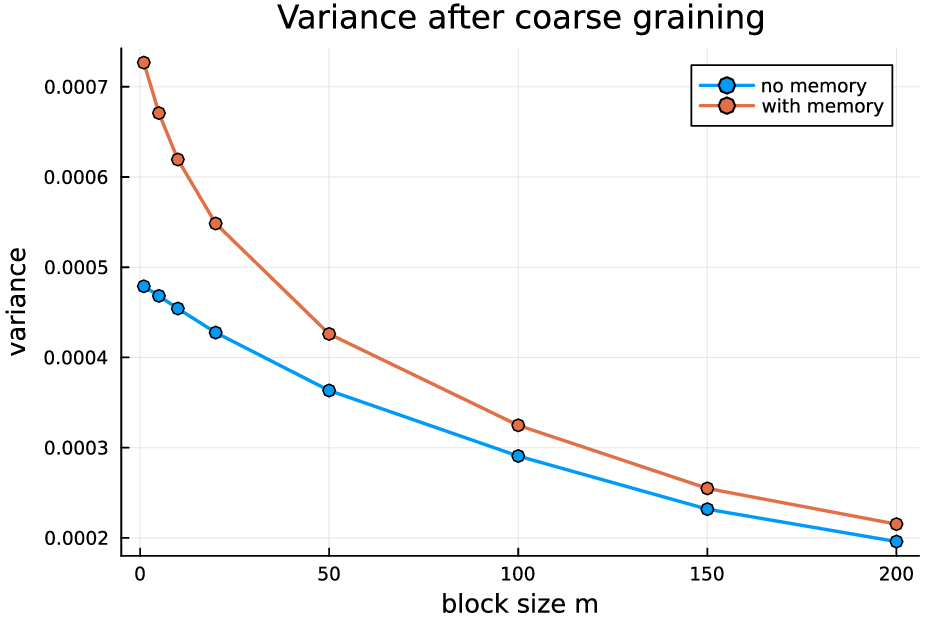}
  \caption{Variance after coarse graining (multiplicative).
  Memory (orange) remains above no-memory (blue) for all block sizes in this example, consistent with increased occupancy of high-\(g(x)\) regions.}
  \label{fig:varcg-mult}
\end{figure}

In the shown run, Figure~\ref{fig:varcg-mult}, $g(x)$ is strongest near the saddle, and the fractional trajectory spends relatively more time in that high noise region, so $\mathbb{P}(Z=\text{saddle})$ increases. Both the within region term and the occupancy weights push the memory curve above the memory-free curve, and this ordering remains after coarse graining. With a different $g(\cdot)$ (for example, larger noise deep inside the wells) the occupancy pattern can flip and so can the ordering, as \eqref{eq:M5} predicts.

\smallskip
\noindent\textit{Two basin special case (for intuition):}
If $Z\in\{\text{L},\text{R}\}$ are the two wells and we ignore the saddle, then
\[
\operatorname{Var}(X)=p\,\operatorname{Var}(X\mid \text{L})+(1-p)\,\operatorname{Var}(X\mid \text{R})
+p(1-p)\bigl(\mu_{\text{L}}-\mu_{\text{R}}\bigr)^{2},\quad p=\mathbb{P}(Z=\text{L}).
\]
Memory can change both the within well variances and the mixing weight $p$, so the global variance can move in either direction with $\alpha$.

\subsubsection*{Barrier crossing under multiplicative noise}
We detect \emph{committed} switches using the same core/hold rule as in the additive section (cores at \(|x|\ge x_{\mathrm{core}}=0.9\), hold time \(\tau_{\mathrm{hold}}=0.2\)).  
When \(g(x)\) peaks near the saddle, memory can produce clusters of passages interleaved with long residence episodes: near-threshold jitter is more frequent, but not every rattle becomes a committed switch.

\begin{figure}[th!]
  \centering
  \includegraphics[width=0.95\linewidth]{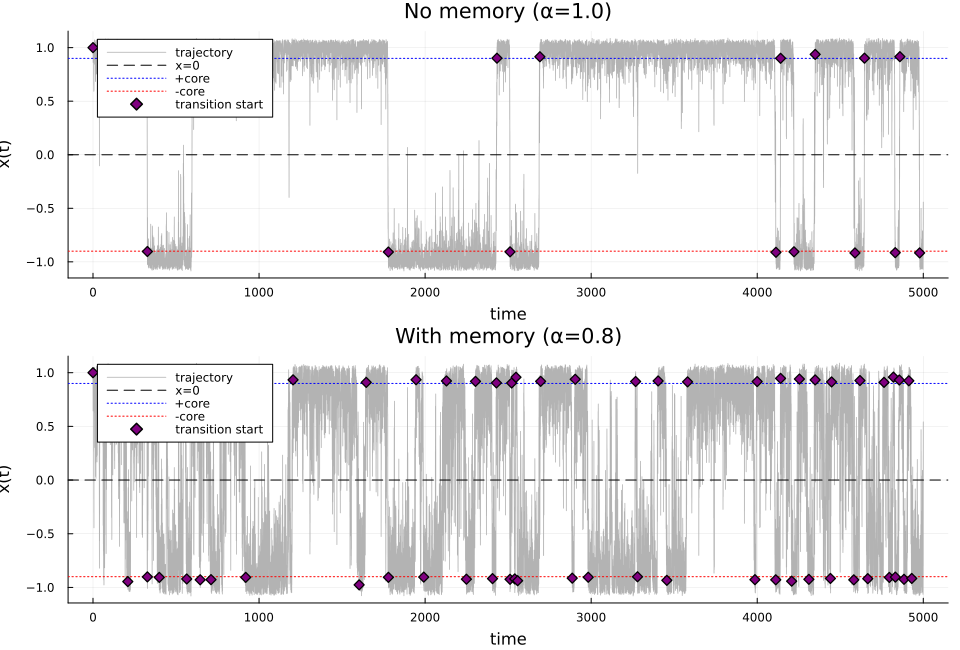}
  \caption{\textbf{Committed transition starts (multiplicative).}
  Gray: \(x(t)\); dotted: \(\pm x_{\mathrm{core}}\); diamonds: starts of alternating committed switches. 
  Top: \(\alpha=1\); bottom: \(\alpha=0.8\).}
  \label{fig:trans-mult}
\end{figure}

In a representative long run (Figure \ref{fig:trans-mult}), the memory case (\(\alpha=0.8\)) registers many more committed switches (59 vs.\ 15) and a longer mean committed dwell (about \(1103\) vs.\ \(674\) steps).  
This coexistence of clustered crossings and heavy-tailed dwells is typical when state-dependent forcing interacts with memory.

\subsubsection*{Numerical details}
Short segments use \(T=1000\), \(h=0.01\), \(x(0)=1\), base multiplicative level \(\sigma=0.2\); switching runs use \(T=5000\) and \(\sigma=0.5\).  
All panels share the same pre-sampled noise path across \(\alpha\); committed detection uses \((x_{\mathrm{core}},\tau_{\mathrm{hold}})=(0.9,0.2)\).  
Halving \(h\) leaves the qualitative comparisons unchanged.  
Unless noted, we adopt the It\^{o}-type discretization (left-endpoint \(g\)); a Stratonovich-like midpoint would add the noise-induced drift \(\tfrac12 gg'\).

Locally, decreasing \(\alpha\) still lengthens correlation and lets more fast content through (Fig.~\ref{fig:multiseg}).  
Globally, variance reflects both within-region variability and region occupancy (mixture identity \eqref{eq:M5}); thus the ordering can change with \(g(x)\) (Fig.~\ref{fig:varcg-mult}). 
Near the saddle, memory plus state-dependent noise yields clustered passages and slow-survival dwell tails (Fig.~\ref{fig:trans-mult}).  

\section*{S5: Fitting a memory-free model to data exhibiting memory}\label{sec: sup5}
Ignoring memory in the formulation of dynamical models can lead to wrong inferences about stability and critical transitions. Here, we quantify those errors in a controlled thought experiment. Synthetic time-series data are generated with a fractional-order version of a classical herbivory model. We then fit the same data with an otherwise identical, memory-free (integer-order) model. Although the memory-free fit reproduces the trajectories, it systematically misplaces tipping points, distorts estimates of resilience and resistance, and reconstructs a static stability landscape where the true landscape is history-dependent.

\subsection*{Numerical implementation}
To generate the data, we use the fractional herbivory model
\begin{equation}\label{eq: frac_herbivor}
\mathcal{D}^{\alpha}x=rx\!\left(1-\frac{x}{K}\right)-\frac{Bx}{A+x}
\end{equation}
where $x$ is vegetation biomass. The parameters are set as follows: intrinsic vegetation growth rate $r=0.8$,
carrying capacity $K=6$, half saturation constant of consumption $A=0.2$. We explore several values of grazing intensity $B=\{0.2, 0.3,..., 1.1, 1.2\}$, and use a fractional derivative order $\alpha=0.8$

The goal is to fit the generated fractional-order data using the integer-order version of the model ($\alpha=1$). In this fitting process, we aim to estimate the parameters while knowing that $B$ is the varying factor. We use a bound-constrained nonlinear least-squares fit using MATLAB’s \texttt{lsqnonlin} (trust-region-reflective because bounds are set). The objective stacks residuals across all runs and both initial conditions, which comes from integrating the ODE with ode45 and evaluating at data times.

\subsubsection*{Memory-dependent feature challenges}

Although fitting a noise-free ODE may appear to be a straightforward task, it is not trivial. Achieving a good fit is challenging because the memory-free model can easily overfit or misrepresent key system features. 
In a noise-free memory-free system, trajectories are unique: once the state lies on a given trajectory, the subsequent evolution is fixed. By contrast, memory endows the dynamics with a dependence on history, so different initial values generate distinct trajectories even when they momentarily share the same state. This distinction matters when assessing how the system reacts to disturbances. and studying stability landscapes. In a memory-free model, the reconstructed stability landscape is independent of where the system begins; trajectories and the resulting landscape overlay are identical no matter the starting point. By contrast, memory ``remembers'' any past trajectory, so the same model parameters can yield different landscapes depending on the initial state.
Without memory, every run collapses onto a single trajectory and landscape; with memory, different starting values produce distinct trajectory shapes and altered landscapes (see e.g. Figure~\ref{fig: initial conditions}). 


\subsubsection*{Parameter fitting under memory constraints}
\begin{figure}[th!]
  \centering
  \includegraphics[width=\linewidth]{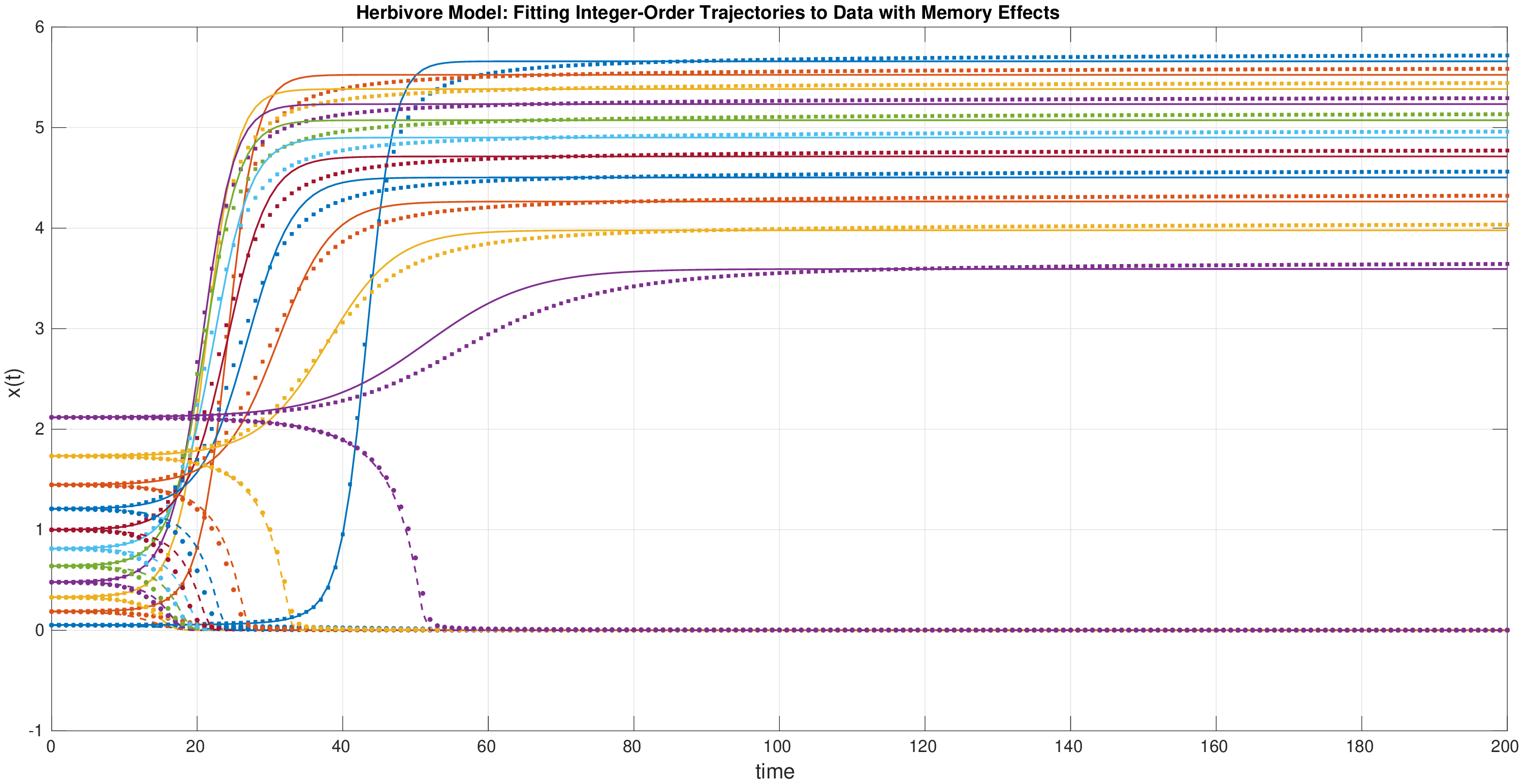}
\caption{
The dotted curves represent trajectories of the herbivory system with memory~\eqref{eq: frac_herbivor}, simulated using parameters 
\(r = 0.8\), \(K = 6\), \(A = 0.2\), \(B = \{0.2,\, 0.3,\, \ldots,\, 1.1,\, 1.2\}\), and fractional order \(\alpha = 0.8\) (memory \(= 0.2\)). 
Initial conditions are chosen near the unstable equilibrium points.  
The solid curves correspond to the fitted integer-order model with estimated parameters 
\(r = 0.5726\), \(K = 5.8242\), and \(A = 0.1112\), and herbivore pressures 
\(B = \{0.0932,\, 0.1656,\, 0.2379,\, 0.3102,\, 0.3825,\, 0.4548,\, 0.5270,\, 0.5990,\, 0.6707,\, 0.7420,\, 0.8124\}\).
}
  \label{fig: fit-data}
\end{figure}

Fitting a memory-driven model with an otherwise identical memory-free one is unreliable, but we can mitigate the problem by enforcing a common starting point: every synthetic trajectory begins near the unstable equilibrium, the hilltop that divides the two basins of attraction, and we restrict parameters to the bistable regime. Fixing this initial state forces the deterministic paths, and therefore the reconstructed stability landscapes, to coincide in both models; working outside the bistable range would eliminate a well-defined separatrix and make both data generation and parameter fitting ill-posed. Under these constraints, we obtain a satisfactory fit, as illustrated in Figure~\ref{fig: fit-data}.
The fitted parameters are \(r = 0.5726\), \(K = 5.8242\), and \(A = 0.1112\).  
The herbivore pressure values are \(B = \{0.0932,\, 0.1656,\, 0.2379,\, 0.3102,\, 0.3825,\, 0.4548,\, 0.5270,\, 0.5990,\, 0.6707,\, 0.7420,\, 0.8124\}\).

\subsection*{Inference failures}
\subsubsection*{Critical-point misestimation when neglecting memory}

When a memory-driven system is fitted using a memory-free model, the numerical fit may appear excellent, yet the resulting bifurcation diagram is systematically distorted. The predicted equilibrium branches, and consequently the tipping points, are displaced along the control-parameter axis. 

In the herbivory model, where the grazing rate \(B\) serves as the bifurcation parameter, we calibrate the memory-free model to data generated by the true (memory-driven) system by varying \(B\) while keeping all other parameters fixed. The fitted model reproduces the qualitative dynamics, including stable and unstable equilibria and the trajectories converging toward them. However, their locations are shifted: e.g., the true system transitions at \(B \approx 1.2\), while the memory-free model predicts the same transition near \(B \approx 0.8\).

\begin{figure}[th!]
  \centering
  \includegraphics[width=\linewidth]{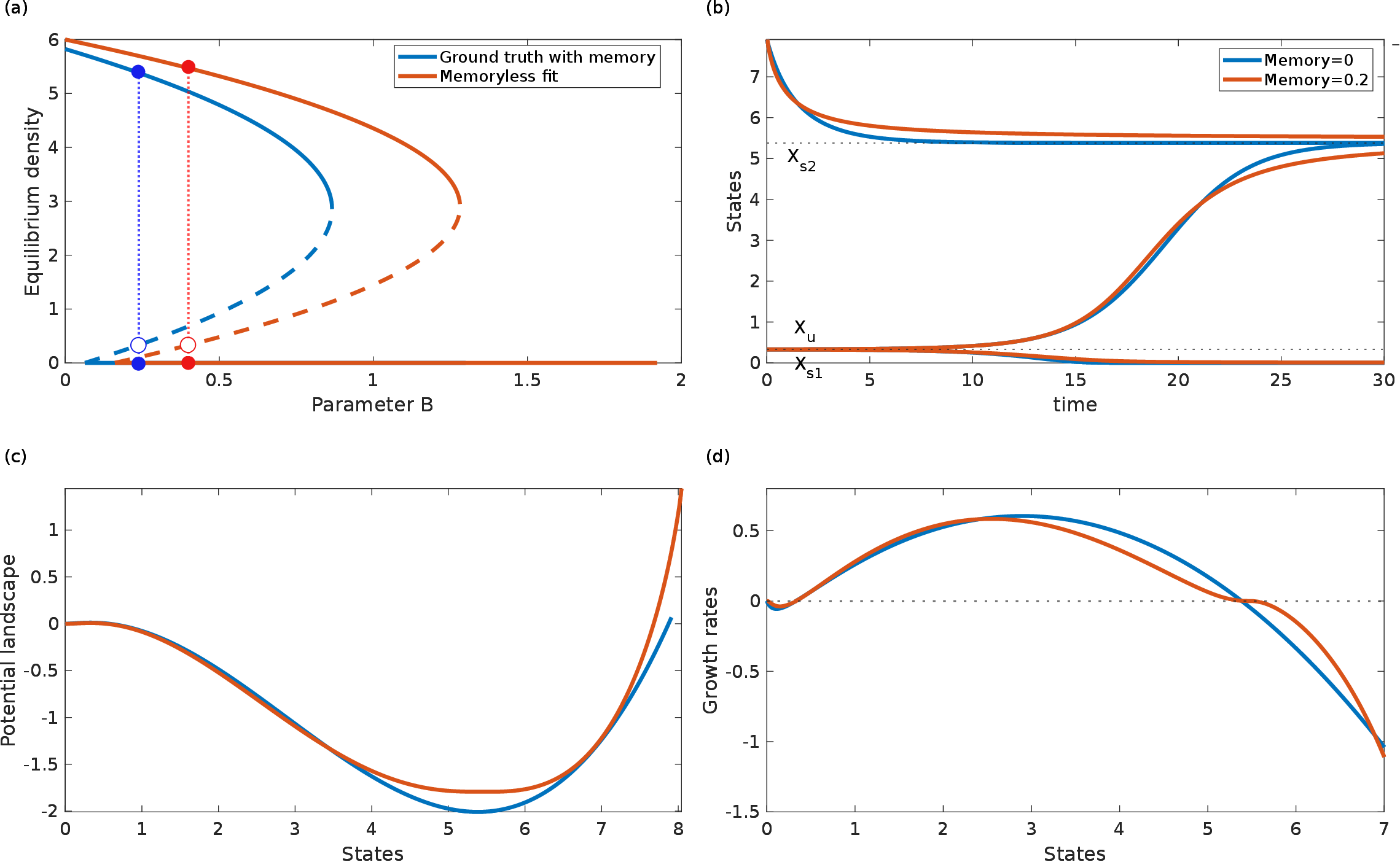}
  \caption{Bifurcations, trajectories, potential landscapes, and potential slopes of the original herbivore model with memory (\(r = 0.8\), \(K = 6\), \(A = 0.2\), \(B =0.4\)), and \(\alpha = 0.8\)) and fitted model without memory (\(r = 0.5726\), \(K = 5.8242\), \(A = 0.1112\), and $B=0.2379$).}
  \label{fig:  land_fit_data}
\end{figure}

Figure~\ref{fig:  land_fit_data} illustrates this displacement: although the reconstructed potential landscapes and trajectories closely resemble those of the memory-driven system, the correspondence occurs at mismatched parameter values—for example, the memory-free model with \(B = 0.2379\) reproduces the behavior of the memory-driven system at \(B = 0.4\).

\begin{figure}[th!]
  \centering
  \includegraphics[width=\linewidth]{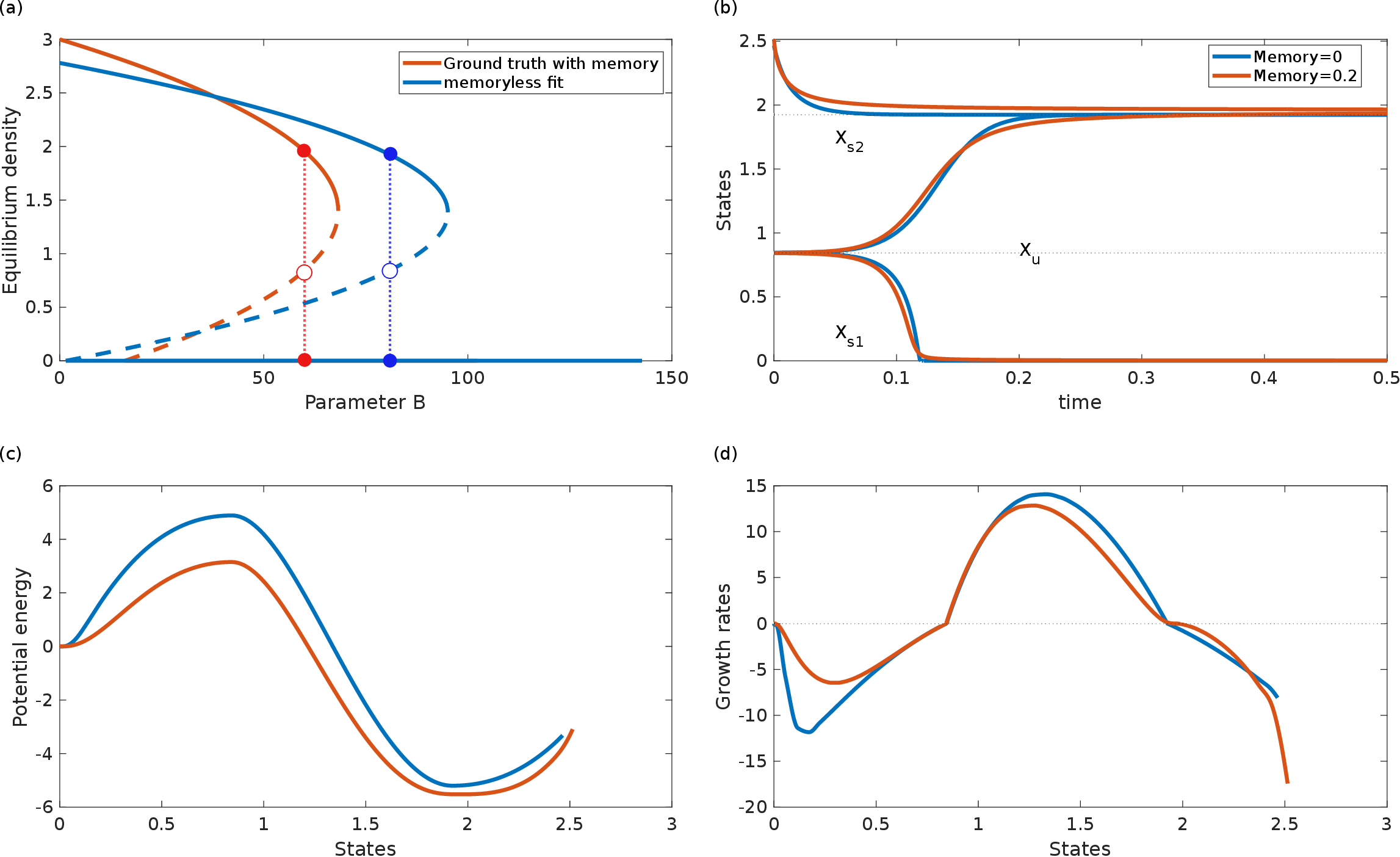}
  \caption{Bifurcations, trajectories, potential landscapes, and potential slopes of the original herbivore model with memory (\(r = 80\), \(K = 3\), \(A = 0.2\), \(B =60\)), and \(\alpha = 0.8\)) and fitted model without memory (\(r = 135.8\), \(K = 2.8\), \(A = 0.01\), and $B=80.89$).}
    \label{fig:  land_fit_data2}
\end{figure}

We repeated the fitting experiment under a different parameter setting to compare the outcomes, and similar misestimations and shifts were observed (see Figure~\ref{fig:  land_fit_data2}).

\subsubsection*{Possible misinterpretation}
As discussed earlier, memory makes the stability landscape dynamic in its response to perturbations. In other words, memory causes the landscape itself to drift, allowing the system to either rollback after a disturbance or collapse long after stress is removed—behaviors that would appear paradoxical in a static model. These effects arise not from noise, but from the delayed reshaping of the landscape: recovery or collapse depends on whether this drift completes before or after the system crosses the saddle. Consequently, post-disturbance analysis must account for slow variables and cumulative effects that reveal how memory reshapes stability.

However, when fitting data that contains memory effects, such unusual behavior is often misattributed to noise rather than to the system’s intrinsic dynamics. This misinterpretation can be problematic, particularly when both endogenous and exogenous perturbations are present—as in the fitting experiment conducted here.

\begin{figure}[th!]
  \centering
  \includegraphics[width=0.7\linewidth]{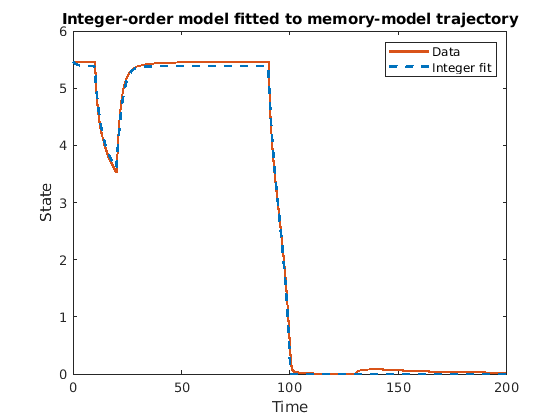}
\caption{Good trajectory fit, wrong parameters. An integer-order (memory-free) model (blue dashed) closely matches a trajectory generated by a memory system (orange solid) under two endogenous pulses, yet the inferred pulse magnitudes (0.6256 and 1.0117) differ from the true values (1 and 1.5), illustrating parameter bias from neglecting memory.}
  \label{fig:  fit data pulse}
\end{figure}

To demonstrate this, we applied two endogenous pulses: first by reducing the grazing parameter to 
\(B = B - 1\), and, after recovery, by applying a second pulse \(B = B - 1.5\). 
Using the previously fitted parameters, we obtained a good trajectory match 
(Figure~\ref{fig:  fit data pulse}), but the inferred parameter deviations were 
\(0.6256\) and \(1.0117\), respectively—noticeably different from the true values. 
This discrepancy is consistent with the shifted bifurcation diagram observed earlier 
(Figure~\ref{fig:  land_fit_data}).

In practice, such misinterpretations can have serious consequences: when control parameters are physically or policy relevant, overlooking memory-driven dynamics can lead to misguided decisions and ineffective management strategies.

To conclude, by contrasting fractional- and integer-order formulations of a system, we showed that (i) parameter estimation remains numerically accurate while qualitative features, location of equilibria, basin geometry, and critical parameter values, are displaced; (ii) the resulting potential landscape becomes time-dependent, invalidating static resilience and resistance metrics; and (iii) seemingly paradoxical behaviors such as slow rollback after pulses or delayed recovery or collapse after stress removal emerge naturally from the interaction between state and memory.

\section*{S6: Limitations, challenges, and future research directions}\label{sec: sup6}

\subsection*{Resilience metrics and threshold effects}
This study has several limitations related primarily to its definitions. For instance, resilience measurements in responding to endogenous perturbations can vary depending on the threshold chosen. Here, thresholds were selected on the basis of the system scale, defined as the distance between equilibrium points and the overall domain of system dynamics. Choosing a very small threshold means the system must exhibit high recovery capability, thereby providing a rigorous measure of resilience. In contrast, a larger threshold indicates weaker recovery, which can yield different interpretations of resilience outcomes. For example, Khalighi \textit{et al.}~\cite{Khalighi2022Ploscb} showed that memory effects can either enhance or diminish resilience depending on the chosen threshold. Therefore, threshold choice should align with specific real-world contexts and research objectives. Given the theoretical nature of this study, simpler models and definitions were preferred to illustrate fundamental concepts without compromising key insights about resilience, which inherently depends on near-complete recovery to the original stable state.

\subsection*{Context-dependent effects on resistance}
Similar to resilience, the definition of resistance in responding to endogenous perturbations assumes the system initially resides at a positive stable equilibrium. Generally, our results show that memory increases resistance, enabling systems with memory to withstand perturbations more effectively compared to systems without memory. This improved resistance typically occurs when the system state is near the basin's bottom or a stable equilibrium and is subsequently perturbed.

\begin{figure}[ht!]
    \centering
    \includegraphics[width=.7\textwidth]{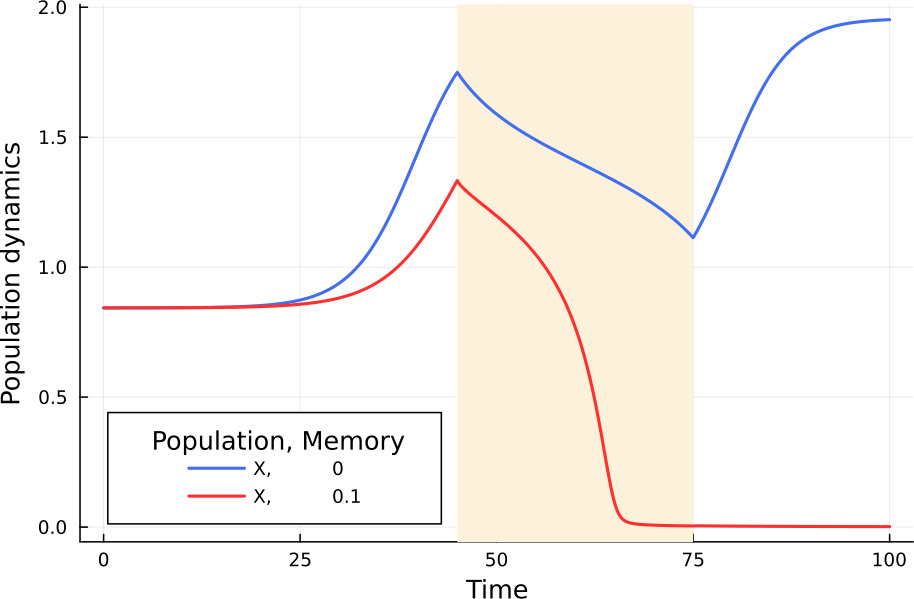}
    \caption{
Memory-induced transition to an alternative stable state. In the presence of memory, slower system dynamics keep the state near the unstable point longer, so even a small perturbation can shift the system to an alternative stable state.
}
    \label{fig: s38}
\end{figure}

However, Figure~\ref{fig: s38} highlights a contrast outcome: if the state is near an unstable point, and in the absence of exogenous perturbations, memory may reduce resistance, potentially leading to transitions to alternative stable states. The dynamics in this scenario occur as follows: when the state is near the unstable point, the system naturally moves toward a stable equilibrium (rolling downhill). However, due to the flattening of the landscape caused by memory, the system initially moves more slowly compared to a memory-free system. Consequently, when a perturbation occurs, the state of the system with memory remains closer to the unstable point, while the state without memory has progressed further towards the basin's bottom.

During the perturbation, the system with memory crosses the unstable point and shifts toward an alternative stable state (e.g., extinction or zero state) due to its proximity to this alternate state. In contrast, the system without memory remains in the original basin when the perturbation ends, allowing it to recover to the initial stable equilibrium.

This scenario demonstrates that memory can reduce resistance under specific conditions by lowering the threshold required to shift states. Therefore, consistent consideration of initial conditions and definitions is crucial for obtaining fair and valid comparisons in studies involving memory effects.

\subsection*{On extending dynamical systems models}
Generalizing from single-variable dynamics to two or more dimensions introduces several challenges. Interactions among state variables mean that the flow is seldom the gradient of a scalar potential; even quasi‑potential reconstructions~\cite{Rodriguez2020EscherStairs, Zhou2012Quasi} can fail, leaving some systems with no meaningful ``landscape''. Nevertheless, landscape–flux decompositions have retrieved early‑warning signals in noisy two‑dimensional shallow‑lake models~\cite{Xu_EWS_2023}, and vector‑field (gradient–curl) decompositions have been used to recover partial landscapes in gene‑regulatory networks~\cite{Zhou2012Quasi} and forest–savanna systems \cite{Xu2021Unifying}.

Despite partial progress, such approaches do not generalize easily, and many systems remain analytically intractable. Accordingly, our framework applies only to models that can be written in the rational form of two polynomials (See Section S1). Many realistic equations, such as consumer–resource systems with non‑rational functional responses, delay differential equations, or piecewise‑defined dynamics, lie outside this class.

\subsection*{Additional stability metrics}
A complementary stability metric is the \textit{mean exit time}~\cite{Arani_exit_time_2021}, defined as the time taken for stochastic fluctuations to drive the system across its basin boundary. Because it combines barrier height, slope, and effective diffusion, it reflects aspects of stability that depth or flatness alone cannot. Arani \emph{et al.},~\cite{Arani_exit_time_2021} showed that exit‑time distributions contract sharply as a system nears a critical transition in memory‑free models. Because fractional memory alters the geometry of the critical points themselves (barrier curvature and height), it is likely to reshape exit-time statistics as well. Quantifying how mean exit time depends on memory strength across our ensemble is thus a promising direction for future work and could help reconcile the apparent trade-off between recovery speed and robustness. Although we did not compute exit times here, we highlight this metric as a valuable bridge between the deterministic geometry analyzed in this study and probabilistic resilience under noise.

\subsection*{Alternative memory types and kernels}
Another important consideration is the type of memory employed in this study. The current work applies fractional derivatives, representing long-term memory effects characterized by a power-law kernel~\cite{kilbas2006theory,podlubny1998fractional}. This approach is particularly valuable as it allows the analytical exploration of stability landscapes through fractional differential equations, providing a framework for understanding how memory influences system dynamics. Additionally, power-law kernels effectively represent real-world scenarios such as anomalous diffusion~\cite{METZLER20001} in spatially structured populations, as discussed in Section S2.

However, other forms of memory, such as exponential kernels or discrete delays, also merit investigation~\cite{Yang2023}. Exploring alternative memory types could provide richer insights into how memory reshapes stability landscapes, particularly by revealing the roles of latency, transient dynamics, and conditions that give rise to distinct forms of memory-dependent stability and resilience.

\section*{S7: Supplementary videos}

\subsection*{Video S1: Anomalous diffusion in spatially structured population agent-based models}
We simulate an agent-based lattice model in Julia with a fixed seed \texttt{MersenneTwister(123)}. Two runs evolve in lock step: a homogeneous \(28\times28\) grid and a porous \(30\times30\) grid with 116 randomly blocked cells. Eight agents start from matched positions within the \(28\times28\) overlap. At each of 300 time steps, agents are randomly ordered, move at most one cell to a free von Neumann neighbor, then attempt reproduction into a free neighbor with probability \(p_{\mathrm{birth}} = r\!\left(1-\frac{N}{K}\right)\) with \(r=0.2\). Local crowding reduces this by \(\exp(-\alpha\, n_{\mathrm{local}})\) with radius \(2\) and \(\alpha=0.1\). In the porous run, agents adjacent to an obstacle become trapped for \(5\) steps before resuming motion. We record positions and \(N(t)\) each step and render a composite animation with two spatial panels and the population time series at \(10\) fps.

\url{https://github.com/moeinkh88/porous-media-simulation/blob/main/videos/comparison_spatial_and_population5.gif}

\subsection*{Video S2: Surpasses the unstable point, then returns to the original basin}
We simulate the fractional herbivory model
\[
\mathcal D^{\alpha}X(t)= r\,X\!\left(1-\frac{X}{K}\right) \;-\; B(t)\,\frac{X}{X+A},
\]
for two derivative orders, \(\alpha=1\) (memory free) and \(\alpha=0.8\) (memory strength \(1-\alpha=0.2\)), with \(r=0.8\), \(K=3\), \(A=0.2\). The state is initialized near the upper stable equilibrium, \(X_0=1.9568\). A single pulse is applied to the control parameter: \(B(t)=0.6\) except on \(t\in[10,20]\) where \(B(t)=0.6+P\) with \(P=0.275\). We integrate on \(t\in[0,100]\) with step \(h=0.01\). The potential \(V(X)\) is reconstructed from the simulated paths by computing \(dX/dt\) with a fourth order central difference and then numerically integrating \(-\!\int (dX/dt)\,dX\). The animation renders a \(2\times 2\) layout: bifurcation diagram, trajectories for \(\alpha=1\) and \(\alpha=0.8\) with the pulse interval shaded, and the evolving \(V(X)\). In this setting, the memory trajectory crosses the unstable point during the pulse and later returns to the original basin.

\url{https://github.com/moeinkh88/porous-media-simulation/blob/main/videos/dynamic_plots.mp4.avi} 

\subsection*{Video S3: Surpasses the unstable point, then shifts to an alternative basin}
The setup matches Video S2 except for a slightly stronger pulse, \(P=0.2755\). All other parameters and numerics are identical: \(r=0.8\), \(K=3\), \(A=0.2\), \(B(t)=0.6\) outside \(t\in[10,20]\), \(X_0=1.9568\), \(\alpha\in\{1,0.8\}\), time window \([0,100]\), and step \(h=0.01\). The same procedure is used to reconstruct \(V(X)\). The animation shows that the memory trajectory again surpasses the unstable point during the pulse. After the pulse ends, it attempts to return toward the unstable equilibrium, but with this slightly stronger disturbance, the system fails to recover and instead transitions into the alternative basin.

\url{https://github.com/moeinkh88/porous-media-simulation/blob/main/videos/dynamic_plots2.mp4.avi}

\subsection*{Video S4: Regime delayed collapse} We simulate the one-dimensional quorum-sensing model with a fractional derivative \(\alpha\in\{1.0,\,0.75\}\):
\[
\mathcal D^{\alpha}A(t)= \frac{V A^2}{K + A^2} + A_0 - d(\rho(t))\,A(t),\qquad 
d(\rho)=0.1+\frac{1-\rho}{\rho(2-\rho)}.
\]
Parameters are \(V=3\), \(K=1\), \(A_0=0.05\), initial state \(A(0)=2.0\), time window \([0,500]\), step \(h=0.01\), and three corrector iterations. White Gaussian noise with standard deviation \(1.3\) is added to the right hand side. The environment schedule is a step change in \(\rho\): \(\rho(t)=0.4\) for \(t<100\) and \(\rho(t)=0.29\) thereafter. We integrate with \texttt{FDEsolver}~\cite{khalighi2024algorithm}, sample the state at integer times, and compute the instantaneous potential
\[
U(A;\rho)= -\Big(V\big(A-\sqrt K\,\arctan(A/\sqrt K)\big)+A_0A-\tfrac{d(\rho)}{2}A^2\Big).
\]
For clarity, the bifurcation diagram and the potential curve \(U(A;\rho(t))\) shown in the background are the same for both \(\alpha=1\) and \(\alpha=0.75\) (they do not depend on \(\alpha\)); with memory the true landscape evolves in time and cannot be represented by a single static curve. We therefore use this instantaneous, memory free surrogate as a common reference and overlay only the moving state marker (“ball”) to illustrate how memory alters trajectories relative to the same geometry. The animation renders, for each \(\alpha\), three panels: equilibria versus \(\rho\) with the current state marked, \(U(A;\rho(t))\) with the state highlighted, and the time series \(A(t)\). Under this schedule, the \(\alpha=0.75\) run collapses later than the \(\alpha=1\) run, illustrating delayed collapse due to memory.

\url{https://github.com/moeinkh88/porous-media-simulation/blob/main/videos/scenario1_collapse.mp4}

\subsection*{Video S5: Regime collapse and delayed recovery} 
The model, numerics, and outputs match Video S4, with a longer window \([0,1200]\) and the piecewise \(\rho\) schedule
\[
\rho(t)=
\begin{cases}
0.335, & t<100,\\
0.289, & 100\le t<300,\\
0.302, & 300\le t<500,\\
0.335, & t\ge 500,
\end{cases}
\]
and noise standard deviation \(1.3\). This sequence drives a collapse and then relaxes conditions back toward baseline. The animation shows that after conditions improve the \(\alpha=0.75\) run returns to the original basin later than the \(\alpha=1\) run, demonstrating delayed recovery driven by memory.

\url{https://github.com/moeinkh88/porous-media-simulation/blob/main/videos/scenario1_delayed_recovery.mp4}

\subsection*{Video S6: Regime collapse and fast recovery} 
The setup is as in Video S4 with time window \([0,400]\), noise standard deviation \(1.0\), and a short pulse in \(\rho\):
\[
\rho(t)=
\begin{cases}
0.4, & t<100,\\
0.25, & 100\le t<120,\\
0.4, & t\ge 120.
\end{cases}
\]
All other parameters are unchanged and integration uses \texttt{FDEsolver} with step \(h=0.01\). The short pulse places the \(\alpha=0.75\) trajectory near the moving ridge; once conditions revert, the barrier reforms behind the state and the system slides back to the original basin sooner than in the \(\alpha=1\) run, illustrating fast recovery with memory.

\url{https://github.com/moeinkh88/porous-media-simulation/blob/main/videos/scenario1_fast_recovery.mp4}

\subsection*{Video S7: Regime collapse and delayed recovery}
The setup is as in Video S4 with initial state \(A(0)=1.5\), time window \([0,500]\), and additive white Gaussian noise of standard deviation \(1.0\) on the right-hand side. The environment schedule is
\[
\rho(t)=
\begin{cases}
0.32, & t<100,\\
0.29, & 100\le t<278,\\
0.32, & t\ge 278.
\end{cases}
\]
In this schedule, the \(\alpha=0.75\) trajectory briefly passes the instantaneous unstable point during the low \(\rho\) interval and, once \(\rho\) returns, rolls back to the original basin as the barrier reforms, demonstrating rollback after a near crossing. By contrast, the memory-free system remains at the lower equilibrium.

\url{https://github.com/moeinkh88/porous-media-simulation/blob/main/videos/scenario2.mp4}

\subsection*{Video S8: Broadened hysteresis under ramps} 
Same model and numerics as Video~S4, where, \(A(0)=1.5\), and additive white noise \(\sigma=1.2\). The control parameter follows a smooth ramp down and up:
\[
\rho(t)=
\begin{cases}
\mathrm{ramp}(t;0,\,T_c,\,0.38,\,0.29), & 0\le t<T_c,\\[2pt]
\mathrm{ramp}(t;T_c,\,1000,\,0.29,\,0.38), & T_c\le t\le 1000,
\end{cases}\qquad T_c=250,
\]
with \(\mathrm{ramp}(t;t_0,t_1,y_0,y_1)=y_0+(y_1-y_0)\,s\), \(s=\mathrm{smoothstep}\!\left(\frac{t-t_0}{t_1-t_0}\right)\), and \(\mathrm{smoothstep}(u)=u^3(10-15u+6u^2)\) for \(u\in[0,1]\). 

\url{https://github.com/moeinkh88/porous-media-simulation/blob/main/videos/scenario3.mp4}


\end{document}